\newif\ifdraft
\documentclass[%
12pt
]{article}
\usepackage{amsmath,amssymb,amscd,verbatim,amsthm}
\ifdraft
  
  \usepackage{showkeys}
\fi
 
\newcommand\set[1]{\left\{\,#1\,\right\}} 
\newcommand\sing[1]{\left\{#1\right\}} 
\newcommand\seq[1]{\left<\,#1\,\right>} 
\newcommand\sseq[1]{\langle#1\rangle} 
\newcommand\krtree{{}^{\kappa>}2}
\newcommand\closure[1]{{#1}^{\wedge}}

\newcommand\is[2]{\text{IS}(#1,#2)}
\newcommand\cof{\operatorname{cf}}
\newcommand\cf{\operatorname{cf}}
\newcommand\restrict{{\upharpoonright}} 
\newcommand{\downto}{\restrict}
\newcommand\otp{\operatorname{otp}}
\newcommand\lex{<_{\text{lex}}}
\newcommand\qle{<_{\text{Q}}}

\newcommand\dom{\operatorname{dom}}
\newcommand\ran{\operatorname{ran}}
\newcommand\clp{\operatorname{clp}}

\newcommand\cat{{}^{\frown}}
\newcommand{\pre}[1]{{^{#1}}}
\newcommand{\prel}[1]{\pre{#1{>}}}
\newcommand{\prele}[1]{\pre{#1{\ge}}}
\newcommand{\bool}[1]{\|#1\|}
\newcommand{\QQ}{\mathbb{Q}}    
\newcommand\qkappa{\mathbb{Q}_\kappa}
\newcommand\compact{small }

\let\meet=\wedge

\DeclareMathOperator{\uf}{\mathcal D}
\newcommand\deq{:=}
\newcommand\cut{|}
\let\meet=\wedge
\newcommand\case[2]{\paragraph{Case #1: #2.}}
\newtheorem{theorem}{Theorem}[section]
\newtheorem{lemma}[theorem]{Lemma}
\newtheorem{cor}[theorem]{Corollary}

\newtheorem{question}[theorem]{Question}

  \newtheorem{claim}{Claim}[theorem]
 
\theoremstyle{remark}

\theoremstyle{definition}
\newtheorem{definition}[theorem]{Definition}

\newcommand\supp{\operatorname{supp}}
\newcommand\rest{\restrict}
\newcommand\forces{\Vdash}

\newcommand\into{\to}
\newcount\skewfactor
\def\mathunderaccent#1#2 {\let\theaccent#1\skewfactor#2
\mathpalette\putaccentunder}
\def\putaccentunder#1#2{\oalign{$#1#2$\crcr\hidewidth
\vbox to.2ex{\hbox{$#1\skew\skewfactor\theaccent{}$}\vss}\hidewidth}}
\def\name{\mathunderaccent\tilde-3 }

\newcommand\cone{\operatorname{Cone}}
\newcommand{\dless}{\lessdot}

\newcommand{\altperm}{P}
\newcommand\ddprime{d^{\,\prime\prime}}
\begin{document}
\title{A partition theorem for a large dense linear order}
\author{M.~D\v{z}amonja%
\thanks{Support by EPSRC through an
Advanced Fellowship is gratefully acknowledged.}\\
\scriptsize{School of Mathematics}\\
\scriptsize{University of East Anglia}\\
\scriptsize{Norwich, NR4 7TJ, UK}\\
J.A.~Larson\thanks{Support by EPSRC and the University of East Anglia
during the period when the project was started, and by the University
of M\"unster, during the writing of the first draft of the paper
is gratefully acknowledged.}
and W.J.~Mitchell$^\dagger$
\thanks{Research on this paper was partly supported by grant number
DMS 0400954 from the US National Science Foundation.}
\\
\scriptsize{Department of Mathematics}\\
\scriptsize{University of Florida-Gainsville}\\
\scriptsize{358 Little Hall, PO Box 118105}\\
\scriptsize{Gainesville, FL 32611--8105, USA}\\
}
\maketitle
\begin{abstract}
Let $\mathbb{Q}_\kappa=(Q,\le_Q)$ be a strongly $\kappa$-dense linear order
of size $\kappa$ for $\kappa$ a suitable cardinal.
We prove, for $2\le m<\omega$, that there is a finite value $t_m^+$
such that the set $[Q]^m$ of $m$-tuples from $Q$ can be divided into
$t_m^+$ many classes, such that whenever any of these classes $C$ is
colored with $<\kappa$ many colors, there is a copy
$\mathbb{Q}^\ast$ of $\mathbb{Q}_\kappa$ such that $[\mathbb{Q}^\ast]^m\cap C$ is
monochromatic. As a consequence we obtain that whenever we color 
$[\mathbb{Q}_\kappa]^m$ with $<\kappa$ many colors, there is a copy of
$\mathbb{Q}_\kappa$ all $m$-tuples from which are colored in at most $t_m^+$
colors. In other words, the partition relation
$\mathbb{Q}_\kappa\rightarrow(\mathbb{Q}_\kappa)^m_{<\kappa/r}$ holds  
for some finite $r=t_m^+$. 

We show that $t_m^+$ is the minimal value with this property.  
We were not able to give a formula for $t_m^+$ 
but we can describe $t_m^+$ as the cardinality of a certain finite set
of types. We
also give an upper and a lower bound on its value
and for $m=2$ we obtain $t_2^+=2$, while for $m>2$ we have $t_m^+> t_m$,
the $m$th tangent number.  
The paper also contains similar positive partition results about 
$\kappa$-Rado graphs.

A consequence of our work and some earlier 
results of Hajnal and Komj\'ath is that
a theorem of Shelah known to follow from a large cardinal assumption
in a generic extension,
does not follow from any large cardinal assumption on its own.
\end{abstract}

\section{Introduction}\label{sec.intro}
For an infinite cardinal $\kappa$,  the set $\krtree$,
ordered by end-extension, $\subseteq$,
is the complete binary tree on $\kappa$ with
root the empty sequence, $\emptyset$. By $s\meet t$ denote the
\emph{meet} of $s$ and $t$, namely the longest initial segment
of both $s$ and $t$.  Call two elements $s$ and $t$ of $\krtree$
\emph{incomparable} if neither is an end-extension of the other.
The lexicographic order for us will be the partial order $\lex$ on $\krtree$ 
defined by $s\lex t$ if
$s$ and $t$ are incomparable and $(s\meet t)\cat\langle 0\rangle\subseteq s$
and $(s\meet t)\cat\langle 1\rangle\subseteq t$.
We also define a linear order $\le_Q$ on $\krtree$ by letting
$s\qle t$ if and only if one of the following conditions holds:
(1) $s=t$;
(2) $t\cat\langle 0\rangle\subseteq s$;
(3) $s\cat\langle 1\rangle\subseteq t$; or
(4) $s$ and $t$ are incomparable and $s\lex t$.

Let us recall the definition of a (strongly) $\kappa$-dense linear order:
it is a linear order $\le$ on a set $L$ in which for every two subsets
$A,B$ of $L$ of size $<\kappa$
satisfying the property that for all $a\in A$ and $b\in B$ the relation 
$a<b$ holds, there is $c\in L$ such that $a<c<b$ for all $a\in A$ and 
$b\in B$. Since $A$ or $B$ may be taken empty here, the definition in 
particular implies that there are no endpoints in the order.
The adjective `strongly' is used to distinguish this type of ordering from
the strictly weaker notion in which for any $a<b$ in $L$ one is required to
have $\kappa$ many $c$ with $a<c<b$. Such orders were first studied by 
Felix Hausdorff in 1908 (\cite{Hausdorff}) and have been of continuous 
interest since. A recent paper on the subject, where one can also find 
a number of further references, is \cite{DzT} by M. D\v zamonja
and Katherine Thompson,
where they give a classification of 
$\kappa$-dense linear orders which are also $\kappa$-scattered. 
In this paper we shall only deal with strongly $\kappa$-dense 
linear orders so we shall omit the adjective `strongly' from our notation.

It is easy to check that for regular $\kappa$
the structure $\mathbb{Q}_\kappa:=(\krtree,\qle)$ is a 
$\kappa$-dense linear order (see \ref{lem.density}). In case that $\kappa^{<\kappa}=\kappa$
then of course this order has cardinality $\kappa$. All $\kappa$ 
that we shall work with will satisfy this additional cardinal 
arithmetic assumption. It is well known and easily proved using 
a back-and-forth argument that in this case the $\kappa$-dense linear 
order of size $\kappa$ is unique up to isomorphism,
and that it is a $\kappa$-saturated homogeneous model 
of the theory of a dense linear order with no endpoints. 
(In fact these properties are equivalent to 
$\kappa=\kappa^{<\kappa}$, as follows from Saharon Shelah's 
classification theory, see \cite{Shc}).

For $\kappa=\omega$, $\mathbb{Q}_\omega=\mathbb{Q}$ 
is a countable dense linear order with no endpoints,
so it has the order type $\eta$ of the rationals. 
An unpublished result of Fred
Galvin as quoted in \cite{HaKo} is that 
$\mathbb{Q}\rightarrow[\mathbb{Q}]^2_{<\omega,2}$, see the
notation below. 
Denis~C. Devlin \cite{DDthesis} proved that 
$$\mathbb{Q}\rightarrow(\mathbb{Q})^n_{<\omega,t_n}\text{ and }
\mathbb{Q}\nrightarrow(\mathbb{Q})^n_{<\omega,t_n-1}$$
where the value of $t_n$ is the $n$-th tangent number. 
The notation here means that
for every $n$, when one colors $[\mathbb{Q}]^n$ with finitely many 
colors
(this is the role of the part $<\omega$ in the subscript),
there is a copy $\mathbb{Q}^\ast$ of $\mathbb{Q}$ with the property that
$[\mathbb{Q}^\ast]^n$ is colored in at most $t_n$ colors. 
At the same time, $t_n$ is the smallest
number for which such a statement holds.

Tangent numbers may be computed using the power series
$\tan(x)=\sum_{1}^{\infty} t_n\frac{x^{2n-1}}{(2n-1)!}$.
Devlin's proof used the language of category theory. 
A proof of this theorem using trees was sketched in 
in the Farah-Todorcevic book \cite{FaTo}; a complete proof 
was given by Vojkan Vuksanovic \cite{V02pams} in which
he uses the special case for the complete binary tree ${}^{\omega>}2$ of
Keith R. Milliken's theorem (\cite{KM79}) about weakly embedded subtrees.
Since many of the notions
used by Vuksanovic generalize to arbitrary infinite $\kappa$ in place
of $\omega$, one may wonder if his proof may be used to obtain a partition
theorem  for $\kappa$-dense linear orders. Here we use some of these
ideas along with new insights to get such a theorem. 
Our work was also inspired by the strong diagonalization of Norbert Sauer in 
\cite{Spreprint}, the approach to similarities in
a triple paper by Claude Laflamme, Sauer and Vuksanovic \cite{LSVpreprint},
and the use of collapses both in the Shelah \cite{Sh288} version of
the Halpern-L\"{a}uchli Theorem and in 
another paper of Vuksanovic \cite{Vrandom}.

One difficulty of the generalization was that the special case of Milliken's
theorem on binary trees was only known to be valid for 
${}^{\omega>}2$, not for an
arbitrary $\kappa$.  In particular, suppose $\kappa$ is uncountable
and $\prec$ is a well-ordering of $\krtree$ with the
property that whenever $s$ is shorter than
$t$ then also $s\prec t$.  Define a coloring of
the height $2$ complete binary trees strongly embedded in
$\krtree$ by $g(S)=0$ if and only if the lexicographic
order and the $\prec$-order agree on the leaves of $S$.
For any strongly embedded copy of $\krtree$, this coloring
on the binary trees generated by pairs of nodes on the
$\omega$th level is essentially the Sierpinski partition,
so has edges of both colors.

Milliken's theorem for weakly embedded subtrees follows
from his theorem for strongly embedded subtrees, and to prove it,
he uses a generalization of the Halpern-L\"{a}uchli Theorem
due independently to Richard Laver and David Pincus (see \cite{KM79}).
Here we are able to generalize part of D. Devlin's theorem starting
with a theorem of Shelah from  \cite{Sh288}, which is a generalization
of the Halpern-L\"{a}uchli theorem \cite{HaLa}.
We  slightly improve Shelah's Theorem to colorings of antichains rather than
only level sets.  
Let us now state our main result for large dense linear orders.

\begin{theorem}\label{thm.partition} 
For every natural number $m$ there is a value $t_m^+<\omega$ such that
for any cardinal $\kappa$ which is measurable
in the generic
extension obtained by adding $\lambda$ Cohen subsets of $\kappa$,
where $\lambda$ is some cardinal satisying $\lambda\to(\kappa)^{2m}_{2^{\kappa}}$, 
the $\kappa$-dense linear order $\mathbb{Q}_\kappa$ 
satisfies 
$$\mathbb{Q}_\kappa\rightarrow(\mathbb{Q}_\kappa)^m_{<\kappa,\, t_m^+}\text{ 
and }
\mathbb{Q}_\kappa\nrightarrow(\mathbb{Q}_\kappa)^m_{<\kappa,\, t_m^+-1}.$$
\end{theorem}

In Theorem \ref{thm.upper.bound.1}, the positive partition relation
is shown to hold for $t_m^+$ the number of \emph{sparse vip $m$-types}
(defined in Section \ref{sec.upper.bnd.1}).  The sparse vip $m$-types
are closely related to those unique $m$-element \emph{strongly diagonal} 
subsets of ${}^{2m-2\ge}2$ that are representatives of the ``essential types''
in \cite{V02pams}.  If we close such a strongly diagonal set
under initial segments and add a \emph{vip level order}, we obtain
a sparse vip $m$-type, and all sparse vip $m$-types are obtained
in this way.
 
In Theorem \ref{thm.lower.bound.1} we show that the negative
partition relation holds for the same value of $t_m^+$.
In Theorem \ref{thm.canon.1} we show that for $\kappa$
as in Theorem \ref{thm.partition} there is a \emph{canonical
partition} $\mathcal{C}=\set{C_0,C_1,\dots,C_{t_m^+-1}}$
of $[\mathbb{Q}_\kappa]^m$.  That is, a partition
whose classes are \emph{persistent} and \emph{indivisible}.
We say $C_j$ is \emph{persistent} if for every $\kappa$-dense $Q^*\subseteq
\mathbb{Q}$ the set $[Q^*]^m\cap C_j$ is non-empty.  We say
$C_j$ is \emph{indivisible} if for every coloring of $C_j$ with
fewer than $\kappa$ many colors, there is a $\kappa$-dense
subset $Q^*\subseteq\mathbb{Q}$ on which $[Q^*]^m\cap C_j$ is
monochromatic.

Richard Rado \cite{Rado64} constructed a (strongly) universal
countable graph in 1964.  That is, he constructed a countable
graph for which every countable graph is an induced subgraph.
By a \emph{$\kappa$-Rado graph} we mean a graph
$G$ of size $\kappa$ with the property that for every two disjoint subsets
$A$, $B$ of $G$, each of size $<\kappa$, there is $c\in G$ connected to
all points of $A$ and no point of $B$. The existence of such $G$ follows
from the assumption $\kappa^{<\kappa}=\kappa$.

Note that a $\kappa$-Rado graph embeds every graph 
with at most $\kappa$-many vertices. 
That is, it is universal for the family of graphs of size
at most $\kappa$.  For $\kappa=\omega$, this graph is also 
called the infinite random graph.

Interest in Rado graphs and the uncountable continues.
Let $G_\omega$ denote the $\omega$-Rado graph.
Recently Gregory Cherlin and Simon Thomas \cite{CT02} 
have shown that for any infinite cardinals $\kappa\le\lambda$,
the assumption $\lambda\le 2^\kappa$ is equivalent to
the existence of a graph $G^*$ of size lambda which is elementarily
equivalent to $G_\omega$ and which has a vertex whose set of 
neighbours has size $\kappa$.

Here is our main result on $\kappa$-Rado graphs.

\begin{theorem}\label{thm.rpartition} 
For every natural number $m$ there is a value $r_m^+<\omega$ such that
for any cardinal $\kappa$ which is measurable
in the generic
extension obtained by adding $\lambda$ Cohen subsets of $\kappa$,
where $\lambda$ is some cardinal satisying 
$\lambda\to(\kappa)^{2m}_{2^{\kappa}}$, 
the $\kappa$-Rado graph $\mathbb{G}_\kappa$
$$\mathbb{G}_\kappa\rightarrow(\mathbb{G}_\kappa)^m_{<\kappa,\, r_m^+}\text{ 
and }
\mathbb{G}_\kappa\nrightarrow(\mathbb{G}_\kappa)^m_{<\kappa,\, r_m^+-1}.$$
\end{theorem}

In Theorem \ref{thm.upper.bound.2}, the positive partition relation
is shown to hold for $r_m^+$ the number of \emph{vip $m$-types}
(defined in Section \ref{sec.upper.bnd.1}).  The vip $m$-types
are closely related to those unique $m$-element \emph{strongly diagonal} 
subsets of ${}^{2m-2\ge}2$ that are representatives of the ``essential types''
used by Laflamme, Sauer and Vuksanovic in
\cite{LSVpreprint} and by Vuksanovic in \cite{Vrandom} 
for the countable Rado graph.  
If we close such a strongly diagonal set
under initial segments and add a \emph{vip level order}, we obtain
a vip $m$-type, and all vip $m$-types are obtained in this way.

In Theorem \ref{thm.lower.bound.2} we show that the negative
partition relation holds for the same value of $r_m^+$.

It is not known if the large cardinal assumptions used in 
the proof of Shelah's Theorem from \cite{Sh288}
are optimal; in Section 
\ref{sec.remarks} we comment more on this as well as on the consistency
strength of these requirements. We note that in conjunction with
a result of Andr\'as Hajnal and P\'eter Komj\'ath from \cite{HaKo} our theorem
gives some necessary indestructibility conditions on $\kappa$ from Shelah's
Theorem. The paper also includes a section with
a proof of the particular variant of Shelah's Theorem that
we need.

For the remainder of the paper an unattributed $m$ will mean a natural
number with $2\le m$ and $\kappa$ a cardinal
satisfying the hypotheses of Theorem \ref{thm.partition} for some
number $m$. In particular,
$\kappa=\kappa^{<\kappa}$ and $\kappa$ is a strong limit.

The paper is organized as follows.  In Section \ref{sec.uni}, we 
define the notion \emph{$\prec$-similarity},
and state the special variant of
Shelah's Theorem
that we will use. 

In Section \ref{sec.upper.bnd.1}, we define the notions of
\emph{diagonal}, vip order and sparse vip $m$-type
and use them together with our variant of Shelah's Theorem to
prove Theorem \ref{thm.upper.bound.1}.

In Section \ref{sec.almost.perfect}, we show that for
any $S\subseteq\krtree$ with $(S,\qle)$, one can build
an \emph{almost perfect} $\kappa$-dense subtree inside
the tree of nodes of $\krtree$
with a $\kappa$-dense set of extensions in $S$.

In Section \ref{sec.lower.bnd.1}, we use the construction
techniques of Section \ref{sec.almost.perfect} to
prove Theorem \ref{thm.lower.bound.1} by showing all
sparse vip $m$-types are embeddable in every sparse diagonal
$\kappa$-dense subset.

In Section \ref{sec.sparse.vip.types}, 
we show the critical numbers $t_m^+$ for the $\kappa$-dense
linear order are bounded below by $t_m$ and bounded above 
by $t_m(m-1)!\left[\prod_{i<m-1}(i!)^2\right]$,
and indicate how one may compute the value of $t_m^+$
recursively.  We conclude the section with a table of
small values of $t_m^+$ quoted from an upcoming paper by Jean Larson 
\cite{Jeanfincomb}.

In Section \ref{sec.Shproof}
we prove the variant of Shelah's Theorem that we use, formulating it
using colorings of antichains. 
In Section \ref{sec.remarks} we comment on the necessity
of the use of large cardinals in our theorem and 
in Shelah's Theorem and the way that the results of this paper
shed light on that question. We also give some open questions.

In Section \ref{sec.upper.bnd.2}, we use  nuanced
diagonalization to prove Theorem \ref{thm.upper.bound.2}
giving the upper bound on the critical values $r^+_m$
for $\kappa$-Rado graphs

In Section \ref{sec.lower.bnd.2}, we prove a reduction
theorem for  the translations to tree form of increasing embeddings 
of the $\kappa$-Rado graph $\mathbb{G}_\kappa$: 
for every such translation with range $D$
there is always a particularly nice form of a diagonalization 
which has range a subset of $D$ and is itself the translation
of an increasing embedding of $\mathbb{G}_\kappa$ into itself.
We use this reduction to prove Theorem \ref{thm.lower.bound.2}.
by showing every vip $m$-type is embeddable in the range
of these particularly nice diagonalizations.

The remainder of this introduction is devoted to background
information, including some definitions, notation, and
the statements of some theorems that will be
used as tools.  
For any cardinal $\lambda$, let $[A]^\lambda$ denote the
collection of all subsets of $A$ of cardinality $\lambda$,
and let $[A]^{<\lambda}$ denote the collection of all 
subsets of $A$ of cardinality less than $\lambda$.

For any tree $T=(T,\subseteq)$, 
a node $s$ is a \emph{leaf of $T$} or \emph{terminal node of $T$}
if for all $t\in T\setminus\{s\}$,
one has $s\not\subseteq t$.  
The notion of the meet of two nodes has already been defined, in the 
first paragraph of the Introduction.  If $s$ is a node of $T$ and
both $s\cat\langle 0\rangle$ and $s\cat\langle 1\rangle$ have
extensions in $T$, then we call $s$ a \emph{splitting node of $T$}.
For any subset $S\subseteq T$,
let $\closure{S}$ denote the \emph{meet closure} of $S$, i.e. the set
$\{u:\,u=s\wedge t\mbox{ for some }s,t\in S\}$. Note that $\closure{S}
\supseteq S$ and that $\closure{S}$ is closed under meets and has a unique
node of the smallest length.

\begin{definition}\label{def.IS}
For any tree $T$ of sequences ordered by end extension,
and any node $s\in T$, define
the \emph{set of immediate successors of $s$ in $T$} as
$$\is{s}{T}:=\set{t\in T: s\subseteq t\land (\forall u)(s\subset
  u\subseteq t\implies u=t)}.$$
\end{definition}

\begin{definition}\label{def.level}
For any tree $T$ of sequences under end extension,
the \emph{$\alpha$th level of $T$}, in symbols $T(\alpha)$, 
is the set of all nodes $t\in T$ for which $\alpha$ is the order
type of the set of predecessors of $t$, namely
$\set{s\in T: s\subsetneq t}$. For $s\in T$ the length $\lg (s)$ is
defined to be $\alpha$ if and only if $s\in T(\alpha)$. Call $T$ an 
\emph{$\alpha$-tree} if each branch of $T$ has order type
$\alpha$.
\end{definition}

The above notion of an $\alpha$-tree in the case that $\alpha$ is a cardinal
differs from the usual notion of an $\alpha$-tree (see e.g. \cite{Kunen}),
which is defined as a tree of height $\alpha$ all of whose levels have size
$<\alpha$. We trust that no confusion will arise, but emphasize that our definition
is given above.
Note that for the complete binary tree $T=\krtree$, the
levels of $T$ are $T(\alpha)={}^\alpha 2$ for $\alpha<\kappa$.
Call $A\subseteq T$ a \emph{level set} if $A\subseteq T(\alpha)$
for some $\alpha$.

Milliken (see Definition 1.2 of \cite{KM79})
defined a notion of strongly embedded tree which
can be simplified in the case of a binary tree. 
The key idea is that there are \emph{splitting levels}
(see the range of $h$ below); 
all nodes of a strongly embedded tree split at splitting
levels as much as is possible, and no splitting occurs
elsewhere. Another way to say this is that
a strongly embedded subtree of height
$\alpha$ is a copy of ${}^{\alpha>}2$ in $\krtree$
where levels are mapped to levels.
Note that a strongly embedded subtree $S$ of $\krtree$
is {\em not} required to be an induced subtree of $\krtree$, in the sense that it is not
required that $S$ is closed under $\subseteq$.

\begin{definition}\label{def.strong.tree}
A subset $S\subseteq\krtree$ is \emph{a strongly embedded subtree of
$\krtree$} if 
\begin{description}
\item{(1)} $(S,\subseteq)$ is an $\alpha$-tree
for some $\alpha\le\kappa$; 
\item{(2)} $S$ has a root and every non-maximal node $s\in S$ has exactly one
extension of
$s\cat\langle 0\rangle$ and one extension of $s\cat\langle 1\rangle$ in
$\is{s}{S}$; 
\item{(3)} there is a \emph{level assignment function}
$h:\alpha\to\kappa$ which is a strictly increasing function 
such that for all $\beta<\alpha$,
$S(\beta)\subseteq {}^{h(\beta)} 2$;
\item{(4)} for any limit $\beta<\alpha$, for any $s\neq t\in S(\beta)$, there is
$\delta<\sup\{h(\gamma):\,\gamma<\beta\}$ such that $s\rest\delta\neq t\rest\delta$.
\end{description}
\end{definition}

Note of course that part (4) above is irrelevant for strongly embedded trees of
height $\le\omega$, as originally studied by Milliken.
As mentioned above, strongly embedded subtrees of $\krtree$
are copies of full binary trees 
of height $\le\kappa$ where levels are mapped to levels. The precise
statement of this is the content of Lemma \ref{lem.strong.embed}.

\begin{definition}\label{def.strong.embed}
For any subsets $S_0$ and $S_1$ of $\krtree$, a function 
$e:S_0\into S_1$ is a \emph{strong embedding}
if it is an injection with the following preservation properties:
\begin{enumerate}
\item (extension)
$s\subseteq t$ if and only if $e(s)\subseteq e(t)$;
\item (level order)
$\lg (s)<\lg (t)$ if and only if $\lg(e(s))<\lg(e(t))$ and 
$\lg (s)=\lg (t)$ if and only if $\lg(e(s))=\lg(e(t))$;
\item (passing number) 
if $\lg (s)< \lg (t)$, then $e(t)(lg(e(s)))=t(\lg (s))$.
\end{enumerate}
\end{definition}

%
%

\begin{lemma}\label{lem.strong.embed}
An $\alpha$-tree $S\subseteq\krtree$ is strongly embedded in $\krtree$
if and only if there is a strong embedding
$e:{}^{\alpha >} 2\to\krtree$ whose range is $S$.
\end{lemma}

\begin{proof}
If $e:{}^{\alpha>} 2\to\krtree$ is a strong embedding, then
there is a level assignment function $h:\kappa\to\kappa$ witnessing
the strong embedding in that
for all $\beta<\alpha$, if $s\in {}^\beta 2$,
then $e(s)\in {}^{h(\beta)}2$ and the immediate
successors of $e(s)$ in the range of $e$ are $e(s\cat\langle 0\rangle)$
and $e(s\cat\langle 1\rangle)$, which are both 
in ${}^{h(\beta+1)}2$.  Thus if $e$ is a strong
embedding its image is a strongly embedded tree.

For the other direction, suppose $S\subseteq \krtree$ 
is a strongly embedded $\alpha$-tree
and $h:\alpha\to\kappa$ is the level assignment function
that witnesses it.  We shall define $e$ on ${}^{\alpha>} 2$ such that
$e$ maps ${}^\beta 2$ into ${}^{h(\beta)} 2$ for all
$\beta<\alpha$, as follows. Let $e(\langle\rangle)$
be the root of $S$. Given $t\in {}^{\alpha>} 2$ non-maximal and of height
$\beta$, let $e(t\cat\langle l\rangle)$
for $l<2$ be the unique extension of $e(t)\cat\langle l\rangle$ in $\is{e(t)}{S}$.
For $\beta<\alpha$ limit and $t\in {}^\beta 2$, note that $\bigcup_{s\subseteq t}
e(s)$ is an element of ${}^{\kappa>}2$ and that it must have an extension in 
${}^{h(\beta)}2$, as $S$ is an $\alpha$-tree. By requirement (4) in 
\ref{def.strong.tree}, this extension is unique and we choose it as $e(t)$.
Now $e$ is a strong embedding whose image is $S$.

\end{proof}

Let us finish the section by proving the above mentioned
fact that for regular $\kappa$
the order $<_Q$ on ${}^{\kappa>}2$ is $\kappa$-dense.

\begin{lemma}\label{lem.density} The order $({}^{\kappa>}2,<_Q )$ is 
$\kappa$-dense if and only if $\kappa$ is a regular cardinal.
\end{lemma}

\begin{proof} Suppose that $\kappa$ is regular and $A,B\subseteq
{}^{\kappa>}2$ are both of size $<\kappa$ and 
such that for all $a\in A$ and $b\in B$ we have
$a<_Q b$ (we write this as $A<_Q B$).
Let $B^\ast:=\{t\in {}^{\kappa>}2:\,(\exists b\in B) t\cat\langle 0\rangle
\subseteq b\}$. If $a\in A$ and $t\in B^\ast$ then there is $b\in B$ such that 
$t\cat\langle 0\rangle
\subseteq b$, in particular $b<_Q t$. By the transitivity of $<_Q$ we obtain $a<_Q t$,
and hence $A<_Q B^\ast$.

Let $\gamma$ be the minimal ordinal such that all elements of $A\cup B$ have length
$<\gamma$, so by the regularity of $\kappa$ we have $\gamma<\kappa$. We define $c\in
{}^{\kappa>}2$ by defining $c(\alpha)$ for $\alpha<\gamma$ by recursion on $\alpha$.
If $c\rest\alpha\in B\cup B^\ast$ we set $c(\alpha)=0$, and we set $c(\alpha)=1$ otherwise.
We claim that $A<_Q c<_Q B$.

If $b\in B$ then the length of $b$ is less than the length of $c$ and hence either
$b\subseteq c$ or $b$ and $c$ are incomparable. In the first case we have defined
$c$ so that $b\cat\langle 0\rangle\subseteq c$. In the second case, if $\alpha$ is the
length of $b\cap c$ and $b(\alpha)=0$, then $t:=b\rest\alpha\in B^\ast$ so we have
defined $c(\alpha)=0$, contradicting the choice of $\alpha$. So $b(\alpha)=1$ and
$c(\alpha)=0$. In both cases we have $c<_Q b$, and hence $c<_Q B$.

If $a\in A$ and $a\subseteq c$ then we claim $a\cat\langle 1\rangle
\subseteq c$. Otherwise,
since the length of $a$ is strictly less than that of $c$ we have $a\cat\langle 0\rangle
\subseteq c$, so $a=c\rest\alpha$ must be a member of $B\cup B^\ast$, a contradiction.
If $a$ and $c$ are incomparable we can similarly show that $a<_{\rm lex} c$. In any
case $a<_Q c$ and we have proved that $A<_Q c$. This proves that $\qkappa$ is
$\kappa$-dense.

Suppose now that $\kappa$ is singular and $\kappa_i$ for $i<\cf(\kappa)$ is an increasing
sequence of limit ordinals with limit $\kappa$, and with $\kappa_0=\kappa$. For all $i$
let $a_i$ be the sequence in ${}^{\kappa_i}2$ which is constantly equal to 1, and let
$A:=\{a_i:\,i<\cf(\kappa)\}$. We then have $|A|<\kappa$, yet we claim that there is no
$c\in {}^{\kappa}>2$ with $A<_Q c$. Namely, let $c\in {}^{\alpha}2$ for some $\alpha<\kappa$
and let $i$ be such that $\alpha\in [\kappa_i, \kappa_{i+1})$.
Then $c\cat\langle 1\rangle\subseteq a_{i+2}$, so $A<_Q c$ cannot hold. This proves that
$\qkappa$ is not $\kappa$-dense.

\end{proof}

An argument similar to the first part of the above proof is presented in the
proof of Lemma \ref{lem.cofinal.dense}.

\section{Uniformization}\label{sec.uni}
In this section we state a particular variant of a theorem 
of Shelah. His theorem was stated as 
a generalization of the Laver-Pincus version of the
Halpern-L\"{a}uchli Theorem, and the variant we state is
a generalization of Milliken's Ramsey theorem for finite 
weakly embedded subtrees.  

Before we can state Shelah's Theorem, we need the definition
of two subsets being $\prec$-similar, where $\prec$ is a level
ordering of $\krtree$: 

\begin{definition}\label{def.level.order}
Say that $\prec$ is a \emph{level ordering} of $\krtree$ 
or alternatively that $\prec$ is \emph{an ordering of the levels}
of $\krtree$, if $\prec$ extends the length order, i.e.
$\prec$ is a linear order of $\krtree$ and $\lg(s)<\lg(t)$
implies $s\prec t$.
\end{definition}

Milliken's Theorem for weakly embedded subtrees 
requires us to recognize subtrees
of the same embedding type.  We will be counting the number
of types of particularly nice trees 
at a later point in the paper, so it is convenient
to relate these types to specific finite examples called \emph{similarity
trees}.

\begin{definition}\label{def.sim.tree}
A \emph{similarity tree} is a finite subtree of ${}^{\omega>}2$
closed under initial segments and such that every
level contains at least one leaf (terminal) node or meet of leaf nodes
(split point). An
\emph{ordered similarity tree} is a similarity tree $t$ with an
ordering $\prec_{t}$ of its levels.


For a finite subset $x$ of $\krtree$ we define as $\clp(x)$ the
subtree $y$ of ${}^{\omega>}2$ that includes the root and is of minimal possible
height such that there is a strong embedding from $\closure{x}$ onto the closure
$\closure{z}$ of the set $z$ of terminal nodes of $y$. If $\prec$ is a given 
order of $\krtree$ then $\prec_x$ is the order on $\clp(x)$ induced by 
the strong embedding from $\clp(x)$ to $x$ and $\prec$.
\end{definition}


To illustrate the definition, consider a specific antichain:
\[
x=\set{\langle 0,0,0,1\rangle, \langle 0,0,1\rangle,
\langle 0,1,0,0,0,1\rangle}.
\]
Further suppose that
$\prec$ is the following ordering of the levels:
$\langle 0,0\rangle\prec \langle 0,1\rangle$; 
$\langle 0,0,0\rangle\prec \langle 0,1,0\rangle
\prec \langle 0,0,1\rangle$; and
$\langle 0,0,0,1\rangle\prec \langle 0,1,0,0\rangle$.
The corresponding similarity tree is $(\clp(x),\subseteq)$
where 
\[
\clp(x)=\set{\emptyset,\langle 0\rangle,
\langle 1\rangle, \langle 0,0\rangle, \langle 0,1\rangle,
\langle 1,0\rangle, \langle 0,0,1\rangle, \langle 1,0,0\rangle,
\langle 1,0,0,0\rangle}.  
\]
The induced ordering of the levels is 
$\langle 0\rangle\prec_x \langle 1\rangle$; 
$\langle 0,0\rangle\prec_x \langle 1,0\rangle
\prec_x \langle 0,1\rangle$; and
$\langle 0,0,1\rangle\prec_x \langle 1,0,0\rangle$.

\begin{definition}\label{def.ordered.sim}
Suppose that $\prec$ is an ordering of the levels of $\krtree$.
Two antichains $x$ and $y$ in $\krtree$ are \emph{similar}
if $\clp(x)=\clp(y)$, and \emph{$\prec$-similar} if the ordering
$\prec$ induces the same ordering $\prec_x=\prec_y$ on the collapsed trees.
In this case we call $(\clp(x),\prec_x)$
the \emph{ordered similarity type} of $(x,\prec)$ and $(y,\prec)$.
\end{definition}

We remark that if $x$ and $y$ are $\prec$-similar antichains,
then $\clp(x)$ and $\clp(y)$ are subtrees of $\krtree$ having
the same \emph{weak embedding type}, a notion which appears in 
various places in the literature, see e.g. \cite{KM79} and 
\cite{V02pams}.  The notion of $\prec$-similarity given above is a translation of
the definition of similarity by Shelah in Definition 2.2 \cite{Sh288}.
It also corresponds to the notion of \emph{strong similarity} from
\cite{LSVpreprint}, which coincides with the notion of
\emph{similarity} of that paper when 
restricted to \emph{strongly diagonal sets}.

Theorem~\ref{thm.SHL} below is a variant of Shelah's 
Lemma 4.1 of \cite{Sh288} together
with Conclusion 4.2 where it is called ``$Pr^f_{ht}(\kappa,m,\sigma)$
(even with (3))".  In fact, 
except for the assertion that $e$ preserves $\prec$ and that we
are dealing with coloring of antichains rather than level sets
(a level
set is a subset of ${}^{\alpha}2$ for some $\alpha<\kappa$), our
Theorem~\ref{thm.SHL} is equivalent to Shelah's Conclusion~4.2 from
that paper with option~(3)(b). The fact that $e$ preserves
$\prec$ implies that option~(3)(a) is automatically satisfied as
well. Lemma~\ref{thm:endhom}, which is the main lemma in the proof of
Theorem~\ref{thm.SHL}, is equivalent to Shelah's Lemma~4.1 with
option~(3)(a). In Shelah's notation
the superscript $f$ of Lemma 4.1 indicates that one gets
a strongly embedded tree $T$, and the subscript $ht$ means that
the coloring restricted to $m$-element level sets of $T$
is \emph{homogeneous}, i.e. is constant on subsets whose collapses
under the level assignment function of $T$ are $\prec$-similar. 
The star orderings, $<_\alpha^*$, in this option, are defined
from the $0$ orderings, $<_\beta^0$.  
We use $\prec$ in place of $<_\beta^0$, and will define
the order of interest at a later point.  We write
$e(\prec)$ for $\set{(e(a),e(b)):a\prec b}$.
Thus given a strong embedding
$e:\krtree\to\krtree$, 
$s<_\alpha^* t$ is translated by
$e(s)e(\prec) e(t)$, which holds 
if and only if $s\prec t$.  We also note that the statement of 
the theorem as we give it only refers to a dense
set of elements $w$. A proof of 
Theorem~\ref{thm.SHL} will be given in section \ref{sec.Shproof}.


The following notation is introduced for convenience.

\begin{definition}
Write $\cone(w)$ for the set of all extensions $t\supseteq w$
in $\krtree$.
\end{definition}

\begin{theorem}\label{thm.SHL}[Shelah \cite{Sh288}]
  Suppose that $m<\omega$
  and $\kappa$ is a cardinal which is measurable
  in the generic
  extension obtained by adding $\lambda$ Cohen subsets of $\kappa$,
  where $\lambda\to(\kappa)^{2m}_{2^{\kappa}}$. 
  Then for any coloring $d$ of the $m$-element antichains of $\krtree$
  into $\sigma<\kappa$ colors, 
  and any well-ordering $\prec$ of the levels of $\krtree$, 
  there is a strong embedding $e\colon\krtree\to \krtree$ 
  and a dense set of elements $w$ such that 
\begin{enumerate}
\item $e(s)\prec e(t)$ for all $s\prec t$ from $\cone(w)$, 
   and 
\item $d(e[a])=d(e[b])$ for all $\prec$-similar $m$-element antichains $a$ and
  $b$ of $\cone(w)$. 
\end{enumerate}
\end{theorem}


\section{Upper bound for dense linear orders}\label{sec.upper.bnd.1}
In this section we prove a limitation of colors
result for $\kappa$-dense linear orders using
Shelah's Theorem \ref{thm.SHL}.

We turn to some ideas from work of 
Sauer, Laflamme and Vuksanovic and others 
for ways to guarantee that the $m$-element 
sets of our yet to be chosen transverse set
have a minimal number of (unordered) weak embedding types.
The notions of \emph{diagonal} set and \emph{strongly diagonal} set
are used in \cite{LSVpreprint}, for
example.  Our definitions below are simplifications of those definitions
to the special case of binary trees.

\begin{definition}\label{def.diagonal}
Call $A\subseteq \krtree$  \emph{diagonal} if it is
an \emph{antichain} (its elements  are pairwise incomparable), 
and its meet closure, $\closure{A}$, is transverse.
Call it \emph{strongly diagonal} if, in addition, 
for all all $t\in A$  and all $s\in \closure{A}\setminus\{t\}$,
the following implication holds:
\[
(\lg(s)<\lg(t)\text{ and }t(\lg(s))=1)\implies
s\subseteq t\text{ or $s$ has no extension in $A$}.
\]
\end{definition}

\begin{lemma}\label{lem.diag}
For any finite diagonal set $a\subseteq \krtree$,
its meet closure has  $|\closure{a}|=2|a|-1$ elements
and $\clp(a)\subseteq {}^{2\left|a\right|-2\ge} 2$.
Thus for positive $m<\omega$, the number of ordered
similarity types of $m$-element diagonal sets is finite.
\end{lemma}

In the case of $\kappa$-dense linear orders, we are particularly
interested in \emph{sparse diagonal} sets. 

\begin{definition}\label{def.sparse}
Call $A\subseteq \krtree$  \emph{sparse} 
if
for all $t\in A$ and $s\in\closure{A}\setminus\set{t}$ 
the following implication holds:
\[
(\lg(s)<\lg(t)\text{ and }t(\lg(s))=1)\implies
s\subseteq t.
\]
\end{definition}

Notice that if $A$ is a sparse diagonal set, then it
is a strongly diagonal set.  However, a strongly
diagonal subset need not be sparse. 

We are interested in a special collection of level orders.  We call
them \emph{$D$-vip orders} since the elements of $D$ are
Very Important Points, with special roles to play in the orders,
and these roles continue to be played even when finite subsets are
collapsed. 

\begin{definition}\label{def.vip}
Suppose $T$ is a subtree of $\krtree$ and $D\subseteq\krtree$.
Call $\prec$ a \emph{pre-$D$-vip order} on $T\subseteq\krtree$
if $D$ is transverse and $\prec$ is a well-ordering of the
levels of $T$ such that for every $d\in D$, $d$ is
the $\prec$-least element of its level, $T(\lg(d))$, and for
all $u,v\in T(\lg(d))\setminus\{d\}$,
$\prec$ satisfies the following
condition:
\begin{enumerate}
\item if $d\meet u\subsetneq d\meet v\ne d$, then $u\prec v$, and
\end{enumerate}
If $D$ is diagonal, call $\prec$ a 
\emph{$D$-vip order} if it is a pre-$\closure{D}$-vip order which
also satisfies the following condition for all $d\in\closure{D}$
and for all $u,v\in D\setminus\{d\}$:
\begin{enumerate}\addtocounter{enumi}{1}
\item if $d\meet u= d\meet v\ne d$ and $u(\lg(d))<v(\lg(d))$,
then $u\prec v$.
\end{enumerate}
\end{definition}

It is not difficult to construct a pre-$D$-vip order 
for a transverse set $D$. 

\begin{lemma}\label{lem.exist.vip}
If $D\subseteq \krtree$ is transverse, then
there is a pre-$D$-vip order of $\krtree$.
If $D\subseteq \krtree$ is a sparse diagonal set and $\prec$
is a pre-$S$-vip order for some $S$ with $\closure{D}\subseteq S$,
then $\prec$ is a $D$-vip order.
\end{lemma}

\begin{proof}
Let $\dless$ be any well-ordering of the levels of $\krtree$. 
Use recursion to define a pre-$D$-vip order $\prec$ by
adjusting $\dless$ separately on each level which has an
element of $D$.

Note that if $D\subseteq \krtree$ is a sparse diagonal set,
and $\prec$ is a pre-$S$-vip order for some $S$ with 
$\closure{D}\subseteq S$, then the second clause never
applies so $\prec$ is a $D$-vip order.
\end{proof}

Now that we have at least some idea of the level order
we will use, 
we define a family of ordered similarity types which
is rich enough to have all the ordered similarity types
of $m$-element subsets of sparse diagonal sets $D\subseteq\krtree$
for a $D$-vip level ordering $\prec$.

\begin{definition}\label{def.sparse.type}
Call $\tau$ an \emph{$m$-type} if it is
a downward closed subtree of ${}^{2m-2\ge}2$ whose 
set $L$ of leaves is an $m$-element strongly diagonal set.
Call $(\tau,\dless)$ a \emph{vip $m$-type} 
if $\tau$ is a $m$-type and $\dless$ is an $L$-vip order.
If $L$ is sparse, then $\tau$ is called a
\emph{sparse $m$-type} and $(\tau,\dless)$ is called 
a \emph{sparse vip $m$-type}. 
\end{definition}

\begin{lemma}\label{lem.sparse.type1}
Assume $D\subseteq\krtree$ is a strongly diagonal set
and $\prec$ is an ordering of the levels of $\krtree$
which is a $D$-vip order.  Then for all
$m$-element sets $x\subseteq D$, 
$(\clp(x),\prec_x)$ is a vip $m$-type,
and if $D$ is sparse, then $(\clp(x),\prec_x)$
is a sparse vip $m$-type. 
\end{lemma}

\begin{proof}
The collapse of any $m$-element strongly diagonal set $x$ is a subtree closed
under initial segments whose set of leaves, $L$, is strongly diagonal.
Moreover, since the order $\prec$ is $x$-vip, it follows that
$\prec_x$ is $L$-vip.  If, in addition, $D$ is sparse, then
$\clp(x)$ is also sparse, by definition of collapse. 
\end{proof}

In order to apply Shelah's Theorem, we shall need an
ordering of the levels of $\krtree$ and a conveniently chosen
subset of $\krtree$ the antichains of which will realize
the smallest possible number of weak embedding types. Toward that
end, we introduce \emph{cofinal transverse}
subsets of $\krtree$.

\begin{definition}\label{def.cofinal}
A subset $S\subseteq \krtree$ is \emph{cofinal above $w$}
if for all $t\in\cone(w)$ there is some $s\in S$
with $t\subseteq s$.   If $w=\emptyset$, we say $S$
is \emph{cofinal}.
\end{definition}

\begin{definition}\label{def.transverse}
Call $A\subseteq \krtree$  \emph{transverse} if distinct elements
of $A$
have different lengths.  
\end{definition}

By recursion one can construct a cofinal transverse subset
of $\krtree$.

\begin{lemma}\label{lem.cofinal.transverse}
There is a cofinal transverse subset of $\krtree$.
\end{lemma}

\begin{lemma}\label{lem.cofinal.dense}
If $S\subseteq\krtree$ is cofinal above $w$ and transverse,
then $(S\cap\cone(w),\qle)$ is $\kappa$-dense.
\end{lemma}

\begin{proof}
Suppose $A,B\subseteq S\cap\cone(w)$ are two disjoint subsets of size less
than $\kappa$ with $a<_Q b$ for all $a\in A$ and $b\in B$.
Use the fact that $\cone(w)$ is $\kappa$-dense to find 
$d\in \cone(w)\setminus(A\cup B)$ with $a\qle d\qle b$ for
all $a\in A$ and $b\in B$.  Let $\alpha$ be a limit ordinal 
larger than the length of any element of $A\cup B$.  Let $d'$ 
be the extension of $d\cat\langle 1\rangle$ by zeros of length 
$\alpha$ and let $c$ be an extension of $d'\cat\langle 0\rangle$ in $S$.  Since
$d'$ and $c$ are longer than any element of $A\cup B$, they are not in
$A$ nor in $B$.
Then $d\qle c\qle d'$ and $d'\qle b$ for all $b\in B$.  
Thus $c$ is the required witness showing $(S,\qle)$ is $\kappa$-dense.
\end{proof}

\begin{definition}\label{def.sparse.diag.map}
Assume $z,w\in\krtree$ and $S\subseteq\krtree$ is cofinal and transverse.
Call $g$ a 
\emph{sparse diagonalization of $\krtree$ into $S\cap\cone(w)$}
if $g:\krtree\to S\cap\cone(w)$ is an injective $\qle$-preserving map
such that $D:=g[\krtree]$ is a sparse diagonal subset,
$\closure{D}\subseteq S\cap\cone(w)$ and the following conditions hold:
\begin{enumerate}
\item\label{item.meet.lift.1} 
for all three element diagonal sets $\set{x,u,v}$,
if $x\meet u=x\meet v$, then
$g(x)\meet g(u)=g(x)\meet g(v)$;
\item\label{item.meet.length.1} 
for all $x,y,u,v\in\krtree$, if $\lg(x\meet y)<\lg(u\meet v)$,
then $\lg(g(x)\meet g(y))<\lg(g(u)\meet g(v))$;
\item for all sparse diagonal $E\subseteq\krtree$,
$\clp(E)=\clp(g[E])$.
\end{enumerate}
If $w=\emptyset$, then we call $g$ a
\emph{sparse diagonalization of $\krtree$ into $S$}.
\end{definition}

\begin{lemma}\label{lem.1.diag.map}[First Diagonalization Lemma]
Assume 
$\kappa=2^{<\kappa}$,
$w\in\krtree$ and $S\subseteq\krtree$ is cofinal and transverse.
Then there  is a sparse diagonalization 
$\varphi$ of $\krtree$ into $S\cap\cone(w)$.
\end{lemma}

\begin{proof}
Let $\prec$ be a well-ordering of the levels of $\krtree$.
Let $\seq{t_\alpha:\alpha<\kappa}$ list the elements of
$\krtree$ in $\prec$-increasing order (recall that $\kappa=
2^{<\kappa}$ is assumed).  

Define functions $\varphi_0,\varphi_1$ and $\varphi$
on $\krtree$ by recursion on $\alpha$.
To start the recursion,
notice that $t_0=\emptyset$, and let $\varphi_0(t_0)$ be
an element of $S\cap\cone(w)$ of minimal length.

If $\varphi_0(t_\alpha)$ has been defined, let
$\varphi_1(t_\alpha)$ be an element of $S\cap\cone(w)$
extending $\varphi_0(t_\alpha)\cat \langle 1\rangle$ and  
let $\varphi(t_\alpha)$ be an element 
of $S\cap\cone(w)$ extending
$\varphi_1(t_\alpha)\cat\langle 0\rangle$.  
Note that $\varphi(t_\alpha)\supseteq\varphi_1(t_\alpha)\supseteq
\varphi_0(t_\alpha)$.

Suppose $\varphi_0(t_\beta)$, $\varphi_1(t_\beta)$ and $\varphi(t_\beta)$
have been defined for all $\beta<\alpha$, and 
$t_\alpha=s_\alpha\cat\langle\delta\rangle$ for some $\delta$.
Let $\gamma_\alpha$ be the least ordinal
$\gamma$ strictly greater than $\lg(\varphi(t_\beta))$ for all 
$\beta<\alpha$.    
Then let $\varphi_0(t_\alpha)$ be an element of $S\cap\cone(w)$
extending the extension by zeros of 
$\varphi_\delta(s_\alpha)\cat\langle\delta\rangle$ of length $\gamma_\alpha$.

Suppose $\varphi_0(t_\beta)$, $\varphi_1(t_\beta)$ and $\varphi(t_\beta)$
have been defined for all $\beta<\alpha$, and 
$\lg(t_\alpha)=\zeta$ is a limit ordinal.
Let $\varphi^-(t_\alpha):=
\bigcup\set{\varphi_0(t_\alpha\restrict\eta):\eta<\zeta}$.
Let $\gamma_\alpha$ be the least ordinal $\gamma$
strictly greater than $\lg(\varphi^-(t_\alpha))$
and strictly greater than $\lg(\varphi(t_\beta))$ for all 
$\beta<\alpha$. Let $\varphi_0(t_\alpha)$ be an element of
$S\cap\cone(w)$ extending the
extension by zeros of $\varphi^-(t_\alpha)$
of length $\gamma_\alpha$.

By construction, $\varphi_0$ and hence $\varphi_1$ and $\varphi$ are 
injective.  Moreover, the ranges of all three are subsets of 
$S\cap\cone(w)$, so for $D:=\varphi[\krtree]$, the meet closure
satisfies $\closure{D}\subseteq \cone(w)$. Below we shall show that
$\closure{D}\subseteq S$.

For all $\alpha<\kappa$, one has
$\lg(\varphi_0(t_\alpha))<\lg(\varphi_1(t_\alpha))
<\lg(\varphi(t_\alpha))<\lg(\varphi_0(t_{\alpha+1}))$ and if
$\alpha$ is a limit then for all $\beta<\alpha$ we have
$\lg(\varphi(t_\beta))<\lg(\varphi(t_\alpha))$.  
Thus if $\alpha<\beta$, then
$\lg(\varphi(t_\alpha))<\lg(\varphi_0(t_\beta))$,
so different elements of the union of the ranges of
$\varphi_0$, $\varphi_1$ and $\varphi$ have different lengths.

Consider two distinct elements,  $t_\alpha\qle t_\beta$.
By a case analysis below, we show that their images are incomparable,
the $\qle$-order between $t_\alpha$ and $t_\beta$ is preserved by $\varphi$
and we show how to express the meet $\varphi(t_\alpha)\meet
\varphi(t_\beta)$ as one of $\varphi(t_\alpha\meet t_\beta)$,
$\varphi_0(t_\beta)$ and $\varphi_1(t_\alpha)$.

For the first case, suppose $t_\alpha$ and $t_\beta$ are incomparable.
Then $t_\alpha\lex t_\beta$.  Let $\gamma$ be such that 
$t_\gamma=t_\alpha\meet t_\beta$.
Then $t_\gamma\cat\langle 0\rangle\subseteq t_\alpha$
and  $t_\gamma\cat\langle 1\rangle\subseteq t_\beta$.
From the definition of $\varphi_0$, $\varphi_1$ and $\varphi$,
it follows that $\varphi_0(t_\alpha\meet t_\beta)=\varphi_0(t_\gamma)
=\varphi(t_\alpha)\meet\varphi(t_\beta)$ and 
$\varphi(t_\alpha)\lex\varphi(t_\beta)$, so 
$t_\alpha\qle t_\beta$ and
$\varphi(t_\alpha)\qle\varphi(t_\beta)$.

For the second case, suppose  $t_\beta\subseteq t_\alpha$.
Then,  since $t_\alpha$ and $t_\beta$ are distinct
and $t_\alpha\qle t_\beta$, it follows that
$t_\beta\cat\langle 0\rangle\subseteq t_\alpha$.
Consequently $t_\alpha\meet t_\beta=t_\beta$,
$\varphi(t_\alpha)\meet\varphi(t_\beta)=
\varphi_0(t_\beta)$,
and $\varphi(t_\alpha)\lex\varphi_1(t_\beta)
\subseteq\varphi(t_\beta)$,
so $\varphi(t_\alpha)\qle\varphi(t_\beta)$.

For the third case, suppose  $t_\alpha\subseteq t_\beta$.
Then,  since $t_\alpha$ and $t_\beta$ are distinct
and $t_\alpha\qle t_\beta$, it follows that
$t_\alpha\cat\langle 1\rangle\subseteq t_\beta$.
Consequently $t_\alpha\meet t_\beta=t_\alpha$,
$\varphi(t_\alpha)\meet\varphi(t_\beta)=
\varphi_1(t_\alpha)$.
Since $\varphi_1(t_\alpha)\cat\langle 0\rangle\subseteq \varphi(t_\alpha)$
and $\varphi_1(t_\alpha)\cat\langle 0\rangle\lex\varphi(t_\beta)$,
it follows that $\varphi(t_\alpha)\qle\varphi(t_\beta)$.

By the above analysis, $\varphi$ preserves the $\qle$-order
and sends distinct elements of $\krtree$ to incomparable 
elements.  Thus the image of $\varphi$ is an antichain.
Moreover, for $\eta<\zeta$, the meet, 
$\varphi(t_\eta)\meet\varphi(t_\zeta)$, is one of 
$\varphi_0(t_\eta\meet t_\zeta)$ and 
$\varphi_1(t_\eta\meet t_\zeta)$,
and the latter occurs only if $t_\eta\cat\langle 1\rangle\subseteq
t_\zeta$.  Thus different elements of the meet closure of the
image of $\varphi$ have different lengths.  
It follows that the image of $\varphi$ is diagonal.
Moreover, by construction, the meet closure of the image $D$
of $\varphi$ is a subset of $S\cap\cone(w)$.

To complete the proof of the lemma, we must show that the
image $D$ is sparse.  First use induction on
$\beta<\kappa$ to show that for all $\alpha<\beta$, the following
two statements hold:
\begin{enumerate}
\item $\varphi(t_\beta)(\lg(\varphi_0(t_\alpha)))=0$;
\item ($\varphi(t_\beta)(\lg(\varphi(t_\alpha)))=1$ or
$\varphi(t_\beta)(\lg(\varphi_1(t_\alpha)))=1$) 
if and only if \newline
($\varphi(t_\beta)(\lg(\varphi(t_\alpha)))=1$ and
$\varphi(t_\beta)(\lg(\varphi_1(t_\alpha)))=1$) 
if and only if \newline
($t_\alpha\cat\langle 1\rangle\subseteq  t_\beta$ and
$\varphi_1(t_\alpha)\cat\langle 1\rangle\subseteq \varphi(t_\beta)$).
\end{enumerate}

Next suppose $s$ and $t$ in the meet closure of $\varphi[\krtree]$
are such that $\lg(t)>\lg(s)$ and $t(\lg(s))=1$.  Without loss
of generality, we assume that $t$ is in the image of $\varphi$,
since we know it has an extension in the image, and let
$t_\beta$ be such that $t=\varphi(t_\beta)$.
By construction, since $t(\lg(s))=1$,
either $s=\varphi_0(t_\alpha)$ or $s=\varphi_1(t_\alpha)$
for some $\alpha\le\beta$.  If $\alpha<\beta$, then
$s\cat\langle 1\rangle\subseteq t$ by the second statement
above since $\varphi_0(t_\alpha)\cat\langle 1\rangle\subseteq
\varphi_0(t_\alpha)$.
If $\alpha=\beta$, then $s=\varphi_0(t_\alpha)=\varphi_0(t_\beta)$
since $\varphi(t_\beta)$ is an extension of 
$\varphi_1(t_\beta)\cat\langle 0\rangle$.
In this case, $s\cat\langle 1\rangle\subseteq t$ as well.
Thus the image of $D$ is sparse.
\end{proof}

\begin{lemma}\label{lem.embed.sparse}
Suppose $\prec$ is ordering of the levels of $\krtree$.  
If  $e:\krtree\to\krtree$ is a strong embedding which preserves
$\prec$ on $\cone(w)$ and $A\subseteq\cone(w)$ 
is an $m$-element sparse diagonal subset, then
$(\clp(A),\prec_A)=(\clp(e[A]),\prec_{e[A]})$.
\end{lemma}

\begin{proof}
Under the hypotheses of the lemma, since $e$ is a strong embedding,
$\clp(e[A])=\clp(A)$.  Since $e$ preserves order on $\cone(w)$,
$\prec_{e[A]}=\prec_A$.
\end{proof}

\begin{theorem}\label{thm.upper.bound.1} 
Let $m\ge 2$ and
  suppose that $\kappa$ is a cardinal which is measurable
  in the generic
  extension obtained by adding $\lambda$ Cohen subsets of $\kappa$,
  where $\lambda\to(\kappa)^{2m}_{2^{\kappa}}$. 
Then for $t_m^+$ equal to the number of vip $m$-types,
\[
\mathbb{Q}_\kappa\rightarrow(\mathbb{Q}_\kappa)^m_{<\kappa,t_m^+}.
\]
\end{theorem}

\begin{proof}
Suppose $d:[\krtree]^m\to\mu$ is a fixed coloring for some $\mu<\kappa$.  
Use Lemma \ref{lem.cofinal.transverse} to find
$S\subseteq \krtree$ cofinal and transverse.
Lemma \ref{lem.exist.vip} to find a pre-$S$-vip order
$\prec$ of the levels of $\krtree$.
Apply Shelah's
Theorem \ref{thm.SHL} to the restriction to antichains to obtain
a strong embedding $e$ and a node $w$ such that $e$ preserves
$\prec$ on $\cone(w)$ and $d$ is constant on $m$-element subsets
of the same $\prec$-ordered similarity type.

Apply the First Diagonalization Lemma \ref{lem.1.diag.map}
to find a sparse diagonalization $\varphi$ of $\krtree$
into $S\cap\cone(w)$.  Let $D=\varphi[\krtree]$.
Note that $\prec$ is a $D$-vip order, since $\prec$ is a pre-$S$-vip order. 
By Lemma \ref{lem.sparse.type1}, all $\prec$-ordered
similarity types of $m$-element subsets of $D$ are
sparse vip $m$-types.  Let $Q=e[D]$.  By Lemma \ref{lem.embed.sparse},
all $\prec$-ordered similarity types of $m$-element subsets
of $Q$ are sparse vip $m$-types.  

Since $\varphi$ is a sparse diagonalization and $e$
a strong embedding, $(D,\qle)$ and $(Q,\qle)$ are both $\kappa$-dense.
Since $d$ is constant on $m$-element subsets of $Q$ of the same
$\prec$-ordered similarity type, it follows that
$d$ takes on no more colors
than the number $t_m^+$ of sparse vip $m$-types.
\end{proof}

\section{An almost perfect subtree}\label{sec.almost.perfect}
In this section, in preparation for
computing some small values of $r_n^+$,
we show that, for an
arbitrary $S\subseteq \krtree$ with $(S,\qle)$ $\kappa$-dense,
the set of all nodes $w$ with a $\kappa$-dense
set of extensions in $S$ forms an \emph{almost perfect subtree}
(defined later in this section).
We use the almost perfect subtree
to construct a $\kappa$-dense diagonal subset of an
arbitrary $S\subseteq \krtree$ with $(S,\qle)$ $\kappa$-dense.

\begin{lemma}\label{lem.grow}
Suppose $\kappa$ satisfies
$\kappa^{<\kappa}=\kappa$ (so $\kappa$ is a regular cardinal),
$S\subseteq \krtree$, and $(S,\qle)$ is a $\kappa$-dense
linear order.  For all $w\in \krtree$,
the set $S\cap\cone(w)$ is either empty, a singleton or $\kappa$-dense.
\end{lemma}

\begin{proof}
Fix $w\in\krtree$.  If $S\cap\cone(w)$ is empty or a singleton,
there is nothing to prove.  So suppose $x\in S$ and $y\in S$
are two different extensions of $w$.   Notice that $w\subseteq x\meet y$.
Without loss of generality, assume $x\qle y$.

If $x$ and $y$ are incomparable and $x\qle z\qle y$,
then either $x\cat\langle 1\rangle\subseteq z$ or
$y\cat\langle 0\rangle \subseteq z$ or
$x\lex z\lex y$.  In all three cases $w\subseteq z$.

If $x$ and $y$ are comparable, then either
$y\cat\langle 0\rangle\subseteq x$ or
$x\cat\langle 1\rangle\subseteq y$.  Hence
for any $z$ with $x\qle z\qle y$, one of
$x$ and $y$ is a subset of $z$.  In either
case $w\subseteq z$.   

Since $S$ is $\kappa$-dense, it has a
$\kappa$-dense subset $S'$ with
$\set{x}<S'<\set{y}$.  By the above
two paragraphs, every element of $S'$
is an extension of $w$, so
$S\cap\cone(w)$ is $\kappa$-dense.
\end{proof}

\begin{definition}\label{def.T(S)}
Suppose $\kappa^{<\kappa}=\kappa$,
$S\subseteq\krtree$ and $(S,\qle)$ is $\kappa$-dense.
Let $T(S)$ be the set of all $t\in\krtree$ for which
$S\cap\cone(t)$,
is $\kappa$-dense.
\end{definition}

We plan to show that $T(S)$ is an almost perfect tree.
The first step is to show that arbitrarily high above every node there
is a \emph{densely splitting node}.

\begin{definition}\label{def.W(S)}
Define $\mathcal{W}:\wp(\krtree)\to\wp(\krtree)$ by
\[
\mathcal{W}(S):=\set{w\in \closure{S}: 
\text{$S\cap\cone(w\cat\langle 0\rangle)$ and 
$S\cap\cone(w\cat\langle 1\rangle)$ are $\kappa$-dense}},
\]
and call the elements of $\mathcal{W}(S)$ \emph{densely splitting nodes of $S$}.
\end{definition}

\begin{lemma}\label{lem.split.continue}
Suppose $\kappa^{<\kappa}=\kappa$,
$S\subseteq \krtree$ and $(S,\qle)$ is a $\kappa$-dense
linear order.  Then $\mathcal{W}(S)$ is non-empty.
\end{lemma}

\begin{proof}
Assume toward a contradiction that for all $u\in\krtree$,
one or both of $(S\cap\cone(u\cat\langle 0\rangle,\qle)$ and
$(S\cap\cone(u\cat\langle 1\rangle),\qle)$ have cardinality less than $\kappa$.
By Lemma \ref{lem.grow}, it follows that for all $u\in\krtree$,
one or both of $S\cap\cone(u\cat\langle 0\rangle)$ and
$S\cap\cone(u\cat\langle 1\rangle)$ have cardinality less than $2$.

Define $\seq{t_\alpha:\alpha<\kappa}$ by recursion and prove
by induction that $S\cap\cone(t_\alpha)$ is $\kappa$-dense for all $\alpha$.

To start the recursion, let $t_0=\emptyset$.  Then
$S\cap\cone(t_0)=S$ which is $\kappa$-dense.  

If $\alpha=\beta+1$ and 
$S\cap\cone(t_\beta)$ is $\kappa$-dense, then let $t_\alpha=t_\beta\cat
\langle 0\rangle$ if $S\cap\cone(t_\beta\cat\langle 0\rangle)$ is
$\kappa$-dense, and $t_\alpha=t_\beta\cat\langle 1\rangle$ otherwise.
Since $S\cap \cone(t_\beta)\subseteq (S\cap \cone(t_\beta\cat\langle 0\rangle))
\cup (S\cap\cone(t_\beta\cat\langle 1\rangle))$,
by Lemma \ref{lem.grow} and a cardinality argument, 
$S\cap\cone(t_\alpha)$ is $\kappa$-dense.  

If $\alpha$ is a limit ordinal,
let $t_\alpha=\bigcup\set{t_\beta:\beta<\alpha}$.  By assumption,
for all $\beta<\alpha$, if 
$t_\beta\cat\langle\delta\rangle\not\subseteq t_{\beta+1}$, 
then $S\cap\cone(t_\beta\cat\langle\delta\rangle)$ is not $\kappa$-dense,
so has cardinality less than $2$.  It follows that
$|S\setminus \cone(t_\alpha)|<\kappa$.  Hence by Lemma \ref{lem.grow},
$S\cap\cone(t_\alpha)$ is $\kappa$-dense.  

Note that $b=\bigcup\set{t_\alpha:\alpha<\kappa}$ is a branch
through $\krtree$.  Every element $s$ of $S$ is either an initial
segment of this branch or there is some $\beta<\kappa$ such that
$s\meet t_{\beta+1}=t_\beta$.
Recall our assumption that  
for all $\beta<\alpha$, if 
$t_\beta\cat\langle\delta\rangle\not\subseteq t_{\beta+1}$, 
then $S\cap\cone(t_\beta\cat\langle\delta\rangle)$ is not $\kappa$-dense,
so it has cardinality less than $2$.  It follows that
$|S|<\kappa$, contradicting the assumption that $S$ is $\kappa$-dense.  
Thus the lemma follows.
\end{proof}

\begin{lemma}\label{lem.room}
Suppose $\kappa^{<\kappa}=\kappa$,
$S\subseteq \krtree$ and
$(S,\qle)$ is a $\kappa$-dense linear order.
Then for every $w\in \mathcal{W}(S)$, every $\delta<2$,
and every $\alpha<\kappa$, there is $u\in \mathcal{W}(S)$ such that 
$w\cat\langle \delta\rangle\subseteq u$ and $\lg (u)\ge\alpha$.
\end{lemma}

\begin{proof}
Fix $w\in \mathcal{W}(S)$, $\delta<2$ and $\alpha<\kappa$.  
Without loss of
generality, assume $\lg(w\cat\langle\delta\rangle)<\alpha$.
Since $|S\cap\cone(w\cat\langle\delta\rangle)|=\kappa$, it follows
that $S\cap\cone(w\cat\langle\delta\rangle)$ is $\kappa$-dense
by Lemma \ref{lem.grow}.  

Since $\kappa^{<\kappa}=\kappa$, the set of nodes of
$\krtree$ of length at most $\alpha$ has cardinality
less than $\kappa$, but $S\cap\cone(w\cat\langle\delta\rangle)$
has cardinality $\kappa$ because it is $\kappa$-dense.

By the pigeonhole principle, there is some $u_0\in {}^\alpha 2$
which has $\kappa$ many extensions in 
$S\cap\cone(w\cat\langle\delta\rangle)$.  It follows that
$w\cat\langle\delta\rangle\subseteq u_0$ and that
$S\cap\cone(u_0)$ is $\kappa$-dense by Lemma \ref{lem.grow}.  
By Lemma \ref{lem.split.continue}
there is an extension $u\supseteq u_0$ in $\mathcal{W}(S)$,
so the lemma follows.
\end{proof}

\begin{lemma}\label{lem.T(S)}
Suppose $\kappa^{<\kappa}=\kappa$,
$S\subseteq\krtree$ and $(S,\qle)$ is $\kappa$-dense.
Then $(T(S),\subseteq)$ is a rooted tree, 
$\mathcal{W}(S)\subseteq T(S)$, and for all $t\in T(S)$, for all $\alpha>\lg(t)$,
there is some $r\in\mathcal{W}(S)$ with $t\subseteq r$ and $\lg(r)\ge\alpha$.
\end{lemma}

\begin{proof}
Since $T(S)$ is closed under initial segments, $(T(S),\subseteq)$
is a rooted tree.  By definition of $\mathcal{W}(S)$, it is a subset of $T(S)$.
By Lemma \ref{lem.room}, 
every element of $T$ has extensions of arbitrarily large
length which are densely splitting nodes.
\end{proof}

To prove that $T(S)$ is an almost perfect tree, we need to
be able to prove that it has certain continuity properties
at limit levels.  Toward that end, we introduce the notion
of a \emph{limit of densely splitting nodes of $S$}.

\begin{definition}\label{def.robust}
Suppose $\kappa$ is a regular cardinal, $S\subseteq\krtree$,
and $(S,\qle)$ is $\kappa$-dense.
Say $t\in \krtree$ \emph{is a limit of densely splitting nodes of $S$} if 
$\lg(t)$ is a limit ordinal and for
unboundedly many $\beta<\lg(t)$, $t\restrict\beta$
is in $\mathcal{W}(S)$.  
If $t$ is a limit of densely splitting nodes of $S$, then say it is
\emph{evenhanded in $S$} if for $\delta=0,1$, the set
$\set{\beta<\lg(t): t\restrict\beta\in \mathcal{W}(S)\land t(\beta)=\delta}$
is unbounded in $\lg(t)$.
\end{definition}

\begin{lemma}\label{lem.robust}
Suppose $\kappa^{<\kappa}=\kappa$,
$S\subseteq \krtree$ and
$(S,\qle)$ is a $\kappa$-dense linear order.
Further suppose that $z$
is a limit of densely splitting nodes of $S$.  
If $z$ is evenhanded in $S$,
then $z$ has an extension in $\mathcal{W}(S)$.
\end{lemma}
 
\begin{proof}
For each $\beta<\lg(z)$ with $z\restrict\beta\in S$, 
let $g(\beta)$ be an element of $S$ extending
$z\restrict\beta\cat\langle 1-z(\beta)\rangle$.
Let $A$ be the set of $g(\beta)$ for $\beta<\lg(z)$
with $z\restrict\beta\in S$ and $z(\beta)=1$,
and let $B=\ran (g)\setminus A$.

Then $A \qle\set{z}\qle B$.

\begin{claim}
If $w\in S$ and $A\qle\set{w}\qle B$, then $z\subseteq w$.
\end{claim}

\begin{proof}
Suppose $w\in S$ and $A\qle\set{w}\qle B$.

Assume toward a contradiction that $\gamma:=\lg(w\meet z)<\lg(z)$.
Use the fact that $z$ is evenhanded in $S$ to choose
$\eta$ and $\theta$ strictly greater than $\gamma:=\lg(w\meet z)$
such that $z\restrict\eta\in \mathcal{W}(S)$, 
$z(\eta)=1$, $z\restrict\theta\in \mathcal{W}(S)$
and $z(\theta)=0$.
Then $g(\eta)\in A$ and $g(\theta)\in B$, so $g(\eta)\qle w\qle g(\theta)$.

By definition of $g$, $z\restrict\eta\subseteq g(\eta)$
and $z\restrict\theta\subseteq g(\theta)$.
Consequently $z\restrict(\gamma+1)\subseteq g(\eta)\meet g(\theta)$,
and $w\meet g(\eta)=w\meet z=w\meet g(\theta)$.
It follows that both $g(\eta)$ and $g(\theta)$ are in the same
$\qle$-relationship with $w$ as is $z$, which is not possible because
either $w\qle g(\eta)$ and $w\qle g(\theta)$ contradicting the
fact that $g(\eta)\qle w\qle g(\theta)$.
Thus $\lg(w\meet z)\ge\lg(z)$.  In other
words, $z\subseteq w$.   
\end{proof}

Since $S$ is $\kappa$-dense, there is a $\kappa$-dense set $C\subseteq S$
with $A\qle C\qle B$.  By the claim, $C\subseteq S\cap\cone(z)$.
It follows that $S\cap\cone(z)$ is $\kappa$-dense, so by Lemma 
\ref{lem.split.continue}, $z$ has an extension in $\mathcal{W}(S)$.
\end{proof}

\begin{definition}\label{def.favor}
Suppose $\kappa^{<\kappa}=\kappa$,
$S\subseteq\krtree$ and $(S,\qle)$ is $\kappa$-dense.
Say $t\in\krtree$ \emph{favors $\delta$ above $s$ in $T(S)$}
if $s\subsetneq t$, either $t\in T(S)$ or $t$ is a limit of densely splitting 
nodes of $S$,
and for all $r\in \mathcal{W}(S)$ with $s\subsetneq r\subsetneq t$,
$r\cat\langle\delta\rangle\subseteq t$.
\end{definition}

\begin{lemma}\label{lem.favor.alt}
Suppose $\kappa^{<\kappa}=\kappa$,
$S\subseteq\krtree$ and $(S,\qle)$ is $\kappa$-dense.
Further suppose $u_0\lex u_1$, $u_0$ favors $1$ above $u_0\meet u_1$
in $T(S)$, and $u_1$ favors $0$ above $u_0\meet u_1$ in $T(S)$.
Then at least one of $u_0$ and $u_1$ is in $T(S)$.
\end{lemma}

\begin{proof}
If one of $u_0$ and $u_1$ is not a limit of densely splitting nodes 
of $S$, then it is in $T(S)$ by definition of favor.  So assume 
both are limits of densely splitting nodes of $S$.

For each $r\in\mathcal{W}(S)$ with $u_0\meet u_1\subsetneq r \subsetneq
u_0$, let $r^-\in S$ be an extension of $r\cat\langle
0\rangle$, and let
$A$ be the set of these elements.
Similarly, for 
each $r\in\mathcal{W}(S)$ with $u_0\meet u_1\subsetneq r \subsetneq
u_1$, let $r^+\in S$ be an extension of $r\cat\langle
1\rangle$, and let
$B$ be the set of these elements.
Then $A\lex\set{u_0}\qle \set{u_0\meet u_1}\qle\set{u_1}\lex B$.

\begin{claim}
If $A\qle\set{w}\qle\set{u_0\meet u_1}$ and $w\ne u_0\meet u_1$,
then $w$ extends $(u_0\meet u_1)\cat\langle 0\rangle$
and either $w\subseteq u_0$ or $u_0\subseteq w$ or $u_0\lex w$.
\end{claim}

\begin{proof}
Suppose $w$ satisfies the hypotheses.  Then by definition of $\qle$,
$(u_0\meet u_1)\cat\langle 0\rangle\subseteq w$.
Assume toward a contradiction that
none of the three conclusions holds.  Then
$w\lex u_0$.
Let $r\subseteq u_0$ be an element of $\mathcal{W}(S)$
with $u_0\meet u_1\subsetneq r$ and $\lg(r)>\lg(u\meet w)+1$.
Then $w\lex r$, so $w\lex r^-$, contradicting $A\qle\set{w}$.
\end{proof}

\begin{claim}
The set of all $w\in S$ such that $A\qle\set{w}\qle\set{u_0\meet u_1}$,
and $u_0\lex w$ has cardinality $<\kappa$.
\end{claim}

\begin{proof}
Otherwise, by the pigeonhole principle and the previous claim,
there is some $r\subsetneq u_0$ with $(u_0\meet u_1)\cat\langle 0\rangle
\subseteq r$ such that the set of $w\in S$ with $u_0\lex w$ and
$u_0\meet w=r$ has cardinality $\kappa$.  It follows that $u_0(\lg(r))=0$,
and, by Lemma \ref{lem.grow}, that $r\cat\langle 1\rangle \in T(S)$.
Thus $r\in \mathcal{W}(S)$ contradicts the assumption that $u_0$ favors
$1$ above $u_0\meet u_1$.
\end{proof} 

The proofs of the next two claims are similar to those
above, so they are left to the reader.

\begin{claim}
If $\set{u_0\meet u_1}\qle\set{w}\qle B$ and $w\ne u_0\meet u_1$,
then $w$ extends $(u_0\meet u_1)\cat\langle 1\rangle$
and either $w\subseteq u_1$ or $u_1\subseteq w$ or $w\lex u_1$.
\end{claim}

\begin{claim}
The set of all $w\in S$ such that $A\qle\set{w}\qle\set{u_0\meet u_1}$,
and $u_0\lex w$ has cardinality $<\kappa$.
\end{claim}

Let $C\subseteq S$ be a $\kappa$-dense subset with
$A\qle C\qle B$.  Since the inequalities $A\qle\set{u_0\meet u_1}\qle B$ hold,
either $C^-:=\set{c\in C:c\qle u_0\meet u_1}$ or
$C^+:=\set{c\in C:u_0\meet u_1\qle c}$ has cardinality
$\kappa$.  If $C^-$ has cardinality $\kappa$, then
$S\cap\cone(u_0)\cat C^-$ has cardinality $\kappa$, so by Lemma
\ref{lem.grow}, $u_0$ is in $T$.  Similarly, if
$C^+$ has cardinality $\kappa$, then $u_1$ is in $T$.
Thus the lemma follows.
\end{proof}

\begin{lemma}\label{lem.min.W(S)}
Suppose $\kappa^{<\kappa}=\kappa$,
$S\subseteq\krtree$ and $(S,\qle)$ is $\kappa$-dense.
For all $t\in T(S)$ there is a minimal extension of
$t$ in $\mathcal{W}(S)$.  That is, there is $r\in\mathcal{W}(S)$ 
such that $t\subseteq r$ and $r\subseteq w$
for all $w\in \mathcal{W}(S)\cap\cone(t)$.
\end{lemma}

\begin{proof}
By Lemma \ref{lem.grow}, $\mathcal{W}(S)\cap\cone(t)$ is non-empty.
By definition of $\mathcal{W}(S)$, the meet of two incomparable
elements of it is also in the set.  It follows that
$\mathcal{W}(S)\cap\cone(t)$ has an element of minimum length,
and this element is the desired minimal extension.
\end{proof}

\begin{lemma}\label{lem.favor.split}
Suppose $\kappa^{<\kappa}=\kappa$,
$S\subseteq\krtree$ and $(S,\qle)$ is $\kappa$-dense.
For all $x\in \mathcal{W}(S)$ and all $\alpha>\lg(x)$,
there is some extension $y$ of $x$ with $\lg(y)=\alpha$
such that $y$ favors one of $0$,  $1$ above $x$.
\end{lemma}

\begin{proof}
Fix $x$ in $\mathcal{W}(S)$.  

For as long as possible, define a $\subseteq$-increasing sequence
$u_\alpha$ of extensions of $x\cat\langle 1\rangle$ which favor $0$ above $x$,
with $\lg(u_\alpha)=\alpha$.
To start the recursions with $\alpha=\lg(x)+1$,
let $u_\alpha =x\cat\langle 1\rangle\in T(S)$.
If $\alpha$ is a limit ordinal and $u_\beta$ has been defined
for $\lg(x)<\beta<\alpha$, then let $u_\alpha=\bigcup\set{u_\beta:
\lg(x)<\beta<\alpha}$.
If $\alpha=\beta+1$, $u_\beta$ has been defined, and $u_\beta\in T(S)$,
then let $u_\alpha=u_\beta\cat\langle 0\rangle$ if $u_\beta\in \mathcal{W}(S)$,
and otherwise let $u_\alpha$ be the one-point extension of $u_\beta$
which is a subset of the extension of $u_\beta$ of minimal length in
$\mathcal{W}(S)$.
If $\alpha=\beta+1$, $u_\beta$ has been defined, and $u_\beta\notin T(S)$,
then $\beta$ is a limit ordinal, $u_\beta$ is a limit of densely 
splitting nodes
of $S$ and the recursion stops with its definition.

Also, for as long as possible, define a $\subseteq$-increasing sequence
$v_\alpha$ of extensions of $x\cat\langle 0\rangle$ which favor $1$ above $x$,
with $\lg(v_\alpha)=\alpha$.
To start the recursion with $\alpha=\lg(x)+1$,
let $v_\alpha =x\cat\langle 0\rangle\in T(S)$.
If $\alpha$ is a limit ordinal and $v_\beta$ has been defined
for $\lg(x)<\beta<\alpha$, then let $v_\alpha=\bigcup\set{v_\beta:
\lg(x)<\beta<\alpha}$.
If $\alpha=\beta+1$, $v_\beta$ has been defined, and $v_\beta\in T(S)$,
then let $v_\alpha=v_\beta\cat\langle 0\rangle$ if $v_\beta\in \mathcal{W}(S)$,
and otherwise let $v_\alpha$ be the one-point extension of $v_\beta$
which is a subset of the extension of $v_\beta$ of minimal length in
$\mathcal{W}(S)$.
If $\alpha=\beta+1$, $v_\beta$ has been defined, and $v_\beta\notin T(S)$,
then $\beta$ is a limit ordinal, $v_\beta$ is a limit of densely
splitting nodes of $S$ and
the recursion stops with its definition.

By construction, for all $\alpha$ with $u_\alpha$ defined,
$u_\alpha$ favors $0$ above $x$ in $T(S)$.
Similarly, for all $\alpha$ with $v_\alpha$ defined, 
$v_\alpha$ favors $1$ above $x$ in $T(S)$.
If one of the recursions continues for all $\alpha<\kappa$,
the lemma follows.

Assume toward a contradiction that $u_\beta$ is defined but 
$u_{\beta+1}$ is not, and that $v_\beta$ is defined but
$v_{\beta+1}$ is not.  Then $u_\beta\notin T(S)$ and
$v_\beta\notin T(S)$,
contradicting Lemma \ref{lem.favor.alt}.
\end{proof}

\begin{definition}\label{def.almost.perfect}
A subset $T\subseteq \krtree$ is an \emph{almost perfect tree}
if it is a rooted induced subtree of $\krtree$ closed under initial
segments such that the following conditions hold:
\begin{enumerate}
\item for all $t\in T$, for all $\alpha<\kappa$, 
there is an extension $w\in T$ of $t$ of length at least $\alpha$
which is a densely splitting node (both $w\cat\langle 0\rangle$ and
$w\cat\langle 1\rangle$ are in $T$);
\item for all $s\in \krtree$, if $s$ is an evenhanded limit
of densely splitting nodes of $T$, then $s$ is in $T$;
\item for all $s,t\in\krtree$, if $\lg(s)=\lg(t)$, both $s$ and $t$ 
are limits of densely splitting nodes of $T$, and for some $x$,
$s$ favors $0$ above $x$ and $t$ favors $1$ above $x$,
then one of $s$ and $t$ is in $T$.
\end{enumerate}
\end{definition} 

%

\begin{lemma}\label{lem.1.almost.perfect}
Suppose $\kappa^{<\kappa}=\kappa$,
$S\subseteq\krtree$ and $(S,\qle)$ is $\kappa$-dense.
Then $T(S)$ is almost perfect.
\end{lemma}

\begin{proof}
Apply Lemmas \ref{lem.T(S)},  \ref{lem.robust} and \ref{lem.favor.split}.
\end{proof}

\begin{lemma}\label{lem.2.almost.perfect}
Suppose $\kappa^{<\kappa}=\kappa$
and $T\subseteq \krtree$ is almost perfect. Let $C$ be
the set of all limit ordinals $\alpha>0$ such that
every node of $T$ of length less than $\alpha$ extends to
a densely splitting node of $T$ of length less than $\alpha$.
Then $C$ is closed unbounded in $\kappa$,
and for all $\alpha\in C$, for all $u\in T\cap {}^{\alpha>}2$, the
set $T_\alpha(u):=\set{t\in T: u\subseteq t\land \lg(t)=\alpha}$ has
cardinality at least $2^{\cof{(\alpha)}}$.
\end{lemma}

\begin{proof}
Use the definition of almost perfect to show that $C$ is
non-empty and unbounded.  It follows immediately
from the definition of $C$ that it is closed.

Fix attention on $\alpha\in C$ and $u\in T\cap {}^{\alpha>}2$.
Let $\lambda=\cof (\alpha)$,
and suppose $\sigma\in {}^\lambda 2$ is a sequence 
such that if $\eta=\theta+k<\lambda$ for $\theta=0$ or $\theta$ limit,
and $k\equiv \delta$ mod $3$ for $\delta<2$, then $\sigma(\eta)=\delta$.

Define $\seq{r_\sigma(\eta):\eta<\lambda}$
by recursion and show by induction that
every element of it has length $<\alpha$.  

To start the recursion, let $r_\sigma(0)$ be the
minimal densely splitting node of $T$ extending $u$. It has length less than
$\alpha$ by the definition of $C$.  

Suppose $0<\eta<\lambda$ and $r_\sigma$ has been defined on elements
smaller than $\eta$.
If $\eta=\zeta+1$, let
$r_\sigma(\eta)$ be a densely splitting node of $T$
which extends $r_\sigma(\zeta)\cat\langle
\sigma(\zeta)\rangle$ of length less than $\alpha$.
Such a node exists by the definition of $C$.

If $\eta$ is a limit ordinal, then 
$r'_\sigma(\eta):=\bigcup\set{r_\sigma(\zeta):\zeta<\eta}$
is a limit of densely splitting nodes of $T$ 
and has length a limit ordinal less than $\alpha$
since $\eta<\lambda=\cof {(\alpha)}$.  The properties of $\sigma$
guarantee that this union is evenhanded, so by the definition of an
almost perfect tree, $r'_\sigma(\eta)$ is in $T$.  
Let $r_\sigma(\eta)$ be a densely splitting node of length less than $\alpha$
extending $r'_\sigma(\eta)$, which exists
by definition of $C$.

This completes the definition of $\seq{r_\sigma(\eta):\eta<\lambda}$.
Let $s_\sigma$ be the union of this sequence.  Then $s_\sigma$ is 
a limit of densely splitting nodes in $S$ and is evenhanded, 
so it is an element of $T$.  Since $r_\sigma(0)$ is an extension
of $u$, so is $s_\sigma$.

If $s_\sigma$ has length $\alpha$, then
let $t_\sigma=s_\sigma$ be this union, and notice
that it is in $T$.  Otherwise let $t_\sigma$ be an
extension of $s_\sigma$ of length $\alpha$ in $T$,
which must exist by the definition of an almost perfect tree.

Notice that if $\sigma,\tau\in{}^\lambda 2$ are two distinct sequences
with the property that  
$\eta=\theta+k<\lambda$ for $\theta=0$ or $\theta$ limit
and $k\equiv \delta$ mod $3$ for $\delta<2$ implies
$\sigma(\eta)=\tau(\eta)=\delta$, then $t_\sigma\ne t_\tau$.

Since no constraints have been placed on $\sigma(\theta+k)$ for
$k\equiv 2$ mod $3$, 
there are $2^\lambda$ sequences $\sigma\in {}^\lambda 2$ with the
special property described above. Thus the set $T_\alpha(u)$ has 
cardinality at least $2^\lambda$, and the lemma follows.
\end{proof}

\begin{lemma}\label{lem.nearly.perfect}
Suppose $\kappa^{<\kappa}=\kappa$,
$S\subseteq\krtree$ and $(S,\qle)$ is $\kappa$-dense.
Let $C(S)$ be the set of all limit ordinals $\alpha>0$ such that
every $t\in T(S)\cap {}^{\alpha>}2$ has proper extensions
in both $S\cap {}^{\alpha>}2$ and $\mathcal{W}(S)\cap {}^{\alpha>}2$.
Then $C(S)$ is closed unbounded in $\kappa$.
\end{lemma}

\begin{proof}
Use the definition of $T$ and Lemmas \ref{lem.grow} 
and \ref{lem.split.continue} to show that $C(S)$ is
non-empty and unbounded.  It follows immediately
from the definition of $C(S)$ that it is closed.
\end{proof}

\begin{lemma}\label{lem.diag.in.S}
Suppose that $\kappa$ is an inaccessible limit of inaccessible cardinals 
and $S$ is a subset of $\krtree$ with $(S,\qle)$ $\kappa$-dense.
Then there is a diagonal set $D\subseteq S$ such that
$(D,\qle)$ is $\kappa$-dense. 
\end{lemma}

\begin{proof}
Let $\prec$ be a total order on $\krtree$ satisfying $\lg(s)<\lg(t)\implies s
\prec t$. Define $r: T(S)\to\mathcal{W}(S)$
by setting $r(t)$ to be the minimal extension of $t$ in $\mathcal{W}(S)$.
By Lemma \ref{lem.min.W(S)}, $r$ is
well-defined.  Let $s:T(S)\to S$ be such that $s(t)$ is an extension
of $t$.

By recursion on $\prec$, define 
functions $\ell:\krtree\to\kappa$, 
$f_0,f_1:\krtree\to \mathcal{W}(S)$, and
$g:\krtree\to S$  as follows.

To start the recursion, note that the $\prec$-least element is $\emptyset$, 
define $\ell(\emptyset):=0$, and let 
$f_0(\emptyset):=r(\emptyset)$.
If $f_0(z)$  has been defined, let 
$f_1(z):=r(f_0(z)\cat\langle 1\rangle)$, and set
$g(z):=s(f_1(z)\cat\langle 0\rangle)$.

To continue the recursion, suppose 
$z\in\krtree$ has $\lg(z)=\alpha$
and for all $x\prec z$, both $\ell(x)$ and
$f_0(x)$ have been defined, with $f_1$ and $g$ defined from them as above.
Let $\ell(z)$ be the least ordinal greater than
$\lg(g(x))$ for all $x\prec z$.

If $\alpha=\beta+1$ is a successor and $z=y\cat\langle \delta\rangle$
for some $y$, then let
$f_0(z)$ be an extension in $\mathcal{W}(S)$
of $r(f_\delta(y)\cat\langle\delta\rangle)\cat\langle 1-\delta\rangle$
of length greater than $\ell(z)$.
Since $T(S)$ is almost perfect and $\mathcal{W}(S)$
is the set of densely splitting nodes of $T(S)$, such a node exists. 

If $\alpha>0$ is a limit ordinal, then 
$f^-(z):=\bigcup\set{f(z\restrict\beta):\beta<\alpha}$
has length a limit ordinal.
Since $f^-(z)$ extends $f_0(z\restrict\eta)$ 
for all $\eta<\alpha$, the definition of $f_0$ 
on nodes of successor length guarantees that $f^-(z)$ 
is a limit of densely splitting nodes of $S$
and evenhanded.  Hence by Lemma \ref{lem.robust}, it
has an extension in $\mathcal{W}(S)$.  It follows that $f^-(z)$
is in $T(S)$.  Let  $f_0(z)$ be an extension of $f^-(z)$ in $\mathcal{W}(S)$
of length greater than $\ell(z)$.

Let $D$ be the range of $g$.  Then $D$ is a subset of $S$,
and the meet closure of $D$ is a subset of union of the ranges 
of $f_0$, $f_1$ and $g$.  Since $f_0(x)\subsetneq f_1(x)\subsetneq
g(x)$, the lengths of these sequences are strictly increasing.
Also, by construction, if $x\prec z$, then $\lg(g(x))<\lg(f_0(z))$.
So different elements of the meet closure of $D$ have
different lengths.

By induction, one can show that $f_0$ preserves $\subseteq$ and
length order.  By definition of $f_0$,
$f_1$ and $g$, we have the  following properties:
\begin{enumerate}
\item if $x$ and $y$ are incomparable with $x\lex y$,
then $g(x)\meet g(y)=f_0(x\meet y)$ and $g(x)\lex g(y)$; 
\item if $x\cat\langle 0\rangle\subseteq y$, then
$g(x)\meet g(y)=f_0(x)$ and $g(y)\lex g(x)$;
\item if $x\cat\langle 1\rangle\subseteq y$, then
$g(x)\meet g(y)=f_1(x)$ and $g(x)\lex g(y)$.
\end{enumerate}

It follows that any two elements of $D$ are incomparable,
that is, $D$ is an antichain.  Hence $D$ is diagonal.
By construction, $D$ is a subset of $S$.
By the above three properties, 
$g$ preserves $\qle$, so $D$ is $\kappa$-dense.
\end{proof}

\section{Lower bound for dense linear orders}\label{sec.lower.bnd.1}
In this section, we show that all sparse vip $m$-types
can be embedded in any sparse diagonal set $D$ with $(D,\qle)$
$\kappa$-dense, and derive a lower bound result
for $\mathbb{Q}_\kappa$. 

\begin{definition}\label{def.compact}
Call an ordering $\prec$ of the levels of $\krtree$ \emph{\compact}
if  $(\pre\alpha2,\prec)$ has order type $2^{\alpha}$
  for each cardinal $\alpha<\kappa$.   
\end{definition}

\begin{lemma}\label{lem.compact.cofinal}
Suppose 
$\prec$ is a \compact ordering of the levels of $\krtree$
and $\alpha<\kappa$ is a cardinal.  For all $s$ with $\lg(s)<\alpha$,
the interval   $\set{t\in\pre\alpha2:s\subseteq t}$ is cofinal in
$(\pre\alpha2,\prec)$.   
\end{lemma}

\begin{proof} For all $s$ with $\lg(s)<\alpha$,
the interval $\set{t\in\pre\alpha2:s\subseteq t}$
has cardinality $2^{\alpha}$. Hence such an interval 
has order type $2^{\alpha}$ under $\prec$.
\end{proof}

\begin{lemma}\label{lem.compact.embed}
Suppose $\kappa$
is a limit cardinal, $\prec^{\,\prime}$ is a \compact ordering of the levels of $\krtree$
and $w\in\krtree$.
For all $n<\omega$ and orderings $\dless$ of the levels of 
$\prele n2$, there is an order preserving 
strong embedding $j$ of $(\prele n2,\dless)$
into $(\cone(w),\prec^{\,\prime})$, i.e. $s\dless t$ implies $j(s)\prec j(t)$.
Furthermore, $j$ may be chosen such that for all $s$, the length
of $j(s)$ is a cardinal.
\end{lemma}

\begin{proof}
Use induction on $n$.  For $n=0$, there is only the empty sequence,
and any embedding of this single point to a point $\cone(w)$
whose length is a cardinal works.  

Suppose the lemma is true for $m$ and $\dless$ is an ordering of
the levels of $\prele {n}2$ for $n=m+1$.  Let $j_0:\prele{m}2\to\krtree$ be
an order preserving strong embedding obtained from the induction
hypothesis.  We define an extension $j$ of $j_0$ as follows.
Let $\seq{s_i:i<2^{m+1}}$ enumerate in
  increasing $\dless$-order the nodes of $\pre{m+1}2$.
  Pick a cardinal $\alpha$ larger than the level at which
  $\pre{m}2$ is embedded, and define by recursion on $i$ nodes
  $j(s_i)=t_i$ in $\pre\alpha2$ such that $j(s_i\restrict
  m)\cat\sseq{s_i(m)}\subseteq t_i$ and $t_{i}\succ t_{i-1}$ if $i>0$.
Lemma \ref{lem.compact.cofinal} guarantees that this recursion
is possible.
\end{proof}

\begin{theorem}\label{thm:HL}
  Suppose that $\kappa$ is a cardinal which is measurable
  in the generic
  extension obtained by adding $\lambda$ Cohen subsets of $\kappa$,
  where $\lambda\to(\kappa)^{6}_{2^{\kappa}}$.  Further suppose
  $\prec$ is well ordering  of the levels of $\prel\kappa2$
and $\prec^{\,\prime}$ is a \compact ordering of the levels
of $\prel\kappa2$.  Then there is a strong embedding
$e$ and a node $w$ such that $e$ preserves $\prec^{\,\prime}$ and
$\prec$ and $\prec^{\,\prime}$ agree on all pairs in  $e[\cone(w)]$. 
\end{theorem}
\begin{proof}
Apply Shelah's Theorem \ref{thm.SHL} to the coloring  
$d:[\krtree]^2\to 2$ defined using the Boolean value operator $\|\cdot\|$
by 
\[
d(\{ s,t\})=\| s\lex t\iff s\prec t\| 
\] 
to get a strong embedding $e:\krtree\to\krtree$ and an element
$w\in \krtree$ such that for all $s,t\in\cone(w)$ with $s\prec^{\,\prime}t$
one has $e(s)\prec^{\,\prime}e(t)$ and 
the color of $d(\{ s,t\})$ 
depends only on the $\prec^{\,\prime}$-ordered similarity type
of the pair $\{ s,t\}$.

If $\lg(s)<\lg(t)$, then $\lg(e(s))<\lg(e(t))$,
so $\prec$ and $\prec^{\,\prime}$ agree on
$e(s),e(t)$, since both are level orders.
Similarly, if $\lg(s)>\lg(t)$, then $\prec$ and $\prec^{\,\prime}$ 
agree on $e(s),e(t)$.

Next consider pairs $\{s,t\}$ with
$\lg(s)=\lg(t)$, $s\lex t$ and $s\prec^{\,\prime}t$.
Note that since Cohen forcing on $\kappa$ does not
add bounded subsets to $\kappa$, our assumptions in particular
imply that $\kappa$ must be a limit cardinal.
Suppose $\alpha$ is a cardinal larger than $\lg(w)$ but $<\kappa$.  
Then $e[\pre\alpha2\cap\cone(w)]\subseteq
\pre\gamma2$ is infinite for some $\gamma<\kappa$.  
Let $\seq{s_n\in\pre\alpha2\cap\cone(w):n<\omega}$ 
be a sequence which is increasing in both the $\lex$ 
and $\prec^{\,\prime}$ orders.
Since $e$ is a strong embedding, the sequence
$\seq{e(s_n)\in\pre\gamma2:n<\omega}$ is $\lex$-increasing,
and cannot be decreasing in $\prec$.  Since $s\prec^{\,\prime} t$
implies $e(s)\prec^{\,\prime}e(t)$, it follows that 
 $\prec^{\,\prime}$ and $\prec$ agree
on pairs $\{e(s),e(t)\}$ with 
$\lg(s)=\lg(t)$, $s\lex t$ and $s\prec^{\,\prime}t$.

Finally consider pairs $\{s,t\}$ with
$\lg(s)=\lg(t)$, $t\lex s$ and $s\prec^{\,\prime}t$.
For $\alpha$ as in the previous case,
let $\seq{t_n\in\pre\alpha2\cap\cone(w):n<\omega}$ 
be a sequence which is decreasing in the $\lex$ order
and increasing in the $\prec^{\,\prime}$ order.
Then $\seq{e(s_n)\in\pre\gamma2:n<\omega}$ is $\lex$-decreasing,
and cannot be decreasing in $\prec$.  Thus by an argument
like that above,
$\prec^{\,\prime}$ and $\prec$ agree
on pairs $\{e(s),e(t)\}$ with 
$\lg(s)=\lg(t)$, $t\lex s$ and $s\prec^{\,\prime}t$.
\end{proof}

For the following definition, recall the notation $T(D)$ from Definition
\ref{def.T(S)}. 

\begin{definition}\label{def.semi-strong}
Suppose $D\subseteq\krtree$ is diagonal and $(D,\qle)$ is
$\kappa$-dense.  A function $f:\krtree\to T(D)$ is a
\emph{semi-strong embedding} if $f$ preserves extension and
lexicographic order, maps levels to levels, and for
every $s\in\krtree$, there is some $v\in D$ such that
\[
f(s)\subseteq f(s\cat\langle 0\rangle)\meet
f(s\cat\langle 1\rangle)\subsetneq v
\]
and $\lg(v)<\lg(f(s\cat\langle 0\rangle))$.
\end{definition}

\begin{lemma}\label{lem.semi-strong}
Suppose $\kappa$ is inaccessible and $D\subseteq\krtree$ is diagonal and $(D,\qle)$ is
$\kappa$-dense.  Then there is a semi-strong embedding
$f:\krtree\to T(D)$. 
\end{lemma}

\begin{proof}
Let $C\subseteq\kappa$ be the set of all 
limit $\alpha>0$ such that for all $t\in T(D)\cap\pre \alpha2$,
$t$ is a limit of densely splitting points of $T(D)$ 
and for all $\beta<\alpha$,
there is $v\in T(D)\cap\prele \alpha2$ such that 
$t\restrict\beta\subseteq v\in D$.
By Lemma \ref{lem.nearly.perfect}, $C$ is closed unbounded.

Define $f$ on $\pre \alpha2$ by recursion on $\alpha<\kappa$, using the
assumption that $\kappa$ is inaccessible.
To start the recursion, let $f(\emptyset)$ be an element of
$T(D)$ on the $\gamma_0$ level of $\krtree$, where $\gamma_0$
is the least infinite cardinal in $C$.

If $\alpha$ is a limit and $f$ has been defined on ${}^{\alpha>}2$,
then let $\gamma_\alpha$ be a cardinal in $C$ greater than
$\sup_{\beta<\alpha}\gamma_\beta$, and for each $s\in \pre \alpha2$,
let $f(s)$ be an element of $T(D)\cap\pre {\gamma_\alpha}2$
which is an evenhanded limit of densely splitting points of $T(D)$ extending
$\bigcup\set{f(s\restrict\beta):\beta<\alpha}$.

Suppose $\alpha=\alpha'+1$ is a successor and $f$ has been defined
on ${}^{\alpha>}2$.  Let $\gamma_\alpha$ be a cardinal in $C$
greater than $\gamma_{\alpha'}$.  Recall that by Lemma \ref{lem.min.W(S)}
every node of $T(D)$ has a minimal extension
to a densely splitting node (one in $\mathcal{W}(S)$).   Also, by Lemma \ref{lem.robust},
every evenhanded  extension of a node of $T(D)$ has an extension
to a densely splitting node.
For $s=t\cat\langle\delta\rangle$ and $u_t$ the minimal length
densely splitting node properly extending $t$, let $f(s)$ be an evenhanded
extension of $u_t\cat\langle\delta\rangle$ of length $\gamma_\alpha$.

Use induction to show that $f$ preserves extension and lexicographic
order.  By its construction and the choice of $C$, the remaining
conditions are satisfied for all $s$. 
\end{proof}

For the following definition, recall the notion of a sparse
$m$-type from Definition
\ref{def.sparse.type}. 

\begin{theorem}\label{thm:Qk}
  Suppose that $\kappa$ is a cardinal which is measurable
  in the generic
  extension obtained by adding $\lambda$ Cohen subsets of $\kappa$,
  where $\lambda\to(\kappa)^{6}_{2^{\kappa}}$.  If $D\subseteq\pre\kappa2$ 
  is a sparse diagonal set with $(D,\qle)$ $\kappa$-dense 
  and  $\prec$  is a $D$-vip order 
  of the levels of $\krtree$, 
  then 
  every sparse vip $m$-type $(\tau,\dless)$ 
  is realized as $(\clp(x),\prec_x)$ for some $x\subseteq D$.
\end{theorem}


\begin{proof}
The assumptions imply that $\kappa$ is inaccessible, so Lemma \ref{lem.semi-strong} applies.
Let $f:\krtree\to T(D)$ be a semi-strong embedding from Lemma
\ref{lem.semi-strong}.  Define $f^0:\krtree\to \closure{D}$
by $f^0(s):=f(s\cat\langle 0\rangle)\meet f(s\cat\langle 1\rangle)$.
Let $f^1:\krtree\to D$ be defined by $f^1(s)=v$ where $v$ is the
minimal extension of $f^0(s)$ in $D$ as guaranteed by Definition
\ref{def.semi-strong}. 

For $t\in \pre \alpha2$ and $i=0,1$, define well-orderings $\prec^i_t$
on $\pre\alpha2$ as follows: let $\beta_i=\lg(f^i(t))$ and
set $s\prec^i_t s'$ if and only if 
$f(s\cat\langle 0\rangle)\restrict\beta_i
\prec
f(s'\cat\langle 0\rangle)\restrict\beta_i$.

Let $\prec^{\,\prime}$ be any \compact well-ordering of the levels of $\krtree$.
Call a triple $\set{s,s',t}$ \emph{local} if 
$\lg(s)=\lg(s')=\lg(t)$, $s\lex s'$, $t\prec^{\,\prime} s$,
and $t\prec^{\,\prime} s'$, 
and $s\meet s'\not\subseteq t$.  
Let $d$ be a coloring of the triples of $\krtree$  defined as follows:
if $\set{s,s',t}$ is not local, let $d(\set{s,s',t}):=(2,2)$ and 
otherwise set
\[
d(\set{s,s',t}):=
(\|s\prec^{\,\prime} s'\iff s\prec^0_t s'\|,
\|s\prec^{\,\prime} s'\iff s\prec^1_t s'\|).
\] 

Apply 
Shelah's Theorem \ref{thm.SHL} to $d$ and $\prec^{\,\prime}$
to obtain a strong embedding $e:\krtree\to\krtree$ and
a node $w$ such that
for triples from $T:=e[\cone(w)]$, the coloring depends only
on the $\prec^{\,\prime}$-ordered similarity type of the triple.  Then
two local triples $\set{s,s',t}$ and $\set{u,u',v}$ of $T$
are colored the same if and only if
\begin{equation*}
s\prec_{t}^{0}s'\iff u\prec_{v}^{0}u'\qquad\text{and}\qquad
s\prec_{t}^{1}s'\iff u\prec_{v}^{1}u'.
\end{equation*}
 
Hence for $t\in T$ and for both ordered similarity types of incomparable
pairs of same length nodes from $T$, the orderings $\prec_t^{\delta}$ 
must always agree with one of  $\prec^{\,\prime}$ and its converse on $T$. 
Since $\prec$ is a well-ordering of the levels of $\krtree$,
the orderings $\prec_t^i$ are well-orderings of $\pre \alpha2$,
so they must agree with $\prec^{\,\prime}$ in both cases.  

Assume that $\tau$ has $n+1$ leaves and
let $L$ be the set of these leaves.
Then $\tau$ is a subtree of $\prele{2n}2$ and every level
of $\prele{2n}2$ has exactly one element of $\closure{L}$.
Extend $\dless$ defined on $\tau$ to 
$\dless^*$ defined on all
of $\prele{2n}2$ in such a way that the extension is still a 
$\closure{L}$-vip order.
 
Apply Lemma \ref{lem.compact.embed} 
to get an order preserving strong embedding $j$
of $(\prele{2n}2,\dless^*)$ into $(\krtree,\prec^{\,\prime})$.

Let $\seq{t_\ell:\ell\le 2n}$ enumerate the elements of $\closure{L}$
in increasing order of length.  Note that $\lg(t_\ell)=\ell$.
For $\ell\le 2n$, define
$\beta_\ell:=\lg(f^i(e(j(t_\ell))))$ where $i=0$ if $t_\ell\notin L$
and $i=1$ if $t_\ell\in L$.

Finally define $\sigma:\tau\to T(D)$ by recursion on $\ell\le 2n$.
For $\ell=0$, let $\sigma(\emptyset)=f^0(e(j(\emptyset)))$.
For $\ell>0$, consider three cases for elements of
$\tau\cap{}^\ell 2$.  If $t_\ell\in L$,
let $\sigma(t_\ell)=f^1(e(j(t_\ell)))$.  
If $t_\ell\notin L$, let $\sigma(t_\ell)=f^0(e(j(t_\ell)))$.  
Note that in both of these cases,  $\beta_\ell=\lg(\sigma(t_\ell))$. 
If $s\in\tau\setminus\closure{L}$ has length $\ell$, then there
is a unique immediate successor in $\tau$, $s\cat\langle 0\rangle$.
In this case, let
$\sigma(s)=f^0(e(j(s))\cat\langle 0\rangle)\restrict\beta_\ell$. 
Since $j$ sends $\dless^*$-increasing pairs to
$\prec^{\,\prime}$-increasing pairs
and $e$ is a $\prec^{\,\prime}$ order preserving strong embedding, 
their composition sends sends $\dless^*$-increasing pairs to
$\prec^{\,\prime}$-increasing pairs.  
Since for $v_\ell=e(j(t_\ell))\in T$, the order
$\prec^{\,\prime}$ agrees with $\prec^i_{v_\ell}$ on $T\cap \pre
\gamma2$ where $\gamma=\lg(v_\ell)$, it follows that $\sigma$
sends $\dless^*$-increasing pairs to $\prec$-increasing pairs.
Since $f$ preserves extension and
lexicographic order, $\sigma$ does as well.  By construction
$\sigma$ sends levels to levels, meets to meets (split nodes) and 
leaves to leaves (terminal nodes).   Let $x=\sigma[L]$ be the image under
$\sigma$ of the leaves of $\tau$.  Then $(\clp(x),\prec_x)=(\tau,\dless)$,
as required.
\end{proof}

\begin{theorem}\label{thm.lower.bound.1}
Let $m$ be a natural number and 
  suppose that $\kappa$ is a cardinal which is measurable
  in the generic
  extension obtained by adding $\lambda$ Cohen subsets of $\kappa$,
  where $\lambda\to(\kappa)^{2m}_{2^{\kappa}}$. 
Then for $r=t_m^+$ equal to the number of sparse vip $m$-types,
the $\kappa$-dense linear order $\mathbb{Q}_\kappa$ 
satisfies 
$$\mathbb{Q}_\kappa\nrightarrow(\mathbb{Q}_\kappa)^m_{<\omega,r-1}.$$
\end{theorem}

\begin{proof}
Let $S\subseteq \krtree$ be a cofinal transverse subset
obtain from Lemma \ref{lem.cofinal.transverse}.
Let $\prec$ be a pre-$S$-vip
order on $\krtree$ obtained from Lemma \ref{lem.exist.vip}. 

Let $\varphi:\krtree\to S$ be a $\qle$-preserving injection
such that $D:=\varphi[\krtree]$ is a sparse diagonal set
with $\closure{D}\subseteq S$.  Notice that $(D,\qle)$ is
$\kappa$-dense, since $\varphi$ is $\qle$-preserving.

Let $(\tau_0,\dless_0)$, \dots , $(\tau_{r-1},\dless_{r-1})$
be an enumeration of the sparse vip $m$-types.
Define $c:[\krtree]^m\to r$ by $c(a)=i$ where for $x:=\varphi[a]$,
$(\clp(x),\prec_x)=(\tau_i,\dless_i)$.

Suppose $A\subseteq \krtree$ is a subset with $(A,\qle)$
$\kappa$-dense and $i<r$.
Since $\varphi$ is $\qle$-preserving, its image
$B:=\varphi[A]\subseteq D$ is a sparse diagonal set with
the property that $(B,\qle)$ is $\kappa$-dense. Thus
by Theorem \ref{thm:Qk}, the sparse vip $m$-type $(\tau_i,\dless_i)$
is realized as $(\clp(x),\prec_x)$ for some $x\subseteq B$.
Since $\varphi$ is injective,  there is an $m$-element subset
$u\subseteq A$ with $\varphi[u]=x$ and $c(u)=i$.
Since $A$ and $i$ were arbitrary, in every $\kappa$-dense
subset $A\subseteq \krtree$, every color $i$ is
realized by some $m$-element subset. Therefore,  the theorem follows. 
\end{proof}

Recall the definition of canonical partition introduced
immediately after the statement of Theorem \ref{thm.partition}.  

\begin{theorem}\label{thm.canon.1}
Let $m$ be a natural number and 
  suppose that $\kappa$ is a cardinal which is measurable
  in the generic
  extension obtained by adding $\lambda$ Cohen subsets of $\kappa$,
  where $\lambda\to(\kappa)^{2m}_{2^{\kappa}}$. 
For $t_m^+$ equal to the number of sparse vip $m$-types,
there is a canonical partition of the $m$-element
subsets of $\mathbb{Q}_\kappa=([\krtree]^m,\qle)$
into $t_m^+$ parts.
\end{theorem}

\begin{proof}
Use Lemma \ref{lem.cofinal.transverse} to find
$S\subseteq \krtree$ cofinal and transverse.
Use Lemma \ref{lem.exist.vip} to find $\prec$ a
pre-$S$-vip order on $\krtree$.

For all $w\in\krtree$, let $\varphi_w$ be a sparse diagonalization
of $\krtree$ into $S\cap\cone(w)$.  Use recursion on $\prec$
to define $\pi:\krtree\to\krtree$ such that for all $t\in\krtree$,
$\pi(t)$ is an extension of $t$ with $\cone(\pi(t))$
disjoint from the union over all $s\prec t$ of
$\varphi_{\pi(s)}[\krtree]$.  Since the order type of
$\set{s\in\krtree:s\prec t}$ is less than $\kappa$ and
each $\varphi_{\pi(s)}[\krtree]$ is a sparse diagonal subset of $S$,
it is always possible to continue the recursion.

Define $h:\krtree\to\krtree$ as follows.  For $t$ with
$\otp \set{s\in\krtree:s\prec t}=\alpha$, $z_\alpha={}^\alpha\{0\}$,
and $u\in \set{z_\alpha}\cup\cone(z_\alpha\cat\langle 1\rangle)$,
let $h(u)=\varphi_{\pi(t)}(u)$.

Let $(\tau_0,\dless_0)$, $(\tau_1,\dless_1)$,
\dots , $(\tau_{r-1},\dless_{r-1})$ enumerate the
sparse vip $m$-types.  Let $C_0$ be the set of all $m$-element
subsets $A$ for which 
$(\clp(h[A]),\prec_{h[A]})$ is either $(\tau_0,\dless_0)$ or
not a sparse vip $m$-type.  For positive $j<r$, let 
$C_j$ be the set of all $m$-element subsets $A$ of $\krtree$ for which 
$(\clp(h[A]),\prec_{h[A]})=(\tau_j,\dless_j)$.
Then 
$\mathcal{C}:=\set{C_0,C_1,\dots,C_{r-1}}$ is a partition of $[\krtree]^m$
into $r$ sets.  

To see that each class of $\mathcal{C}$ is indivisible, 
suppose $d:C_j\to\mu$ is a fixed coloring for some $2\le \mu<\kappa$.  
Extend $d$ to all of $[\krtree]^m$ by setting $d(A)=0$ if
$j>0$ and $A\notin C_j$ or by setting $d(A)=1$ if $j=0$
and $A\notin C_j$.  Apply Shelah's
Theorem \ref{thm.SHL} to the restriction to antichains to obtain
a strong embedding $e$ and a node $w$ such that $e$ preserves
$\prec$ on $\cone(w)$ and $d$ is constant on $m$-element subsets
of the same $\prec$-ordered similarity type.
Let $\alpha$ be the order type of $\set{s\in\krtree: s\prec w}$.
Then $\cone(z_\alpha\cat\langle 1\rangle)$ is a $\kappa$-dense
subset.  Since $h$ agrees with $\varphi_{\pi(w)}$ on
$\cone(z_\alpha\cat\langle 1\rangle)$ and $\varphi_{\pi(w)}$
is $\qle$-preserving, $D:=h[\cone(z_\alpha\cat\langle 1\rangle)]$ 
is a $\kappa$-dense subset of $\cone(\pi(w))\subseteq\cone(w)$.
Since $\varphi_{\pi(w)}$ is a sparse diagonalization, the set
$D$ is a sparse diagonal set with $\closure{D}\subseteq
\closure{S}\cap\cone(\pi(w))$.  
Thus by Lemma \ref{lem.sparse.type1}, all $\prec$-ordered
similarity types of $m$-element subsets of $D$ are
sparse vip $m$-types.  Let $K:=e[D]$.  By Lemma \ref{lem.embed.sparse},
all $\prec$-ordered similarity types of $m$-element subsets
of $K$ are sparse vip $m$-types.  It follows that
$[K]^m\cap C_j$ is $d$-monochromatic.  Hence each $C_j$
is indivisible.

To see that each class of $\mathcal{C}$ is persistent, 
suppose $K\subseteq\krtree$ is $\kappa$-dense and $j<r$.
By Lemma \ref{lem.split.continue}, there is a node $z^*$
in $\mathcal{W}(K)$, the set of densely splitting nodes of $K$.
Thus $z^*\cat\langle 1\rangle$ has a $\kappa$-dense
set of extensions in $K$.  In other words, $K\cap\cone(z^*)$
is $\kappa$-dense.  Let $z_\alpha$ be the longest
initial segment of $z^*\cat\langle 1\rangle$ consisting only of zeros.
Then $\cone(z^*)\subseteq\cone(z_\alpha\cat\langle 1\rangle)$.
Let $t$ be the $\alpha$th element of $\krtree$ in the $\prec$
order.  Since $h$ agrees with $\varphi_{\pi(t)}$
on $\cone(z_\alpha\cat\langle 1\rangle)$, it follows
that $h[K\cap\cone(z^*)]$ is a sparse diagonal
set with the property that $(h[K\cap\cone(z^*)],\qle)$
is $\kappa$-dense.  Thus
by Theorem \ref{thm:Qk}, the sparse vip $m$-type $(\tau_j,\dless_j)$
is realized as $(\clp(x),\prec_x)$ for some $x\subseteq h[K\cap\cone(z^*)]$.
Since $h$ is injective,  there is an $m$-element subset
$u\subseteq K\cap\cone(z^*)$ with $h[u]=x$, and $u\in C_j$
as required.
\end{proof}

\section{Sparse vip types}\label{sec.sparse.vip.types}

In this section we give closed form upper and lower
bounds for the number of sparse vip $m$-types that
facilitate comparisons with D. Devlin's theorem
for $(\mathbb{Q},<)$.  We then describe a recursive
procedure for computing the number of  sparse vip $m$-types.

In Figure~\ref{fig.tau*}, we give a picture 
of a specific example of a sparse vip $5$-type, 
which we will call $(\tau^*,<)$, 
so we can use it in later examples to illustrate
a variety of definitions.  To translate the figure
into a representation in which each node is a
sequence of $0$'s and $1$'s, note that the root is the
empty sequence and only
line segments with positive slope represent $1$'s.
We have circled the
nodes in the meet closure of the set $L$ of leaves of $\tau^*$; these
nodes are the designated elements for any ordering
of the levels which makes $\tau^*$ a sparse vip $5$-type.
For only three pairs of nodes
does the requirement that $(\tau^*,<)$
be a sparse vip $5$-type fail to specify the order, and
between each such pair we have indicated the order. 

\setlength{\unitlength}{5mm}
\begin{figure}[htb]
\begin{center}
\begin{picture}(9,9)(0,1)
\thinlines
\put(4,1){\line(2,1){2}}
\put(4,1){\line(-2,1){2}}
\put(2,2){\line(1,1){1}}
\put(2,2){\line(-1,1){1}}
\put(6,3){\line(1,1){1}}
\put(6,3){\line(-1,1){1}}
\put(7,6){\line(1,1){1}}
\put(7,6){\line(-1,1){1}}
\put(6,2){\line(0,1){1}}
\put(1,3){\line(0,1){2}}
\put(3,3){\line(0,1){1}}
\put(5,4){\line(0,1){4}}
\put(7,4){\line(0,1){2}}
\put(8,7){\line(0,1){2}}
\put(4,1){\makebox(0,0){$\bullet$}}
\put(4,1){\makebox(0,0){$\bigcirc$}}
\put(2,2){\makebox(0,0){$\bullet$}}
\put(2,2){\makebox(0,0){$\bigcirc$}}
\put(6,2){\makebox(0,0){$\bullet$}}
\put(1,3){\makebox(0,0){$\bullet$}}
\put(3,3){\makebox(0,0){$\bullet$}}
\put(6,3){\makebox(0,0){$\bullet$}}
\put(6,3){\makebox(0,0){$\bigcirc$}}
\put(1,4){\makebox(0,0){$\bullet$}}
\put(3,4){\makebox(0,0){$\bullet$}}
\put(3,4){\makebox(0,0){$\bigcirc$}}
\put(5,4){\makebox(0,0){$\bullet$}}
\put(7,4){\makebox(0,0){$\bullet$}}
\put(1,5){\makebox(0,0){$\bullet$}}
\put(1,5){\makebox(0,0){$\bigcirc$}}
\put(5,5){\makebox(0,0){$\bullet$}}
\put(7,5){\makebox(0,0){$\bullet$}}
\put(5,6){\makebox(0,0){$\bullet$}}
\put(7,6){\makebox(0,0){$\bullet$}}
\put(7,6){\makebox(0,0){$\bigcirc$}}
\put(5,7){\makebox(0,0){$\bullet$}}
\put(6,7){\makebox(0,0){$\bullet$}}
\put(6,7){\makebox(0,0){$\bigcirc$}}
\put(8,7){\makebox(0,0){$\bullet$}}
\put(5,8){\makebox(0,0){$\bullet$}}
\put(5,8){\makebox(0,0){$\bigcirc$}}
\put(8,8){\makebox(0,0){$\bullet$}}
\put(8,9){\makebox(0,0){$\bullet$}}
\put(8,9){\makebox(0,0){$\bigcirc$}}
\put(2,3){\makebox(0,0){$<$}}
\put(6,4){\makebox(0,0){$<$}}
\put(6,5){\makebox(0,0){$>$}}
\end{picture}\\
\end{center}
\caption{Pictorial representation of $(\tau^*,<)$}\label{fig.tau*}
\end{figure}

\begin{lemma}\label{lem.char.m.type}
If $\tau$ is a sparse $m$-type and $L$ is the set of its leaves, then 
$\tau=\set{x\restrict i:x\in L\land i<2m-1}$.
\end{lemma}

\begin{proof}
By definition, $\tau$ is closed under initial segments.
\end{proof}

Before introducing a lemma on properties of the
leaves of a sparse $m$-type, 
we list the leaves of $\tau^*$ (see Figure~\ref{fig.tau*} on 
page~\pageref{fig.tau*}) with the lexicographically least
one at the top of the stack, and continuing in increasing order down
the stack.  
\begin{quote}
$\langle
0,0,0,0
\rangle$\newline
$\langle
0,1,0
\rangle$\newline
$\langle
1,0,0,0,0,0,0
\rangle$\newline
$\langle
1,0,1,0,0,0
\rangle$\newline
$\langle
1,0,1,0,0,1,0,0
\rangle$\newline
\end{quote}

\begin{lemma}\label{lem.small.prop}
Suppose $\tau$ is a sparse $m$-type whose set of
leaves is $D=\set{d_0,d_1,\dots,d_{m-1}}$ listed in
$\lex$-increasing order.
Then
\begin{enumerate}
\item $\closure{D}=D\cup\set{d_i\meet d_{i+1}:i<m-1}$;
\item $d_0$ is a sequence of all zeros;
\item if $i<m-1$, then $d_i(\lg(d_i\meet d_{i+1}))=0$ and
$d_{i+1}(\lg(d_i\meet d_{i+1}))=1$;
\item if $i<m-1$, then for all $p$ with $\lg(d_i\meet d_{i+1})<p<\lg(d_{i+1})$,
$d_{i+1}(p)=0$.
\end{enumerate}
\end{lemma}

\begin{proof}
The first item follows from the fact that $D\cup\set{d_i\meet d_{i+1}:
i<m-1}$ is a subset of $\closure{D}$ of size $m+(m-1)=2m-1=|\closure{D}|$.
Since $D$ is diagonal, all its elements have
different lengths from among $0,1,\dots , 2m-2$.
Since $\closure{D}$ has $2m-1$ elements, all of different lengths,
it follows that
for all $p<2m-1$, there is some $x=x_p\in \closure{D}$ with $\lg(x)=p$.

Now the second item holds, since
for all $p<\lg(d_0)$, $d_0(p)=d_0(\lg(x_p))=0$,
either because $x_p=d_i\meet d_{i+1}\subseteq d_0$
and $d_0\lex d_{i+1}$ or because $\tau$ is sparse.

Then the third item follows from the definition
of lexicographic order in a binary tree.

For the fourth item, fix attention on some $i<m-1$.
Notice that if $\ell<i$ and $d_\ell\meet d_{\ell+1}\subseteq d_{i+1}$,
then 
$(d_\ell\meet d_{\ell+1})\cat\langle 1\rangle\subseteq d_i\meet d_{i+1}$,
since $d_{\ell+1}\lex d_i$.
Also, if $i\le \ell$ and $d_\ell\meet d_{\ell+1}\subseteq d_{i}$,
then $(d_\ell\meet d_{\ell+1})\cat\langle 0\rangle\subseteq d_i$,
since $d_i\lex d_{\ell+1}$.  Now the fourth item follows from the
previous two statements and the fact that $\tau$ is sparse.
\end{proof} 

Next we associate with each sparse $m$-type 
a sequence which is characteristic.

\begin{definition}\label{def.alternate}
Suppose $\tau$ is a sparse $m$-type whose set of
leaves is $D=\set{d_0,d_1,\dots,d_{m-1}}$ listed in
$\lex$-increasing order.
Define $\altperm(\tau): (2m-1)\to (2m-1)$ by
\[
\altperm(\tau)(k)=\begin{cases}
\lg(d_i),&\text{if $k=2i$,}\\
\lg(d_i\meet d_{i+1}),&\text{otherwise.}\\
\end{cases}
\]
\end{definition}

\begin{lemma}\label{lem.pi.one-one}
If $\tau$ and $\tau'$ are  sparse $m$-types
and $\altperm(\tau)=\altperm(\tau')$,
then $\tau=\tau'$.
\end{lemma}

\begin{proof}
Let $D=L(\tau)$ and $E=L(\tau')$ be the sets of leaves 
of $\tau$ and $\tau'$.  List $D$ and $E$ in increasing lexicographic order
as $d_0$, $d_1$, \ldots , $d_{m-1}$ 
and 
$e_0$, $e_1$, \ldots , $e_{m-1}$.  Then for all $i<m$,
$\lg(d_i)=\lg(e_i)$ by definition of $\altperm$.
Also, for $i<m-1$, $\lg(d_i\meet d_{i+1})
=\lg(e_i\meet e_{i+1})$.  Use induction on
$i<m$ and the previous lemma to show $d_i=e_i$.
\end{proof}

\begin{definition}\label{def.perm}
A function $\altperm:(2m-1)\to (2m-1)$ is an
\emph{alternating permutation}
if it is a permutation, for all
even $i<2m-2$, $\altperm(i)>\altperm(i+1)$ and for all odd $i<2m-2$,
$\altperm(i)<\altperm(i+1)$.
\end{definition} 
 
The study of alternating permutations and the alternating group dates back to
Andr\'{e} \cite{A1881} in 1881, and continues
to be an active area of investigation (see \cite{MSY96}
for example).

\begin{lemma}\label{lem.pi.corr}
For all positive integers $m$,
$\altperm$ is a bijection between the collection of sparse $m$-types
and the the set of alternating permutations on $(2m-1)$.
\end{lemma}

\begin{proof}
Properties of meets guarantee that $\altperm(\tau)$ is an alternating
permutation whenever $\tau$ is 
a sparse $m$-type.  By the previous
lemma, $\altperm$ is one-to-one.  

\newcommand{\pau}{p}
To see that it is onto, suppose
$\pau:(2m-1)\to(2m-1)$ is an alternating permutation.
Let $d_0$ be the sequence of zeros of length $\pau(0)$.
Continue the recursion by defining $d_1$, $d_2$, \ldots ,
$d_{m-1}$, by setting $d_{i+1}$ to be the sequence of
of length $\pau(2i+2)$ which extends $(d_i\restrict \pau(2i+1))
\cat\langle 1\rangle$ with all zeros.  Since $\pau$ is
a permutation, all elements of $D:=\set{d_0,\dots,d_{m-1}}$
have different lengths, and $D$ is an antichain listed
in $\lex$-increasing order.  Also by construction,
$\lg(d_j\meet d_{j+1})=\pau(2j+1)$, so $D$ is diagonal.
By construction, $d_0(p)=0$ for all $p<\lg(d_0)$.
Use induction on $i$ to show that for all $i<m-1$ and
all $p<\lg(d_{i+1})$, if $d_{i+1}(p)=1$,
then for some $\ell\le i$, $p=\lg(d_\ell\meet d_{\ell+1})$
and $d_\ell\meet d_{\ell+1}\subseteq d_{i+1}$.
It follows that $D$ is sparse.  Let $\tau=\clp(D)$.  Then
the set of leaves of $\tau$ is $D$ and $\tau$ is a sparse $m$-type.
By construction $\altperm(\tau)=\pau$.  Since $\pau$ was
arbitrary, the mapping $\altperm$ is onto.
\end{proof}

The above proof was inspired by the counting of
\emph{Joyce trees} (named by Ross Street \cite{Street02})
in a paper by Joyce \cite{Joycepreprint} applying category theory 
to physics.  The meet closure of an $m$-element
diagonal subset of ${}^{(2m-1)>}2$ is an example
of a Joyce tree.  They are counted up to
a certain equivalence, and each equivalence class
has exactly one example whose closure under initial segment
is a sparse $m$-type.  We also used ideas from a
counting argument by Vuksanovic
in \cite{V02pams}.

\begin{lemma}\label{lem.size.alt}
The number of alternating permutations on $2m-1$ is  
$t_m$, where $t_m$ is the $m$-th tangent number, defined
recursively by $t_1=1$ and 
\[
t_n=\sum_{i=1}^{n-1}\binom{2n-2}{2i-1}t_it_{n-i}.
\]
\end{lemma}

\begin{proof}
An exponential generating function for the sequence $a_n$, 
where $a_n$ is the number of alternating
permutations on $n$, is $\sec(x)+\tan(x)$
(see Stanley's
\emph{Enumerative Combinatorics I} \cite{Stanley1}).
Since the
terms corresponding to $\sec(x)$ in this series 
have even powers and
the terms corresponding to $\tan(x)$ have odd powers,
an exponential generating function for $b_m$, where 
$b_m$ is the number of alternating permutations on $(2m-1)$ 
is $\tan(x)$.

It is not difficult to show directly that the number of 
alternating permutations on $2m-1$ satisfies
the above recurrence.   For notational convenience, let $K$
denote the set $2m-1$.  Given a $2i-1$ element subset $I\subseteq
K\setminus\set{2i-1}$ and alternating permutations $p_0$
on $2i-1$ and $p_1$ on $2(m-i)-1=2m-2i-1$, 
use the unique order preserving maps $e_0$ from $(2i-1)$ to $I$ and 
$e_1$ from $(2m-2i-1)$ to $K\setminus (I\cup\set{2i-1})$ to 
$p=p_0\circ e_0^{-1} \cup \set{ (2i-1,0)}\cup p_1\cup e_1^{-1}$.
(See \cite{Stanley1} for details). 
\end{proof}

\begin{cor}\label{cor.size.alt}
The set of sparse $m$-types has cardinality 
$t_m$, where $t_m$ is the $m$-th tangent number,
and $t_m\le t_m^+$ where $t_m^+$ is the number of sparse vip $m$-types.
\end{cor}

The next goal is to provide a closed form upper bound for
the number of sparse vip $m$-types as a product of $t_m$
times a constant factor.  We start by looking at the
size of levels of $\clp(A)$ for $A$ an $m$-element diagonal set.

\begin{definition}\label{def.row.size}
Suppose $A\subseteq{}^{(2m-1)>}2$ is a diagonal set.
Enumerate $\closure{A}$ in increasing order of length
as $a_0,a_1,\dots,a_{2m-2}$ and
define 
\[
\ell_i(A):=i+1-2\left|
\set{j<i:a_j\in \closure{A}\setminus {A}}\right|.
\]
\end{definition}

Note that for $m$-element diagonal subsets of ${}^{(2m-1)>}2$,
the lengths of the elements of the meet closure are
$|a_0|=0,|a_1|=1,\dots,|a_{2m-2}|=2m-2$.

\begin{lemma}\label{lem.level.size}
For $m\ge 2$, and any $m$-element diagonal subset 
$A\subseteq {}^{(2m-1)>}2$
whose meet closure is listed in increasing order as
$a_0,a_1,\dots,a_{2m-2}$, the cardinality of level $i$
of $\clp(A)$ is
$$ |\set{a_j\restrict{i}:i<2m-1}|=\ell_i(A)=
(2m-1)-i-2|\set{j\ge i:a_j\in\closure{A}\setminus A}|.$$
Moreover, if $i< m$, then
$\ell_i(A)\le i+1$ and $\ell_{2m-2-i}\le i+1$.
\end{lemma}

\begin{proof}
Use induction to compute the size of 
$L_i:=\set{a_j\restrict{|a_i|}:j<2m-1}$ in two different ways.
To start the inductions, note that 
for $i=0$ and $i=2m-2$, $L_0=\{a_0\}$ and $L_{2m-2}=\{a_{2m-2}\}$
are singletons, so $\ell_0(A)=1$ and $\ell_{2m-2}(A)=1$.
Also, since $a_0$ is meet of two elements and $a_0,a_1$ are the shortest 
elements, we have $\ell_1(A)=|L_1|=2$.  Since $a_{2m-2},a_{2m-3}$
are the longest elements, we have 
$L_{2m-3}=\{a_{2m-2}\restrict_{{2m-3}},a_{2m-3}\}$ is also
a doubleton, thus $\ell_{2m-3}(A)=2$.
Counting up from $i=0$, we get the values in the formula
for $\ell_i(A)$, since the value for $|L_{i+1}|$ is one more
than the value for $|L_i|$ if $a_i\in \closure{A}\setminus A$ and 
one less if $a_i\in A$.  Thus for $0<i<m$, $\ell_{i}(A)\le \ell_{i-1}(A)+1\le
i+1$.
Counting down from $i=2m-2$, we get the values
in the displayed formula, since the value for $|L_{i-1}|$ is one more
than the value for $|L_i|$ if $a_{i-1}\in A$ and 
one less if $a_{i-1}\in \closure{A}\setminus A$.
Thus for $0<i<m$, $\ell_{2m-2-i}(A)\le \ell_{2m-2-(i-1)}+1\le i+1$.
%
\end{proof}

The values 
$\ell_i(A)= i+1$ and $\ell_{2m-2-i}= i+1$ are achieved if
$A$ is a comb whose leaves listed in increasing order
as $a_0,a_1,\dots,a_{m-1}$ satisfy $\lg(a_i)=m-1+i$
and for $i<m-1$, $\lg(a_i\meet a_{i+1})=i$.

\begin{definition}\label{def.V(A)}
For any $m$-element diagonal subset 
$A\subseteq {}^{(2m-1)>}2$,  
let $V(A)$ be the
number of pairs $(\clp(A),\prec)$ where $\prec$ is a
$\closure{A}$-vip order of $\clp(A)$.
\end{definition}

\begin{lemma}\label{lem.V(A).bnd}
For any $m$-element diagonal subset 
$A\subseteq {}^{(2m-1)>}2$, 
the number of $A$-vip orders of
$\clp(A)$ is bounded above: 
\[
V(A)\le \prod_{0<i<2m-2}(\ell_{i}(A)-1)!\le (m-1)!\prod_{i<m-1}(i!)^2.
\] 
\end{lemma}

\begin{proof}
Recall that every level of $\clp(A)$ has an element
of $\closure{A}$.  For the first inequality, note
that in every $A$-vip order, the element of $\closure{A}$
on each level is the least element of the level. 
The righthand side of the  first inequality is 
the count of level orders that satisfy this constraint.
For the second inequality, use the estimates of
Lemma \ref{lem.level.size} and the fact that $0!=1$.
\end{proof}

\begin{lemma}\label{lem.up.set}
The number of sparse vip $m$-types
at bounded above by $t_m^+\le t_m(m-1)!\prod_{i<m}(i!)^2$.
\end{lemma}

\begin{proof}
By Lemma \ref{cor.size.alt}
the number of sparse $m$-types is $t_m$, so the
lemma follows from Lemma \ref{lem.V(A).bnd}.
\end{proof}

Using notions from finite combinatorics, such as reverse Raney
sequences, and a fine analysis of sparse $m$-types, it is
possible to calculate a closed form lower bound for the values of
$t_m^+$. This is done in an upcoming paper \cite{Jeanfincomb}
by J. Larson, one of the conclusions of which is

\begin{theorem}\label{thm.up.sparse}[Larson \cite{Jeanfincomb}] For all $m\ge 2$, $t_m+(2^{m-1})(-1+
\Pi_{i<m}i!)\le t_m^+$.
\end{theorem}

Figure \ref{fig.Rvalue} summarizes 
the  calculation from \cite{Jeanfincomb} of values of
$t_m^+$ for $m\le 5$.  A comparison with $t_m$ is also included.

\begin{figure}[t]
\[
\begin{array}{|r|r|r|}
\hline
m&t_m^+&t_m\\
\hline
1&1&1\\
\hline
2&2&2\\
\hline
3&20&16\\
\hline
4&776&272\\
\hline
5&151,184&7936\\
\hline
\end{array}
\]
\caption{Some small values of $t_m^+$ and $t_m$.}\label{fig.Rvalue}
\end{figure}
%

\section{A proof of Shelah's Theorem}\label{sec.Shproof}

In this section we proof of Theorem \ref{thm.SHL} based on Shelah's
proof  \cite{Sh288} together with ideas from \cite{Sh289}.
The major part of the proof is dedicated to the proof of
Lemma~\ref{thm:endhom} below; following this we show that this lemma
implies Theorem~\ref{thm.SHL}.

\begin{lemma}[End Homogeneity]
  \label{thm:endhom}
  Assume $m\ge 2$ and 
  that $\kappa$ is measurable in the model obtained by adding
  $\lambda$ Cohen subsets of $\kappa$, where 
  $\lambda\rightarrow(\kappa)^{2m}_{2^\kappa}$.   Then for any  
  well
  ordering $\prec$ of the levels of $\krtree$ and 
  coloring  $d\colon \bigcup_{\alpha<\kappa}[{^\alpha}2]^m\to\sigma$
  of the $m$-element level sets of $\krtree$ 
  with $\sigma<\kappa$ colors
  there is a 
  strong embedding $e\colon\krtree\cong T\subseteq\krtree$
  such that 
  whenever $s\deq (s_0,\dots, s_{m-1})\in [{^\alpha}2]^m$ and
  $\beta<\alpha$ are such that the members of 
  $s\restrict\beta\deq(s_0\restrict\beta,\dots,s_{m-1}\restrict\beta)$
  are distinct, and $s$ and $s\restrict\beta$ are ordered the same way
  by $\prec$, then we have $d(e[s])=d(e[s\restrict\beta])$.
\end{lemma}
We actually use a slightly stronger version of Lemma~\ref{thm:endhom},
although we will indicate how this can be avoided at the cost of a
slight strengthening of the hypothesis:
\begin{lemma}\label{thm:endhomx}
  Suppose that $\kappa$ and $m$ are as in Lemma~\ref{thm:endhom}, and
  that for each $\xi<\kappa$ we have a coloring $d_\xi$ of the
  $m$-element level sets of 
  $\krtree$ in fewer than $\kappa$ colors.   Then there is a strong
  embedding $e$ such that whenever $s$ and $\beta$ are as in
  Lemma~\ref{thm:endhom} we have $d_\xi(e[s])=d_\xi(e[s\restrict\beta])$
  for each $\xi<\kappa$
\end{lemma}
We will give the proof of Lemma~\ref{thm:endhom}, and will indicate in
footnotes how this should be modified to prove
Lemma~\ref{thm:endhomx}.

 

\begin{proof} Let $m,\kappa, d$ and $\prec$ be as in
  Lemma~\ref{thm:endhom}, and let $P$ be the the forcing notion 
  adding $\lambda$ many  Cohen subsets of $\kappa$.  Thus a condition in
  $P$ is a function $p$ with domain 
  $\supp(p)\in[\lambda]^{<\kappa}$ and with values 
  $p(\nu)\in\prel\kappa2$.
  We abuse notation by writing $p'\supseteq p$ if $p'$ is stronger than
  $p$, that is, if $\supp(p')\supseteq\supp(p)$ and $p'(\nu)\supseteq p(\nu)$
  for each $\nu\in\supp(p)$.

  For $i<\lambda$ let $\name\eta_i$ be a name for the $i$th Cohen
  subset of $\kappa$, and let $\name\uf$ be a name for the
  $\kappa$-complete ultrafilter 
  $\uf$ on $\kappa$ assumed to 
  exist in the generic extension.    For any $u\in[\lambda]^{<\omega}$
  there is a unique $j<\sigma$ such that 
  $\set{\xi<\kappa:d(\set{\eta_i\restrict\xi:i\in u})=j}\in\uf$; we write
  $d(u)$ for this $j$ and let $\name d(u)$ be a name for
  $d(u)$.\footnote{For Lemma~\ref{thm:endhomx}, this becomes $\name
  d_{\xi}(u)$ for each $\xi<\kappa$.}

  If $W\subseteq\lambda$ then we write $P\cut W$ for $\set{p\in
    P:\supp(p)\subseteq W}$.   If $H\colon W\to W'$ is an order preserving
  map between subsets of $\lambda$ then $h_H\colon P\cut W\to P\cut W'$ is
  the map defined by $h_{H}(p)(H(i))=p(i)$, and if 
  $\otp(W)=\otp(W')$ then we define 
  $h_{W,W'}=h_H$ where $H\colon W\cong W'$ is
  the unique order preserving map.

  \begin{claim}
    \label{thm:Wu}
    There is a set $Z\in[\lambda]^\kappa$ and a function
    $W\colon[Z]^{\le m}\to[\lambda]^{\le\kappa}$ satisfying the
    following conditions:
    \begin{enumerate}
    \item \label{item:w1}If $u\in[Z]^{\le m}$ then $u\subseteq W(u)$, and
      $P\cut W(u)$ contains a maximal antichain of
      conditions $p$ 
      deciding the value of $\name d(u)$.\footnote{For
      Lemma~\ref{thm:endhomx}, $P\cut W(u)$ contains a maximal
      antichain of conditions deciding the value of $\name d_{\xi}(u)$
      for each $\xi<\kappa$.} 
    \item\label{item:w2} If $u,u'\in[Z]^{\le m}$ and $|u|=|u'|$, then
      $\otp(W(u))=\otp(W(u'))$, the function
      $h_{W(u),W(u')}$ maps $u$ to $u'$, 
      and for all $p\in 
      P\cut W(u)$ and $j<\sigma$, 
      \begin{equation}\label{eq:1}
        p\forces \name d(u)=j\iff h_{W(u),W(u')}(p)\forces \name d(u')=j.
      \end{equation}
    \item\label{item:w3} If $u'\subseteq u\in[Z]^{\le m}$ then
      $W(u')\subseteq W(u)$, 
      and if $u,u'\in[Z]^{\le m}$ then $W(u\cap u')=W(u)\cap W(u')$.
    \end{enumerate}
  \end{claim}
  \begin{proof}
    Because $P$ has the $\kappa^+$-chain condition, there are sets
    $W'(u)\in [\lambda]^{\le\kappa}$ such  that clause~1 is satisfied
    when $W'$ is substituted for $W$.
    We  can also arrange that $W'(u')\supseteq W'(u)$ whenever
    $u'\supseteq u$.

    Now define an equivalence relation on $[\lambda]^{\le 2m}$ as
    follows: two sets 
    $u$ and $u'$ in $[\lambda]^{\le 2m}$ are equivalent if they
    satisfy the following two conditions: 
    \begin{enumerate}\renewcommand\theenumi{\alph{enumi}}
    \item $|u|=|u'|$, and if $|u|\le m$ then clause~2 above holds with
      $W'(u)$ and $W'(u')$ substituted for $W(u)$ and $W(u')$. 
    \item \label{item:b}
      Suppose that $u_1,u_2\subseteq u$, with $u_1,u_2\in[\lambda]^{\le m}$,
      and that $u'_1$ and $u'_2$ are the
      subsets  of $u'$ such that
      $(u,u_1,u_2,<)\cong(u',u'_1,u'_2,<)$.  Then this isomorphism
      extends to an isomorphism
      \[(W'(u),W'(u_1),W'(u_2),<)\cong(W'(u'),W'(u'_1),W'(u'_2),<).\]
    \end{enumerate}

    There are only $2^{\kappa}$ equivalence classes, and therefore
    the
    assumption that $\lambda\to(\kappa)^{2m}_{2^\kappa}$ implies that
    there is 
    $Z'\in[\lambda]^\kappa$ such that any pairs $u,u'\in[Z']^{k}$
    are equivalent for all $k\le 2m$.   It follows that clauses~1 and~2
    are satisfied  when $W'$ and
    $Z'$ are substituted for $W$ and $Z$. 

    Now define 
    \begin{equation}\label{eq:S}
      W(u)\deq\bigcup\set{\bigcap_{v\in X}W'(v):
        X\subseteq[Z']^{\le m}\And \bigcap X\subseteq u},
    \end{equation}
    and let $Z=\set{\gamma_{\omega\nu}:\nu<\kappa}$ where
    $\seq{\gamma_\nu:\nu<\kappa}$ is the increasing enumeration of
    $Z'$.   We claim that $W$ and $Z$ are as required.

    First note that we only need to consider finite sets $X$ in
    \eqref{eq:S}.  To see this, suppose $\bigcap X\subseteq u$, 
    fix a set $u_0\in X$, and pick $X'\subseteq X$ such that $u_0\in
    X'$ and for each $\xi\in u_0\setminus u$ there is a set $u_{\xi}\in X'$
    such that $\xi\notin u_{\xi}$.  Then $\bigcap X'\subseteq u$, $|X'|\le
    1+{|u_0|}-| u|\le 1+m$, and $\bigcap_{v\in X'}W'(v)\supseteq
    \bigcap_{v\in X}W'(v)$. 

    Now note that condition~(\ref{item:b}) can be 
    extended to arbitrary finite sets: If $\set{u_i:i<k}$ and
    $\set{u'_i:i<k}$ are subsets of $[Z']^{<\omega}$ such that
    \begin{equation}
      \label{eq:u}
      \left(\bigcup_{i<k}u_i,u_0,\dots,u_{k-1},<\right)\cong
      \left(\bigcup_{i<k}u'_i,u'_0,\dots,u'_{k-1},<\right)  
    \end{equation}
    then
    \begin{multline}\label{eq:4}
      \left(\bigcup_{i<k}W'(u_i),W'(u_0),\dots,W'(u_{k-1}),<\right)\\\cong
      \left(\bigcup_{i<k}W'(u'_i),W'(u'_0),\dots,W'(u'_{k-1}),<\right).
    \end{multline}
    To see this, note that~condition~(\ref{item:b}) implies that
    $\otp(W'_{u})=\otp(W'_{u'})$ whenever $u,u'\in[Z']^{s}$ for some
    $s\le2m$.  Let $W'(u_i)=\set{\gamma_{i,\nu}:\nu<\xi_{i}}$ and
    $W'(u'_i)=\set{\gamma'_{i,\nu}:\nu<\xi_i}$ be the increasing
    enumerations.   
    It follows that 
    $\bigcup_{i<k}W'(u_i)=\set{\gamma_{i,\nu}:i<k\And\nu<\xi_i}$; and
    for each $i,i'<k$ and $\nu<\xi_i$ and $\nu'<\xi_{i'}$ we have
    $\gamma_{i,\nu}=\gamma_{i',\nu'}\iff\gamma'_{i,\nu}=\gamma'_{i',\nu'}$,
    and
    $\gamma_{i,\nu}<\gamma_{i',\nu'}\iff
    \gamma'_{i,\nu}<\gamma'_{i',\nu'}$.  This easily
    implies~\eqref{eq:4}.

    A consequence of this is that the intersection $\bigcap_{v\in
    X}W'(v)$ in~\eqref{eq:S} does not depend on $X$, but only on the
    isomorphism type of $(\bigcup X,u,u_0,\dots,u_{k-1})$ where
    $X=\set{u_i:i<k}$.   The proof of this proceeds by induction on the size
    of the set of ordinals in $\bigcup X$ on which two sets $X$ and
    $X'$ with the same isomorphism type differ.  Let
    $\xi\in\bigcup X$
    and $\xi'\in\bigcup X'$ be the least corresponding pair of
    ordinals which differ between $X$ and $X'$.  Assuming without loss
    of generality that $\xi<\xi'$, let $X''$ be the
    set obtained by replacing 
    $\xi'$ by $\xi$ in $X'$.  Then $\bigcap_{v\in
    X}W'(v)=\bigcap_{v\in X''}W'(v)$ by the induction hypothesis.   To
    see that $\bigcap_{v\in X'}W'(v)=\bigcap_{v\in X''}W'(v)$, pick
    some $v_0\in X'$ with $\xi\notin v_0$.  Then
    $\bigcap_{v\in X'}W'(v)\subseteq W'(v_0)$, and the members of
    $W'(v_0)$ are not moved in the isomorphism~\eqref{eq:4} between
    $\bigcup_{v\in X'}W'(v)$ and $\bigcup_{v\in X''}W'(v)$.

    This implies that the union in~\eqref{eq:S} can be taken to be
    finite, since there are only finitely many such isomorphism
    classes, and hence $|W(u)|\le\kappa$.
    If $u,u'\in[Z]^{s}$ for some $s\le m$ then it follows that
    $\otp(W(u))=\otp(W(u'))$, since the fact that $Z$ contains only
    limit members of $Z'$ implies that for each set $X$ contributing an
    intersection to $W(u)$ there is an isomorphic set $X'$ contributing an
    intersection to $W(u')$.   Furthermore this isomorphism between
    $W(u)$ and $W(u')$ preserves the isomorphism between $W'(u)$ and
    $W'(u')$, so the equivalence~\eqref{eq:1} in
    clause~\ref{item:w2} of Claim~\ref{thm:Wu} holds for $W$ as well
    as for $W'$.   Finally, 
    the definition of $W(u)$ easily implies clause~\ref{item:w3} of
   Claim~\ref{thm:Wu}. 
  \end{proof}

  We are now ready to construct the promised strong embedding $e$. 
  For sets $u\in[Z]^{\le m}$ define
  $W^*(u)=W(u)\setminus\bigcup_{u'\subsetneq u}W(u')$.  
  
  For $\alpha<\kappa$ let $R_{\alpha}$ be the
  set of  one to one functions $s$ with $\dom(s)\in[\pre\alpha2]^{\le
  m}$ and $\ran(s)\subseteq Z$.   If $s\in R_\alpha$,  $\beta<\alpha$,
    and the 
  members of $\set{x\restrict\beta:x\in\dom(s)}$ are all distinct then
  we will abuse notation by writing $s\downto\beta$ for the function
  $s'\in R_{\beta}$ with $\dom(s')=\set{x\restrict\beta:x\in \dom(s)}$
    defined by 
  $s'(x\restrict\beta)=s(x)$.
  The sets $R_{\alpha}$ include the empty function $\emptyset$ as a
  member, and we will abuse the notation by taking the function
    $\emptyset\in R_{\alpha}$ to be different from $\emptyset\in
    R_{\beta}$ whenever $\alpha\not=\beta$.

  We use recursion on $\alpha<\kappa$ to define an ordinal
  $\zeta_\alpha<\kappa$ and a map
  $e\restrict\pre\alpha2\colon\pre\alpha2\to\pre{\zeta_\alpha}2$,
  along with conditions $p_s\in P|W^*(\ran(s))$ for each $s\in R_{\alpha}$. 
  Several times in the following construction we will use the
  observation that this, together with the fact from
  Claim~\ref{thm:Wu} that $W(u)\cap W(u')=W(u\cap u')$ for all
  $u,u'\in [Z]^{\le m}$,
  implies that $p_s$ and $p_{s'}$ are compatible whenever $s\cup s'$
  is a function.   Hence $\bigcup\set{p_{s}:s\in R_{\alpha}\And
  s\subseteq\tau}$ is a condition for any one to one function
  $\tau$ with $\dom(\tau)\subseteq\pre\alpha2$ and
  $\ran(\tau)\subset Z$.     In particular this is true for all
  $\tau\in R_{\alpha}$.  

  The construction will satisfy the following induction hypotheses:
  \begin{enumerate}
  \item\label{item:psorder}
    If
    $g\colon\ran(s)\to Z$ preserves order, then $p_{g\circ
      s}=h_{W^*(\ran(s)),W^*(\ran(g))}(p_s)$. 
  \item \label{item:0}
    If $\beta<\alpha$ and the sets
    $\sing{x\restrict\beta:x\in\dom(s)}$ are distinct then
    $p_s\supseteq p_{s\downto\beta}$.   
  \item \label{item:1}
    $\bigcup_{t\subseteq s}p_t\forces \name d(\ran(s))=
    d(e[\dom(s)])$ whenever $s$ is an order isomorphism between
    $(\dom(s),\prec)$ and $(\ran(s),\in)$. \footnote{For
    Lemma~\ref{thm:endhomx} this should hold for $d_{\xi}$ for all
    $\xi<\alpha$.} 
  \item\label{item:2}   
    If $\dom(s)=\sing x$ and $s(x)=\eta\in Z$ then $p_{s}(\eta)=e(x)$.
  \end{enumerate}

  Clause~\ref{item:psorder} asserts that $p_s$ depends, up to
  isomorphism, only on $\dom(s)$ and the order which $s$ induces on
  that set.  
  
  Clause~\ref{item:1} will imply the required end-homogeneity, since
  if $\beta<\alpha$ then 
  $\bigcup_{t\subseteq s}p_t$ forces that $\name d(\ran(s))=\name
  d(\ran(s\downto\beta))$, and consequently $d(e[\dom(s)])=d(e[\dom(
  s\restrict\beta)])$ whenever $\prec$ orders $\dom(s)$ and
  $\dom(s\restrict\beta)$ alike.

  For $n$ small enough that $2^{n}<m$ we can define $\zeta_n=n$, along
  with 
  $e(s)=s$ and $p_s=\emptyset$ for all $s\in R_{n}$.   The
  construction for $\alpha$
  with $2^{\alpha} \ge m$ consists of three steps.  The first step
  defines conditions $\bar p_s$ which will satisfy
  clauses~\ref{item:psorder} and~\ref{item:0} (with $\bar p_s$ instead
  of $p_s$) and 
  satisfies clause~\ref{item:1} to the extent that 
  \begin{equation}
    \label{eq:2}
    (\forall s\in R_{\alpha} )   (\exists j_s<\sigma)\;\bigcup_{t\subseteq
    s}\bar 
    p_{t}\forces \name 
    d(\ran(s))=j_s. 
  \end{equation}
  The second step will define
  $\zeta_\alpha$ and define
  $e\restrict\pre\alpha2\colon\pre\alpha2\to\pre{\zeta_\alpha}2$ so
  that if the ordering given to $\dom(s)$ by $s$  agrees with the
  given ordering 
  $\prec$ then $d(e[\dom s])$ has the value $j_s$ determined in the
  first step.
  The final step consists of setting $p_s=\bar p_s$ except for the
  adjustments necessary to satisfy clause~\ref{item:2}. 

  The first step is divided into two cases, depending on whether or not
  $\alpha$ is a limit ordinal.

  \case{1}{$\alpha$ is a limit ordinal}
  Let $\bar\zeta:=\sup_{\beta<\alpha}
  \zeta_\beta$. For $x\in {}^\alpha 2$ let
  $\bar e(x):=\bigcup_{\beta<\alpha}(e(x\rest\beta))$, and for
  $s\in R_{\alpha}$, let 
  \begin{equation*}
    \bar p_s:= \bigcup
    \set{p_{s\downto\gamma}:
      \beta<\alpha\And\sing{x\restrict\beta:x\in\dom(s)} 
      \text{ are pairwise distinct} }.
  \end{equation*}
  Note that the  induction hypothesis implies that  $\bar p_s\in P\cut
  W^*(\ran(s))$.  Furthermore $\bar p_s$ satisfies 
  clause~\ref{item:0}, and by the induction hypothesis this implies
  that $\bar p_s$ satisfies~\eqref{eq:2}. 

  \case{2}{$\alpha$ is a successor ordinal}
  Let $\alpha=\beta+1$,  set $\bar\zeta=\zeta_\beta+1$, and set $\bar
  e(x)=e(x\restrict\beta)\cat\sseq{x(\beta)}\in\pre{\bar\zeta}2$ for
  each $x\in\pre\alpha2$. 
  For each $s\in R_{\alpha}$ such that the members of
  $\set{x\restrict\gamma:x\in\dom(s)}$ are  
  distinct set 
  $p^{0}_{s}=p_{s\restrict\beta}$.    If the members of
  $\set{x\restrict\gamma:x\in\dom(s)}$ are not distinct then set
  $p^{0}_s=\emptyset$.  

  Let $\seq{s_i:i<\gamma}$ enumerate the set of
  $s\in R_{\alpha}$ such that $\ran(s)$ is an initial segment of $Z$
  of length $m$.   
  Thus for every $s\in R_\alpha$ there is a unique ordinal $i<\gamma$ and order
  preserving map $g$ such that $s=g\circ s_i$.
  Now define $p^{i}_s$
  for each    $i\le\gamma$ and  
  $s\in R_{\alpha}$ by  recursion on $i\le\gamma$:
  If $i$ is a limit ordinal then $p^{i}_s=\bigcup_{i'<i}p^{i'}_s$.
  For a successor ordinal $i+1$, suppose that  $p^{i}_s$ is defined for
  all $s\in R_{\alpha}$ and set  
  $q'\deq\bigcup_{t\subseteq s_i}p^{i}_{t}$.
  Now choose $q\supseteq q'$ in $W(\ran(s_i))$ so that $q$ decides the
  value\footnote{For
    Lemma~\ref{thm:endhomx},   $q$ decides $\name
    d_{\xi}(\ran(s_i))$ for each $\xi<\alpha$.} of $\name
  d(\ran(s_i))$, and for each $t\subseteq s_i$ set 
  $p^{i+1}_{t}=q\restrict W^*(\ran(t))$.
  To finish up, define $p^{i+1}_{g\circ
  t}=h_{W^*(\ran(s_i)),W^*(\ran(g))}(p_t)$, as in
  clause~\ref{item:psorder}, for all $t\subseteq s_i$ and all order
  preserving functions $g\colon\ran(t)\to Z$; and set $p^{i+1}_s=p^{i}_s$ for
  all other $s\in R_{\alpha}$.

  Now complete the first step of
  the construction by setting $\bar p_s=p^{\gamma}_s$.

  \smallskip

  For the second step of the construction let
  $\tau\colon\pre{\alpha}2\to Z$ be a one to one map which 
  preserves the given ordering $\prec$ on $\pre\alpha2$.   Such a map
  exists since $|\pre\alpha2|<\kappa=|Z|$.
  Now set
  $q\deq\bigcup\set{p_{\tau\restrict y}:y\in[\pre\alpha2]^{\le m}}$.
  This is a condition, and $q$ decides the value of 
  $\name d({\tau[y]})$ for each $y\in[\pre\alpha2]^{m}$. 
  
  Now, since $\name\uf$ is forced to be a $\kappa$-complete ultrafilter in the
  generic extension, there is a condition $q'\supseteq q$ and an ordinal
  $\xi\ge\bar\zeta$ such that 
  \begin{equation}
    \label{eq:3}
    q'\forces
    d(\set{\eta_\nu\restrict\xi:\nu\in\tau[y]})=\name
    d(\tau[y]))\footnote{For Lemma~\ref{thm:endhomx} this should hold
    for $d_{\xi}$ for each $\xi<\kappa$.} 
  \end{equation}
  for each $y\in[\pre\alpha2]^{m}$.   We can assume without loss of
  generality that for each  $x\in\pre{\alpha}2$, we have
  $\dom(q'(\tau(x)))\supseteq\xi$.   Set $\zeta_\alpha=\xi$ and
  complete the definition of
  $e\restrict\pre\alpha2$ by setting 
  $e(x)=q'(\tau(x))\restrict\xi$.   Finally complete the definition of
  $p_s$ by setting $p_s=\bar p_{s}$ 
  unless $s=\tau\restrict\sing{x}$ for some $x$, in which case define
  $p_{s}\supseteq \bar 
  p_s$ by setting $p_s(\tau(x))=q'(\tau(x))\restrict\xi+1$. 

  This completes the definition of the embedding $e$, and hence of
  the proof of the End Homogeneity Lemma~\ref{thm:endhom}.
\end{proof}

\begin{cor}\label{thm:presE}
  Suppose that $\prec$ and $\set{\prec_\xi:\xi<\sigma}$ are well orderings of 
  $\krtree$.  Then there is a strong embedding 
  $e\colon\krtree\to\krtree$ 
  such that for a dense set of nodes $t\in\krtree$, we have
  $s_0\prec s_1\iff e(s_0)\prec_{\xi}e(s_1)$ for all $s_0,s_1\supseteq
  t$ and $\xi<\sigma$.
\end{cor}
\begin{proof}
  Define $d\colon\bigcup_{\alpha<\kappa}[\pre\alpha2]^2\to\pre\sigma2$
   by setting $d(s_0,s_1)(\xi)$ equal to the boolean value
   $\bool{s_0\lex s_1\iff s_0\prec_{\xi}s_1}$. 
  
  Apply Lemma~\ref{thm:endhom} to get a strong embedding
  $e\colon\krtree\to\krtree$ so that if $s_0\lex s_1$ are in
  $\pre\alpha2$ then $d(e(s_0),e(s_1))$ depends only on whether or not
  $s_0\prec s_1$.     We want to find, for any given node $t$,  a node
  $t'\supset t$ such that
  $d(e(s_0),e(s_1))(\xi)=\bool{s_0\prec s_1}$ for all
  $s_0,s_1\in\pre\alpha2$ for $\alpha>\lg(t)$ such that $s_0\lex s_1$,
  $t'\subseteq s_0\land s_1$.
  It will be sufficient to show that we can do this for any one
  $\xi<\sigma$ and for one of the two cases  $s_0\prec s_1$ or
  $s_1\prec s_0$.   We will 
  do it for $s_0\prec s_1$, and will indicate the change for $s_1\prec
  s_0$. 

  Suppose, for the sake of contradiction, that for every $t'\supset t$
  there is a level set $\sing{s_0,s_1}$ such that $s_0\lex s_1$,
  $t'\subseteq s_0\land s_1$ and  $s_0\prec s_1$ but
  $e(s_0)\succ_{\xi}e(s_1)$.    By the end-homogeneity the same will be
  true of any $s'_0\prec s'_1$ such that $s'_i\supseteq s_i$.

  Define an infinite sequence of pairs  $\seq{(u_n,v_n):n<\omega}$ of
  nodes as follows:  Set $v_0=t$.
  If $v_n$ is defined then let $u_n,v_{n+1}$ be nodes of the same
  length, extending $v_n$, 
  such that $u_{n}\lex v_{n+1}$ and 
  $e(s_0)\succ_{\xi}e(s_1)$ for all $s_0\supseteq u_{n}$ and
  $s_1\supseteq v_{n+1}$ such that $s_0\prec s_1$.

  Now fix $\alpha<\kappa$ such that $\alpha>\lg(u_n)$ for each
  $n<\omega$, and choose $w_n\supseteq u_n$ for each $n<\omega$.
  The nodes $w_n$ have the properties that $w_n\lex w_{n'}$ for each 
  $n<n'$, and $e(w_n)\succ_{\xi}e(w_{n'})$ for each $n<n'$ such that
  $w_n\prec w_{n'}$.  Now apply Ramsey's theorem to get
  an infinite subset $X\subseteq\omega$ such that the boolean values
  $\bool{w_n\prec w_{n'}}$ are
  constant for $\sing{n,n'}\in[X]^{2}$.    Since $\prec$ is a well
  order we must have $w_{n}\prec w_{n'}$ for all $n<n'$ in $X$, but
  then $\seq{e(w_{n}):n\in X}$ is an infinite descending
  $\prec_{\xi}$-sequence, contradicting the assumption that
  $\prec_{\xi}$ is a well order.

  The proof for the case $s_1\prec s_0$ is the same except that
  the pair $\sing{u_{n},v_{n+1}}$ satisfies $v_{n+1}\lex u_{n}$ 
  for $s_0\supseteq u_n$ and $s_1\supseteq v_{n+1}$. 

\end{proof}

\begin{cor}
  Suppose that $\kappa$, $m$, $\prec$, and $d_\xi$ are as in
  Lemma~\ref{thm:endhomx}.   Then there is a strongly embedded tree
  $T\subseteq\krtree$ of height $\kappa$ such that
  $d_\xi(s)=d_\xi(s')$ whenever
  $s$ and $s'$ are $\prec$-similar members of $[\pre\alpha2\cap
  T]^m$. and $\xi$ is smaller than the number of split levels of $T$ below
  $\alpha$.  
\end{cor}
\begin{proof}
  Let $e\colon\krtree\to\krtree$ be the strong embedding
  given by Lemma~\ref{thm:endhomx}.   By Corollary~\ref{thm:presE}
  there is a strong embedding $e'\colon\krtree\to\krtree$ and a
  $t\in\krtree$ such that,
  for all $s,s'\supseteq t$ we have $e'(s)\prec e'(s')\iff e'(s)
  \mathrel{e'[\prec]} e'(s')\iff s \mathrel{e^{-1}[\prec]} s'$.
  Then $T=e'e[\krtree]$ is as required.
\end{proof}

\begin{lemma}\label{thm:SHLvar}
  Suppose that $\kappa$ is a cardinal which is measurable
  in the generic
  extension obtained by adding $\lambda$ Cohen subsets of $\kappa$,
  where $\lambda\to(\kappa)^{2m}_{2^{\kappa}}$. 
  Then for any coloring $d$ of the $m$-element antichains of $\krtree$
  into $\sigma<\kappa$ colors, 
  and any well-ordering $\prec$ of the levels of $\krtree$, 
  there is a strong embedding $e\colon\krtree\cong
  T\subseteq\krtree$ such 
  that 
  $d(a)=d(b)$ for all $\prec$-similar $m$-element antichains $a$ and
  $b$ of $T$. 
\end{lemma}
First we show that this lemma implies Theorem~\ref{thm.SHL}:
\begin{proof}[Proof of Theorem~\ref{thm.SHL} from Lemma~\ref{thm:SHLvar}]
  Let $e\colon\krtree\to\krtree$ satisfy the conclusion of
  Lemma~\ref{thm:SHLvar}.   By Corollary~\ref{thm:presE} there is a
  strong embedding $e'\colon\krtree\to\krtree$ and an dense set of
  nodes $w\in\krtree$ such that $s\prec t$ if and only if $e(s)\prec
  e(t)$ for all $s,t\in\bigcup_{\alpha<\kappa}\cap\cone(w)$.   The composition $e'\circ e$ satisfies the
  conclusion of Theorem~\ref{thm.SHL}.
\end{proof}

\begin{proof}[Proof of Lemma~\ref{thm:SHLvar}]
\newcommand{\prect}{\prec_t}
  The maximum height of a similarity tree with $m$ terminal
  nodes is equal to the number of meets plus the number of terminal
  nodes, which is $(m-1)+m=2m-1$.   
  We will define a sequence $T_0\subseteq T_1\subseteq\dots \subseteq
  T_{2m-2}\subseteq T_{2m-1}=\krtree$ of strongly embedded subtrees,
  using a reverse recursion on $k<2m-1$ starting
  with $T_{2m-1}=\krtree$. 

  For $k<2m-1$
  we assume as a recursion hypothesis that $T_{k+1}$ has the following 
  homogeneity property:

  \begin{quote}
    Suppose that $x,y\in[T_{k+1}]^{m}$ are $\prec$-similar antichains, with
    the common collapse $\clp(x,\prec)=\clp(y,\prec)=(t,\prect)$, and
    suppose that $i\restrict\prele
    k2=j\restrict\prele k2$ where 
    $i\colon(t,\prect)\cong (\closure{x}, \prec)$ and
    $j\colon(t,\prect)\cong(\closure{y},\prec)$ are the maps witnessing
    this similarity.    Then $d(x)=d(y)$.
  \end{quote}

  Note that this puts no constraint on $T_{2m-1}$, and it implies
  that $T_0$ satisfies the conclusion of Lemma~\ref{thm:SHLvar}.

  For each $\alpha<\kappa$, each antichain $a\in[T_{k+1}\cap\prel\alpha2]^{<m}$, and
  each ordered similarity tree $(t,\prect)$ consider the following
  coloring 
  $d_{k,a,t}$ of the ${\le}m$-element level sets
  $x\in[T_{k+1}\cap\pre\alpha2]^{\le m}$:   If $x\cup a$ is an 
  antichain, and $\clp(x\cup a)=(t\cap\prele i2,\prect)$, then
  $d_{k,a,t}(x)=d(y)$ where $y\in[T_{k+1}]^{m}$ is any antichain  such
  that 
  $x=\set{\nu\restrict\alpha:\nu\in
  y}$. 
  If $x\cup a$ is not an antichain,  if $\clp(x\cup
  a)\not=(t\cap\prele i2,\prect)$, 
  or if no such antichain $y$ exists then 
  set $d_{i,a,t}(x)=0$. 
  Lemma~\ref{thm:SHLvar} implies that there is  a strongly embedded
  subtree $T_{k}\subset T_{k+1}$ on which 
  $d_{k,a,t}(x)=d_{k,a,t}(\set{\nu\cap\beta:\nu\in x})$ whenever
  the nodes $\nu\cap\beta$ for $\nu\in x$ are distinct, and the sets 
  $x$ and $\set{\nu\cap\beta:\nu\in x}$ are both in $T_{k}$ and are
  ordered the same way by $\prec$.  Thus $T_{k}$ satisfies the
  recursion hypothesis. 
\end{proof}

We note that the proof of
Theorem~\ref{thm.SHL} does not require that the ultrafilter $\uf$ be
normal, and hence is valid for $\kappa=\omega$. This is essentially
Harrington's proof of the Halpern-L\"{a}uchli
theorem, which may be found in the Farah-Todorcevic book \cite{FaTo}.
Harrington's proof served as a starting point for Shelah in 
\cite{Sh288}.

\section{Further remarks on partition theorems}\label{sec.remarks}
There are a number of remaining open questions suggested by the
results presented so far so we comment on some of them. 

\begin{question}
 Is Theorem~\ref{thm.SHL} (or Lemma~\ref{thm:endhom}) consistent with
 the GCH? Is the conclusion of Theorem~\ref{thm.SHL} compatible with $L$,
 namely can there be an uncountable cardinal $\kappa$ in $L$ which satisfies
 that conclusion?
\end{question}
Note that the hypothesis of lemma~\ref{thm:endhom} implies that the
GCH fails at almost every $\alpha<\kappa$. Indeed for almost every
$\alpha<\kappa$  the power set of $\alpha$ in $V$ is a generic
extension obtained by adding $\lambda_{\alpha}$ Cohen subsets of
$\alpha$ to some ground model $M_\alpha$, where
$M_\alpha\models\lambda_a\rightarrow(\alpha)^{2m}_{(2^{\alpha})^{M_\alpha}}$,
since this is true in the generic extension of $V$ obtained by adding
$\lambda$ Cohen subsets to $\kappa$, and this extension does not add
bounded subsets to $\kappa$. 

The latter part of the question was asked by Michael Hru\v sak. See below
for an argument showing that such a $\kappa$ must be weakly compact.

\begin{question}
  Can Theorem~\ref{thm.SHL} be strengthened to require that
  $e[\prec]\subseteq{\prec}$?
\end{question}
It would be sufficient to prove this for level sets in
$[\krtree]^{2}$, as this would imply the corresponding modification of
Corollary~\ref{thm:presE}.   A positive solution would be likely to
also give a positive answer to the next question:

\begin{question}
  Is the relation $\krtree\xrightarrow{\text{str
      embedding}}(\krtree)^{<\omega}_{\sigma,{<}\omega}$ consistent?
  That is,  is it consistent that for any
  $d\colon[\krtree]^{<\omega}\to\sigma$ for some $\sigma<\kappa$, there is a
  strongly embedded tree $T$ such that $d[[T]^{m}]$ is finite
  for   each $m<\omega$?  
\end{question}
Shelah observes that the corresponding variation of
Lemma~\ref{thm:endhom} follows from the assumption that $\kappa$ is
measurable after adding $\lambda$ Cohen subsets of $\kappa$, where
$\lambda\rightarrow(\kappa)^{<\omega}_{2^{\kappa}}$.   This may be
easily seen by examining the proof of Lemma~\ref{thm:endhom} given
here, which requires only minor changes.
However the top-down proof of
Theorem~\ref{thm.SHL} from Lemma~\ref{thm:endhom} given here
does not work with 
the superscript ${{<}\omega}$.

It should be noted that in spite of the use of the ultrafilter $\uf$
in the construction, the set of splitting levels of the homogeneous
strongly embedded subtree $T$ is not in any sense a member of $\uf$;
indeed it is far from being even stationary.
This follows from the fact that the set of splitting levels of $T$
need not contain any fixed points, even for $m=1$, which in turn
follows from the observation that $\prel\alpha2$ contains only $\alpha$ 
many strongly embedded subtrees, each of which has  $2^{\alpha}$ many
branches.
Thus we can color, for each cardinal $\alpha$, the members of
$\pre\alpha2$ with $2^{\alpha}$ colors in such a way that  any strongly
embedded subtree of $\prel\alpha2$ contains branches with every
color. 

\begin{question}
  What is the large cardinal strength of the conclusion of
  Theorem~\ref{thm.SHL}? 
\end{question}

It is easy to see that a supercompact cardinal whose supercompactness
was made indestructible to $(<\kappa)$-directed closed
forcing using Laver's method \cite{Laver}, satisfies the assumptions of
the theorem. A much better upper bound is available. Namely,
work of Gitik in \cite{GiSh}
(mentioned also as \cite{Gitik} in \cite{Sh288}),
together the Erd\H{o}s-Rado theorem, implies that a 
model satisfying the hypothesis of 
Theorem~\ref{thm.SHL} can be constructed by forcing over a model of
$\text{GCH}+o(\kappa)=\kappa^{+2m+2}$.

On the other hand, it is easy to see
that Theorem~\ref{thm.SHL} implies that $\kappa$ is weakly compact.
To show that $\kappa\rightarrow
(\kappa)^2_2$, let the function $g:[\kappa]^2\to 2$ be given and
define a coloring  
$h$ of the two element antichains in $\krtree$ by  
by $h(\set{s,t}) = g(\set{\lg(s),\lg(s\meet t)})$ 
if $\lg(s)=\lg(t)$ and $h(\set{s,t})=0$ otherwise.
If $T$ is the strongly embedded tree whose existence is
guaranteed by Shelah's Theorem \ref{thm.SHL} then 
there are two possible $\prec$-similarity classes of $2$-element
level sets, 
one in which the lexicographic order of the pair agrees
with the $\prec$-order and the other in which these
two orders disagree.
By a Sierpinski argument, at the $\omega$th level of $T$,
there are pairs of both classes with meets on the same level.  
Since the coloring of these level sets depends only on
the lengths of the pair and their meet, both 
$\prec$-similarity classes receive the same color.
It follows that $g$ is monochromatic on pairs from 
the set of $\alpha$ for which $\alpha$ is a splitting level of $T$. 

As a consequence of Theorem~\ref{thm.partition} and an earlier result of
Hajnal and Komj\'{a}th (\cite{HaKo}) we obtain the following
theorem which shows that 
the
conclusion of Theorem~\ref{thm.SHL} does not follow from any large
cardinal hypothesis. This suggests that the use of an ultrafilter in a
generic extension, as opposed to one in $V$, is a necessary part of the
proof.

\begin{theorem}\label{thm.largecard} The
conclusion of Theorem~\ref{thm.SHL} does not follow from any large
cardinal hypothesis on $\kappa$.
\end{theorem}

\begin{proof}
Hajnal and Komj\'{a}th \cite{HaKo} show that there is a forcing
of size $\aleph_1$ which  adds an order type $\theta$ of size
$\aleph_1$ with the property that 
$\psi\not\rightarrow[\theta]^2_{\omega_1}$ for every order type
$\psi$, regardless of its size.
That is, 
there is a coloring of the
pairs of $\psi$ into $\aleph_1$ many colors
such that every suborder of type $\theta$ gets all the colors.
The conclusion of our Theorem~\ref{thm.partition} states, in contrast,
that 
the ordering $\QQ_{\kappa}$
has, for any coloring of its pairs, a subset of the full order type
$\QQ_{\kappa}$ which gets only 
finitely many colors.   This subset contains subsets of every order
type of size less than $\kappa$, since $\QQ_{\kappa}$ is a $(<\kappa)$-universal
linear order. Hence Theorem~\ref{thm.partition} cannot hold in the 
Hajnal-Komj\'{a}th extension, and so Shelah's theorem from which
it is derived, cannot either.
\end{proof}

\section{Upper bound for $\kappa$-Rado graphs}\label{sec.upper.bnd.2}
In this section we prove a limitation of colors
result for $\kappa$-Rado graphs orders using
Shelah's Theorem \ref{thm.SHL}. 
By a $\kappa$-Rado graph we mean a graph
$G$ of size $\kappa$ with the property that for every two disjoint subsets
$A$, $B$ of $G$, each of size $<\kappa$, there is $c\in G$ connected to
all points of $A$ and no point of $B$. The existence of such $G$ follows
from the assumption $\kappa^{<\kappa}=\kappa$ is regular.

In order to apply Shelah's Theorem, we need a method
of embedding $\kappa$-Rado graphs into $\krtree$
and a  well-ordering of the levels of $\krtree$
that is compatible with that embedding.  
We observe that any $\kappa$-Rado graph
is isomorphic to one whose universe is $\kappa$, and
generalize the approach used by
Erd\H{o}s, Hajnal and P\'osa \cite{EHP75} to embed
such a graph into $\krtree$.

\begin{definition}\label{def.translate}
Given a $\kappa$-Rado graph $\mathbb{G}=(\kappa,E)$,
the \emph{tree embedding of $\mathbb{G}$ into $\krtree$} is the function
$\sigma_{\mathbb{G}}:\kappa\to\krtree$ defined by
$\sigma_{\mathbb{G}}(0)=\emptyset$, and for $\alpha>0$,
$\sigma_{\mathbb{G}}(\alpha):\alpha\to 2$ is defined by   
$\sigma_{\mathbb{G}}(\alpha)(\beta)=1$ if and only if
$\set{\alpha,\beta}\in E$.
\end{definition}

\begin{lemma}\label{lem.sigma.cofinal}
For any $\kappa$-Rado graph $\mathbb{G}=(\kappa,E)$,
the range of the tree embedding $\sigma$ of $\mathbb{G}$
into $\krtree$ is a cofinal transverse subset of $\krtree$.
\end{lemma}

\begin{proof}
By definition of $\sigma$, for all $\alpha<\kappa$,
$\lg(\sigma(\alpha))=\alpha$, so 
the range of $\sigma$ is transverse.

To see that the range is cofinal, suppose $s\in\krtree$.
Let $\alpha=\lg(s)$ and let
$A$ be the set of all $\beta<\lg(s)$ with $s(\beta)=1$.
Since $\mathbb{G}$ is a $\kappa$-Rado graph, 
there is an element $\gamma>\alpha$ such that 
$\set{\beta,\gamma}\in E$ for all $\beta\in A$
and $\set{\beta,\gamma}\notin E$ for all $\beta\in\alpha\setminus A$.
It follows that $s\subseteq\sigma(\gamma)$.  Thus $\sigma[\kappa]$
is cofinal in $\krtree$.
\end{proof}

Our next goal in this section is a translation
of questions about isomorphisms of 
$\kappa$-Rado graphs 
into themselves to questions about $\krtree$.  Toward that end,
we define \emph{passing number preserving maps}.  This notion
was used in the proof of the limitation of colors result
by Laflamme, Sauer and Vuksanovic \cite{LSVpreprint}
for the countable Rado graph, which is also known as the
(infinite) random graph. 

\begin{definition}\label{def.pnp}
For $s,t\in\krtree$
with $|t|>|s|$, call $t(|s|)$ the \emph{passing number}
of $t$ at $s$. 
Call a function $f:\krtree\to \krtree$  \emph{passing number preserving}
or a \emph{pnp map} if it preserves 
\begin{enumerate}
\item length order: $\lg(s)<\lg(t)$
implies $\lg(f(s))<\lg(f(t))$; 
and 
\item passing numbers: $\lg(s)<\lg(t)$
implies $f(t)(\lg(f(s)))=t(\lg(s))$.
\end{enumerate}
\end{definition}

The first lemma states that any induced subgraph of a
$\kappa$-Rado graph $\mathbb{G}=(\kappa,E)$
has an induced subgraph which is 
isomorphic to the $\kappa$-Rado graph by an isomorphism
that also preserves $<$.

\begin{lemma}\label{lem.Rado.increasing}
For any cardinal $\kappa$ with $\kappa^{<\kappa}=\kappa$, any
$\kappa$-Rado graph $\mathbb{G}=(\kappa,E)$
and any $H\subseteq\kappa$ with $\mathbb{G}\cong (H,E\restrict H)$
there is a $<$-increasing map $g:\kappa\to H$ with
$\mathbb{G}\cong(g[\kappa],E\restrict g[\kappa])$.
\end{lemma}

\begin{proof}
Fix attention on a particular $\kappa$-Rado graph $\mathbb{G}=
(\kappa,E)$ and a specific induced subgraph $(H,E\restrict H)$
isomorphic to $\mathbb{G}$.  Let $h:\kappa\to H$ be the isomorphism.

Since $\kappa=\kappa^{<\kappa}$, by the mapping extension property,
for any $\gamma<\kappa$ and any
subset $A\subseteq \gamma$, there are cofinally many $\zeta$
with $\set{\delta<\gamma: \{\delta,\zeta\}\in E}=A$.

Define $z:\kappa\to\kappa$ and 
$g:\kappa\to H$ by recursion.  Let $z(\emptyset)=\emptyset$
and $g(\emptyset)=h(\emptyset)$.
Suppose $z\restrict\alpha$ and 
$g\restrict\alpha$ have been defined such that $z$ is increasing,
for all $\beta<\alpha$, $g(\beta)=h(z(\beta))$, and  
$g\restrict\alpha$ is an increasing
isomorphism of $(\alpha,E\restrict\alpha)$ into $(H,E\restrict H)$. 
Let $\gamma>\alpha$ be so large that if $h(\eta)<
\sup\set{\lg(g(\beta))+1:\beta<\alpha}$, then
$\gamma>\eta$.  
Let $A_\alpha:=\set{z(\beta): \beta<\alpha\land \{\beta,\alpha\}\in E}$.
Let $z(\alpha)\ge\gamma$ be such that 
$\set{\delta<\gamma: \{\delta,\zeta\}\in E}=A_\alpha$.
Let $g(\alpha)=h(z(\alpha))$.  Since $h$ is an isomorphism,
a pair $\{ h(z(\beta)),h(z(\alpha))\}$ is in $E$ if
and only if the pair $\{ z(\beta),z(\alpha)\}$ is in $E$.
It follows that $\set{\beta<\alpha: \{
  g(\beta),g(\alpha)\}\in E}=A_\alpha$.
Therefore by induction, $g$ is the desired increasing isomorphism
into $H$.
\end{proof}

\begin{lemma}\label{lem.Rado.iso.1}
Suppose $\mathbb{G}=(\kappa,E)$ is a $\kappa$-Rado graph 
with tree embedding $\sigma$ and $S=\sigma[\kappa]$.
For any $<$-increasing map $g:\kappa\to \kappa$ with
$\mathbb{G}\cong(g[\kappa],E\restrict g[\kappa])$,
the composition $\sigma\circ g\circ\sigma^{-1}:S\to S$
is a pnp map.
\end{lemma}

\begin{proof}
Let $f:=\sigma\circ g\circ\sigma^{-1}$ for some
$<$-increasing isomorphism of $\mathbb{G}$ into itself.
Suppose $s,t\in S$ and $\beta:=\lg(s)<\lg(t)=\alpha$.  
Then $\sigma^{-1}(s)=\beta$ and $\sigma^{-1}(t)=\alpha$.
Since $g$ is $<$-increasing, $g(\beta)<g(\alpha)$.
Hence $\lg(f(s))=g(\beta)<g(\alpha)
=\lg(f(t))$.  Moreover, 
$t(\lg(s))=t(\beta)=1$ if and only if $\{\beta,\alpha\}\in E$.
Since $g$ is an isomorphism, $t(\lg(s))=1$ if and only if
$\{g(\beta),g(\alpha)\}\in E$.  By definition of tree embedding,
it follows that $t(\lg(s))=1$ if and only if $f(t)(\lg(s))=1$.
Thus $f$ is a pnp map from $S$ into $S$.
\end{proof}

By much the same reasoning, one can show the converse.

\begin{theorem}\label{thm.1.translate}[Translation Theorem]
Suppose $\mathbb{G}=(\kappa,E)$ is a $\kappa$-Rado graph 
with tree embedding $\sigma$ and $S=\sigma[\kappa]$.
For any  pnp map $f:S\to S$,   
the composition $g:=\sigma^{-1}\circ f\circ\sigma:\kappa\to \kappa$
is an $<$-increasing map with 
$\mathbb{G}\cong(g[\kappa],E\restrict g[\kappa])$.
\end{theorem}

\begin{proof}
Let $g:=\sigma^{-1}\circ f\circ\sigma$ for some
pnp map $f:S\to S$.
Suppose $\beta<\alpha<\kappa$.  Then $\lg(\sigma(\beta))=\beta
<\alpha=\lg(\sigma(\alpha))$.  Since $f$ is a pnp map,
$\sigma^{-1}\circ f\circ\sigma(\beta)=\lg(f(\sigma(\beta)))<
\lg(f(\sigma(\alpha)))=\sigma^{-1}\circ f\circ\sigma(\alpha)$,
so $g$ is a $<$-increasing map.

By the definition of the tree embedding, $\{\beta,\alpha\}$
is an edge of $\mathbb{G}$ if and only if 
$\sigma(\alpha)(\lg(\sigma(beta)))=1$.
Since $f$ is a pnp map, it follows that 
$\{\beta,\alpha\}$
is an edge of $\mathbb{G}$ if and only if
$f(\sigma(\alpha))(\lg(f(\sigma(beta))))=1$.
Apply the definition of tree embedding to
$g(\alpha)=\sigma^{-1}(f(\sigma(\alpha)))$
and 
$g(\beta)=\sigma^{-1}(f(\sigma(\beta)))$, to see that
$\{\beta,\alpha\}$
is an edge of $\mathbb{G}$ if and only if
$\{g(\beta),g(\alpha)\}$ is an edge.

Thus $g$ is a $<$-increasing isomorphism of $\mathbb{G}$
into itself.
\end{proof}

The next definition identifies sufficient conditions
for a map to carry a strongly diagonal set to one
of the same $m$-type.  

\begin{definition}\label{def.polite}
Call a map $f:\krtree\to \krtree$ \emph{polite} 
if it satisfies the 
following 
conditions for all $x,y,u,v$: 
\begin{enumerate}
\item (preservation of lexicographic order) if 
$x$ and $y$ are incomparable and $x\lex y$, then 
$f(x)$ and $f(y)$ are incomparable and $f(x)\lex (y)$;
\item (meet regularity) if $\set{x,u,v}$ is
diagonal and 
$x\meet u=x\meet v$, then $f(x)\meet f(u) =f(x)\meet f(v)$;
\item (preservation of meet length order) if $\lg(x\meet y)<\lg(u\meet v)$,
then $\lg(f(x)\meet f(y))<\lg(f(u)\meet f(v))$.
\end{enumerate}
Call it \emph{polite to strongly diagonal sets} if it is a pnp map 
which satisfies the above conditions for all $x,y,u,v$ 
with $\set{x,y,u,v}$ a strongly diagonal set.
\end{definition}

The next lemma follows immediately from the above definition.

\begin{lemma}\label{lem.preserve.polite}
Strong embeddings are polite 
and polite to strongly diagonal sets.  The collection of polite embeddings
is closed under composition as is the collection of embeddings
polite to strongly diagonal sets.
\end{lemma}

\begin{lemma}\label{lem.diag.to.diag}
Suppose $\phi:\krtree\to\krtree$ is a map 
which is polite to strongly diagonal sets and
whose image is a strongly diagonal set.
Then for any strongly diagonal set $A$, $\clp(A)=\clp(\phi[A])$ and
there is a pnp map $\overline{\phi}:
\closure{A}\to\closure{(\phi[A])}$ such that for all
$x,y$ in $A$, $\overline{\phi}(x\meet y)
=\phi(x)\meet\phi(y)$. 
\end{lemma}

\begin{proof}
Fix a strongly diagonal set $A$.  
Let $\overline{\phi}(a)=\phi(a)$ for $a\in A$ and let
$\overline{\phi}(a\meet b)=\phi(a)\meet\phi(b)$ for $a,b\in A$.
The proof of the lemma for $A$ proceeds by a series of claims.

\begin{claim}
The map $\overline{\phi}$ 
is well-defined.
\end{claim}

\begin{proof}
We must show that
for all $x,y,u,v\in A$, if $x\meet y=u\meet v$, then $\phi(x)\meet \phi(y)
=\phi(u)\meet\phi(v)$. 
If $x=y$ or $u=v$, then $x=y=u=v$ since $A$ is diagonal.  In this case
the claim follows immediately, so assume $x\ne y$ and $u\ne v$.
If $\set{x,y}\cap\set{u,v}$ is non-empty, then the claim follows
from the assumption of meet regularity.  So assume 
$\set{x,y}\cap\set{u,v}$ is empty.  Since $\set{x,y,u}$ is
a three element diagonal set and $x\meet y$ is an initial segment
of all three elements, either $x\meet y=x\meet u$ or $x\meet y=y\meet u$.
Hence by meet regularity, either $\phi(x)\meet\phi(y)
=\phi(x)\meet\phi(u)=\phi(u)\meet\phi(v)$ or
$\phi(x)\meet\phi(y)
=\phi(y)\meet\phi(u)=\phi(u)\meet\phi(v)$, and the claim follows.
\end{proof}

\begin{claim}
The map $\overline{\phi}$ preserves length order.
\end{claim}

\begin{proof}
Assume $s,t\in \closure{A}$ satisfy $\lg(s)<\lg(t)$.
Let $x,y\in A$ be such that $s=x\meet y$ and $u,v\in A$ be such that
$t=u\meet v$.  By preservation of meet length order, since 
$\lg(x\meet y)<\lg(u\meet v)$, it follows that $\lg(\overline{\phi}(s))
=\lg(\phi(x)\meet\phi(y))<\lg(\phi(u)\meet\phi(v))=\lg(\overline{\phi}(t))$.
\end{proof}

\begin{claim}
The map $\overline{\phi}$ is a pnp map
such that for all
$x,y$ in $A$, $\overline{\phi}(s\meet t)
=\phi(s)\meet\phi(t)$. 
\end{claim}

\begin{proof}
By the definition of $\overline{\phi}$ and the previous claims, 
it is enough to show that $\overline{\phi}$
preserves passing numbers.  Suppose $s$ and $t$ are in $\closure{A}$
and $\lg(s)<\lg(t)$.
Since any element $t'$ of the closure of $A$ is either in $A$ or has a 
proper extension $t$ in $A$ with 
$\overline{\phi}(t')\subseteq \overline{\phi}(t)=\phi(t)$, 
we may assume without loss of generality that $t$ is in $A$.  
If $s\in A$, then the conclusion follows since $\phi$ is a pnp map.

So suppose $s=x\meet y$ for $x$ and $y$ distinct elements of $A$.
Consider $s\meet t$.  Either it has the same length as $s$
or it is shorter.

If $\lg(s\meet t)<\lg(s)$, then 
$\lg(\overline{\phi}(s\meet t))<
\lg(\overline{\phi}(s))$ by the previous claim.  In this case,
$\overline{\phi}(t)(\lg(\overline{\phi}(s)))=0
=t(\lg(s))$, since $A$ and its image under $\phi$ are both
strongly diagonal.

If $\lg(s\meet t)=\lg(s)$, then $s\subseteq t$. 
Let $w\in A$ be such that $s=t\meet w$.  Then the value of
$t(\lg(s))$ and $\overline{\phi}(t)(\lg(\overline{\phi}(s)))$
are determined by the lexicographic order of the pairs $t,w$
and ${\phi}(t),{\phi}(w)$.  Since $\phi$ preserves lexicographic
order, $\overline{\phi}(t)(\lg(\overline{\phi}(s)))=t(\lg(s))$,
as required.
\end{proof}
%
%

\begin{claim} $\clp(A)=\clp(\phi[A])$.
\end{claim}

\begin{proof}
Enumerate $\closure{A}$ in increasing order of length 
as $\seq{a_\alpha:\alpha<\mu}$ for some $\mu<\kappa$.  
Let $B=\phi[A]$.  Then $B$ is a strongly
diagonal set since it is a subset of a strongly diagonal set.  
Enumerate $\closure{B}$ in increasing order of length as
$\seq{b_\beta:\beta<\nu}$ for some $\nu<\kappa$.  

Since $\overline{\phi}$ is a pnp map that carries $\closure{A}$
onto $\closure{B}$,
it is a bijection from $\closure{A}$ to $\closure{B}$.
Since the order type of $\set{a\in \closure{A}:\lg(a)<\lg(a_\alpha)}$
is $\alpha$ and the order type of 
$\set{b\in \closure{B}:\lg(b)<\lg(b_\beta)}$
is $\beta$, it follows that 
$\overline{\phi}(a_\alpha)=b_\beta$ and $\mu=\nu$.

For $\alpha<\mu$, let $A_\alpha:=\set{a\restrict\lg(a_\alpha): a\in A}$
and $B_\alpha:=\set{b\restrict\lg(b_\alpha): b\in B}$.  Let $A_\mu=A$
and $B_\mu=B$.  Use induction to prove that for all positive $\alpha\le\mu$,
$\clp(A_\alpha)=\clp(B_\alpha)$.  To start the induction, observe that
$|A_0|=|B_0|=1$, so $\clp(A_0)=\set{\emptyset}=\clp(B_0)$.  For the limit
case, assume $\alpha$ is a limit ordinal and for all $\beta<\alpha$,
$\clp(A_\beta)=\clp(B_\beta)$.  In this case, $\clp(A_\alpha)=\clp(B_\alpha)$,
since 
$\clp(A_\alpha)=\bigcup\set{\clp(A_\beta):\beta<\alpha}$ and
$\clp(B_\alpha)=\bigcup\set{\clp(B_\beta):\beta<\alpha}$.
For the successor case, assume $\alpha=\beta+1$ and $\clp(A_\beta)
=\clp(B_\beta)$.  Let $\sigma=\sigma_\alpha$ be the increasing
enumeration of $\set{\lg(a_\gamma):\gamma<\alpha}$ and
let $\tau=\tau_\alpha$ be the increasing
enumeration of $\set{\lg(b_\gamma):\gamma<\alpha}$.
Then $\clp(A_\alpha)$ is the similarity tree whose set of leaves
is $\set{a\circ \sigma:a\in A_\alpha}$
and $\clp(B_\alpha)$ is the similarity tree whose set of leaves
is $\set{b\circ \sigma:b\in B_\alpha}$.
Moreover $\clp(A_\beta)$ is the similarity tree whose set of leaves
is $\set{a\circ (\sigma\restrict\beta):a\in A_\beta}$
and $\clp(B_\beta)$ is the similarity tree whose set of leaves
is $\set{b\circ (\sigma\restrict\beta):b\in B_\beta}$.
For $a\in A_\alpha\setminus A_\beta$, $a\circ\sigma(\beta)$
is the passing number of $a$ at $\lg(a_\beta)$.
Since $\overline{\phi}$ preserves passing numbers and
carries elements of $A$ to elements of $B$ and
elements of $\closure{A}\setminus A$ to elements
of $\closure{B}\setminus B$, it follows that $\clp(A_\alpha)=\clp(B_\alpha)$
in this case as well.  Therefore by induction, 
$\clp(\closure{A})=\clp(\closure{B})$, and since 
$\clp(A)=\clp(\closure{A})$ and $\clp(B)=\clp(\closure{B})$,
the claim follows.
\end{proof}

Since $A$ was an arbitrary strongly diagonal set,
by the last two claims, the lemma follows.
\end{proof}

Our next goal is to prove the existence of a \emph{pnp diagonalization}.

\begin{definition}\label{def.thin.diag.map}
Suppose $w\in\krtree$ and $S\subseteq\krtree$.
Call $f$ a 
\emph{pnp diagonalization into $S\cap\cone(w)$}
if $f$ is a polite injective $\qle$-preserving pnp map
whose range is a strongly diagonal subset $D$ with
$\closure{D}\subseteq S\cap\cone(w)$.  
Call $f$ a \emph{pnp diagonalization} if it is a
pnp diagonalization into $\krtree\cap\cone(\emptyset)$.
\end{definition}

An extra quality we desire for our diagonalization is \emph{level
  harmony}, which will be used in the section on lower bounds.

\begin{definition}\label{def.harmony}
Suppose $f:\krtree\to\krtree$ is an injective map.
Define $\hat{f}:\krtree\to\krtree$ by
$\hat{f}(s)=f(s\cat\langle 0\rangle)\meet f(s\cat\langle 1\rangle)$.
The function \emph{$f$ has level harmony} if
$\hat{f}$ is an extension and $\lex$-order preserving map such that
for all $s,t\in\krtree$, the following conditions hold:
\begin{enumerate}
\item $\hat{f}(s)\subsetneq f(s)$;
\item $\lg(s)<\lg(t)$ implies $\lg(f(s))<\lg(\hat{f}(t))$;
\item $\lg(s)=\lg(t)$ implies $\lg(\hat{f}(s))<\lg(f(t))$.
\end{enumerate}
\end{definition}

\begin{lemma}\label{lem.2.diag.map}[Second diagonalization lemma]
Suppose $w\in\krtree$ and that $S\subseteq\krtree$ is cofinal and transverse.
Then there is a pnp diagonalization into $S\cap\cone(w)$ which
has level harmony. 
\end{lemma}

\begin{proof}
Our plan is to approach the problem in pieces
by using recursion to define three functions, $\varphi_0,\varphi_1,\varphi:
\krtree\to S$ so that $\varphi$ is the desired diagonalization,
$\hat{\varphi}=\varphi_0$, $\varphi_1(t)$ is the minimal
extension in $S$ of $\varphi_0\cat\langle 1\rangle$, 
$\varphi(t)\meet\varphi(t\cat\langle 1\rangle)=\varphi_1(t)$,
and $\varphi(t)\lex\varphi(t\cat\langle 1\rangle)$.

For notational convenience, if $\varphi$ has been defined
on ${}^{\alpha>}2$, then we let $\ell_0(\alpha)$ be the
least $\theta$ such that $\lg(\varphi(t))<\theta$ for
all $t\in {}^{\alpha>}2$.  Also, if $\varphi_1$ has been defined
on ${}^{\alpha\ge}2$, then we let $\ell_1(\alpha)$ be the
least $\theta$ such that $\lg(\varphi_1(t))<\theta$ for
all $t\in {}^{\alpha\ge}2$.  

Let $\prec$ be a well-ordering of the levels of $\krtree$.
We use recursion on $\alpha<\kappa$ to define  the restrictions
to ${}^\alpha 2$ of $\varphi_0$, $\varphi_1$
and $\varphi$ so that the following 
properties hold:
\begin{enumerate}
\item extension and lexicographic order: 
\begin{enumerate}
\item the restriction of $\varphi_0$ to ${}^{\alpha\ge} 2$ 
is extension and $\lex$-order preserving;
\item for all $s\in{}^\alpha 2$, $\varphi_1(s)$ is the minimal
extension in $S\cap\cone(w)$ of $\varphi_0(s)\cat\langle 1\rangle$;
\item for all $s\in{}^\alpha 2$, $\varphi(s)$ is 
an extension in $S\cap\cone(w)$ of $\varphi_1(s)\cat\langle 0\rangle$;
\item for all $s\in{}^{\alpha>}2$, 
$\varphi_0(s\cat\langle 0\rangle)$ is 
an extension of $\varphi_0(s)\cat\langle 0\rangle$ and
\newline
$\varphi_0(s\cat\langle 1\rangle)$ is 
an extension of $\varphi_1(s)\cat\langle 1\rangle$;
\end{enumerate}
\item length order: 
\begin{enumerate}
\item for all $t\in{}^\alpha 2$ and $s\in{}^{\alpha\ge}2$, 
if $s\prec t$, then 
\newline
$\ell_0(\lg(s))\le\lg(\varphi_0(s))<\lg(\varphi_0(t))$ and 
\newline
$\ell_1(\lg(s))\le\lg(\varphi(s))<\lg(\varphi(t))$;
\end{enumerate}
\item passing number: 
\begin{enumerate}
\item for all $t\in{}^\alpha 2$ and $s\in{}^{\alpha\ge}2$, 
if $s\prec t$ and $s\not\subseteq t$, then 
\newline
$\varphi_0(t)(\lg(\varphi_0(s))=0$ and $\varphi_0(t)(\lg(\varphi_1(s))=0$;
\item for all $t\in{}^\alpha 2$ and $s\in{}^{\alpha\ge}2$, 
if $s\prec t$ and $s\not\subseteq t$, then 
\newline
$\varphi(t)(\lg(\varphi_0(s))=0$ and $\varphi(t)(\lg(\varphi_1(s))=0$;
\item for all $t\in{}^\alpha 2$ and $s\in{}^{\alpha>}2$, 
$\varphi_0(t)(\lg(\varphi(s))=t(\lg(s))$;
\item for all $s,t\in {}^\alpha 2$, if $s\prec t$, then 
$\varphi(t)(\lg(\varphi(s))=0$.
\end{enumerate}
\end{enumerate}

Suppose $\alpha<\kappa$ is arbitrary and for all $\beta<\alpha$,
the restrictions to ${}^\beta 2$ of $\varphi_0$, $\varphi_1$
and $\varphi$ have been defined. To maintain length order,
we first define $\varphi_0$ and $\varphi_1$ 
by recursion on $\prec$ restricted to level ${}^\alpha 2$.
So suppose $\lg(t)=\alpha$ and for all $s\prec t$, 
$\varphi_0$ and $\varphi_1(s)$ have been defined.

Use extension and $\lex$-order properties to identify
an element $\varphi_0^-(t)$ of which $\varphi_0(t)$ is
to be an extension.  If $\alpha=0$, set $\varphi^-_0(t)=\emptyset$.
If $\alpha$ is a limit ordinal, let $\varphi^-_0(t)$ be
$\bigcup\set{\varphi_0(t\restrict\beta):\beta<\alpha}$.
If $\alpha$ is a successor ordinal and $t=t^-\cat\langle \delta\rangle$,
let $\varphi^-_0(t)=\varphi_\delta(t^-)\cat\langle\delta\rangle$.

Next determine an ordinal
$\gamma_0(t)$ sufficiently large that if $\varphi_0(t)$
is at least that length, it will satisfy the
length order property.
If $t$ is the $\prec$-least element of length $\alpha$,
let $\gamma_0(t)=\ell_0(\alpha)$.  If $t$ has a $\prec$-immediate
predecessor $t'$ of length $\alpha$, let $\gamma_0(t)=
\lg(\varphi_1(t'))+1$.  If $t$ is a $\prec$-limit  of
elements of length $\alpha$, then let $\gamma_0(t)$ be
the supremeum of $\lg(\varphi_1(s))+1$ for $s$ of length $\alpha$
with $s\prec t$.  

Next define an extension $\varphi^+_0(t)$ of $\varphi^-_0(t)$
of length $\gamma_0(t)$ so that the passing number properties
are satisfied by $\varphi^+_0(t)$.  If $\alpha=0$,
then $t=\emptyset$, $\varphi^-_0(t)=\emptyset$,
$\gamma_0(0)=0$ and $\varphi^+_0(t)=\emptyset$.  If $\alpha>0$
is a limit ordinal, then by induction on $\beta<\alpha$,
$\ell_0(\beta)$ is an increasing sequence.  Moreover,
the limit of this sequence is the length of $\varphi^-_0(t)$.
It follows that $\varphi^-_0(t)$ satisfies the passing
number properties for $s\in {}^{\alpha>}2$.  Let $\varphi^+_0(t)$
be the sequence extending $\varphi^-_0(t)$ by zeros, as needed,
to a length of $\gamma_0(t)$.  If $\alpha$ is a successor ordinal
and $t=t^-\cat\langle\delta\rangle$, then let $\varphi^+_0(t)$
be the extension of $\varphi^-_0(t)$ of length $\gamma_0(t)$
such that for all $\eta$ with $\lg(\varphi^-_0(t))\le\eta<\gamma_0(t)$,
$\varphi^+_0(\eta)=\delta$ if $\eta=\lg(\varphi(s))$ for
some $s$ with $\lg(s)+1=\alpha$, and 
$\varphi^+_0(\eta)=0$ otherwise.

Next let $\varphi_0(t)$ be an extension in $S\cap\cone(w)$ of 
$\varphi^+(t)$
and let $\varphi_1(t)$ be an extension in $S\cap\cone(w)$ of 
$\varphi_0(t)\cat\langle 1\rangle$ as required by
the extension and lexicographic order properties.
The careful reader may now check that the various
properties hold for the restrictions of $\varphi_0$ and $\varphi_1$
to ${}^\alpha 2$.

Use a similar process to define the restriction of $\varphi$
to ${}^\alpha 2$ by recursion on $\prec$ restricted to ${}^\alpha 2$. 
Suppose that $\lg(t)=\alpha$ and for all $s\prec t$,
$\varphi(s)$ has been defined.  Let $\varphi^-(t)=
\varphi_1(t)\cat \langle 0\rangle$.  

If $t$ is the $\prec$-least element of length $\alpha$,
let $\gamma_1(t)=\ell_1(\alpha)$.  If $t$ has a $\prec$-immediate
predecessor $t'$ of length $\alpha$, let $\gamma_1(t)=
\lg(\varphi(t'))+1$.  If $t$ is a $\prec$-limit  of
elements of length $\alpha$, then let $\gamma_1(t)$ be
the supremum of $\lg(\varphi(s))+1$ for $s$ of length $\alpha$
with $s\prec t$.  

Next define an extension $\varphi^+(t)$ of $\varphi^-(t)$
of length $\gamma_1(t)$ so that the passing number properties
are satisfied by $\varphi^+(t)$.  If $\alpha=0$,
there are no passing number properties that need be checked, and we
set $\varphi^+(t)=\varphi^-(t)$.  If $\alpha>0$,
then let $\varphi^+(t)$ be the extension by zeros
of $\varphi^-(t)$ of length $\gamma_1(t)$.  
Since $\varphi_0(t)$ and $\varphi_1(t)$ satisfy
the passing numbers properties, it follows that
$\varphi^+(t)$ does as well, since all passing numbers
longer than $\lg(\varphi^-(t))$ will be zero.

Finally let $\varphi(t)$ be an extension in $S\cap\cone(w)$ of 
$\varphi^+(t)$.
The careful reader may now check that the various
properties hold for the restriction of $\varphi$
to ${}^\alpha 2$.

This completes the recursive construction of $\varphi_0$,
$\varphi_1$ and $\varphi$.  By induction, the various
properties hold for all $\alpha<\kappa$.

Thus $\varphi_0=\hat\varphi$ is extension and
$\lex$-order preserving, and by the length order
property, different elements of the
union of the ranges of $\varphi_0$, $\varphi_1$
and $\varphi$ have different lengths.  It also follows
that these three maps are injective.  Moreover
the union of their ranges is a subset of $S\cap\cone(w)$.
By the extension and lexicographic order properties
and the length order property, $\varphi$
has level harmony.

By the passing number properties, $\varphi$
is a pnp map.  By the extension and lexicographic order
properties, $\varphi$ preserves $\qle$-order.

By the extension and lexicographic order
properties, $\varphi$ 
carries incomparable elements into incomparable elements
and preserves $\lex$-order.
By the length order property, $\varphi$
preserves meet length order.  Since $\hat{\varphi}=\varphi_0$
preserves extension, $\varphi$ satisfies meet regularity.
Thus $\varphi$ is polite.

From the extension and lexicographic order 
properties, it follows that the meet closure of 
$D:=\ran(\varphi)$ is the union of the ranges
of $\varphi$, $\varphi_0$ and $\varphi_1$ and
all elements of the range of $\varphi$ are incomparable.
Hence $D$ is an antichain and $\closure{D}\subseteq S\cap\cone(w)$ 
is transverse, so $D$ is diagonal.
Note that passing numbers of $1$ were introduced
only to keep $\varphi_0$ extension and $\lex$-order
preserving, to ensure $\varphi(t)\lex\varphi_1(t)$ so
that $\qle$-order is preserved, 
and to ensure $\varphi$ is a pnp map.
It follows that $\closure{D}$ is strongly diagonal.

Therefore, $\varphi$ is the required pnp diagonalization into
$S\cap\cone(w)$ with level harmony.
\end{proof}

\begin{lemma}\label{lem.poly.diag}
Suppose $S\subseteq\krtree$ is cofinal and transverse.
Then there are a diagonal set $D$, maps $\seq{\varphi_t: t\in \krtree}$
and a pre-$S$-vip order $\prec$ such that 
for all $t\in\krtree$, the following conditions hold:
\begin{enumerate}
\item $\varphi_t$ is a pnp diagonalization with level harmony
into $S\cap\cone(t)$; 
\item the meet closure of the set 
$D_t:=\ran(\varphi_t)$ is a subset of $S$ 
disjoint from $\closure{D}_s$ for all $s\prec t$;
and 
\item $\prec$ is a $D_t$-vip order.
\end{enumerate}
\end{lemma}

\begin{proof}
Use Lemma \ref{lem.exist.vip} to find $\prec^*$, a
pre-$S$-vip order on $\krtree$.  By Lemma \ref{lem.cofinal.dense},
$S$ is dense.

Apply the Second Diagonalization Lemma \ref{lem.2.diag.map}
to each $t\in\krtree$ to obtain $\varphi^*_t$, a  pnp diagonalization 
into $S\cap\cone(t)$ which has level harmony. 
Use recursion on $\prec$
to define $\pi:\krtree\to\krtree$, $\seq{\varphi_t:t\in\krtree}$
and $\seq{D_t:t\in \krtree}$ 
such that for all $t\in\krtree$, $D_t:=\ran(\varphi_t)$,
$\pi(t)$ is an extension of $t$ with $\cone(\pi(t))$
disjoint from the union over all $s\prec t$ of
$\closure{D}_s$.  Since the order type of
$\set{s\in\krtree:s\prec t}$ is less than $\kappa$ and
each $D_s$ is a strongly diagonal set whose meet closure is a of $S$,
it is always possible to continue the recursion.

Use induction on the recursive construction to show
that the meet closures of the sets $D_t$ are disjoint.  

Let $D=\bigcup\set{\closure{D}_t:t\in \krtree}$.  
Then $D$ is transverse since it is a subset of $S$
and $S$ is transverse. Let $\prec$ agree with $\prec^*$ on 
all pairs from different levels, and use recursion on 
$\alpha<\kappa$ to define $\prec$ from $\prec^*$ as follows. 
If there is no element
of $D$ in $\pre\alpha2$ then $\prec$ and $\prec^*$ agree
on $\pre\alpha2$.  

So suppose $d\in\closure{D}_t$ and $\lg(d)=\alpha$.
For each $\beta<\alpha$, let $C(\beta)$ be the set of all
$x\in \pre\alpha2$ such that $x\restrict\beta=d\restrict\beta$
and $x(\beta)\ne d(\beta)$.  Since $\prec^*$ is a pre-$S$-vip
order, if $\beta<\gamma<\alpha$, then $C(\beta)\prec^* C(\gamma)$
in the sense that for every element $x$ of $C(\beta)$
and $y$ of $C(\gamma)$, one has $x\prec^* y$.
Use the fact that $D_t$ is a diagonal set to
partition $\pre\alpha2=\{ d\}\cup A_\alpha(0)\cup A_\alpha(1)\cup
A_\alpha(2)$ into disjoint pieces where for $\delta<2$,
$A_\alpha(\delta):=\set{u\restrict\alpha:u\in D_t\land u(\alpha)=\delta}$.
Let the restriction of $\prec$ to $\pre\alpha2$ be such that 
for each $\beta<\alpha$, 
\[
C(\beta)\cap A_\alpha(0)\prec
C(\beta)\cap A_\alpha(1)\prec
C(\beta)\cap A_\alpha(2)
\] 
and otherwise $\prec$ agrees with $\prec^*$.
Then the restriction of $\prec$ to $\pre\alpha2$ is a well-order,
since the restriction of $\prec^*$ is and because $\prec^*$
is a pre-$S$-vip order.  

Since the restriction of $\prec$ to each level is a well-order,
it follows that $\prec$ is a well-ordering of the levels of $\krtree$.
For each $t\in\krtree$, 
since $\prec^*$ is a pre-$S$-vip order,
it follows that 
$\prec$ is a pre-$\closure{D_t}$-vip order,
so by construction, $\prec$ is a $D_t$-vip order.
\end{proof}

\begin{theorem}\label{thm.upper.bound.2} 
Let $m\ge 2$ and
  suppose that $\kappa$ is a cardinal which is measurable
  in the generic
  extension obtained by adding $\lambda$ Cohen subsets of $\kappa$,
  where $\lambda\to(\kappa)^{2m}_{2^{\kappa}}$. 
Then for $r_m^+$ equal to the number of vip $m$-types,
any $\kappa$-Rado graph $\mathbb{G}=(\kappa,E)$
satisfies 
\[
\mathbb{G}\rightarrow(\mathbb{G})^m_{<\kappa,r_m^+}.
\]
\end{theorem}

\begin{proof}
Let $\sigma:\kappa\to\krtree$ be the tree embedding and set
$S=\ran\sigma$.  Then by Lemma \ref{lem.sigma.cofinal},
$S$ is cofinal and transverse.  
Apply Lemma \ref{lem.poly.diag} to obtain a pre-$S$-vip
order $\prec$ and a sequence $\seq{\varphi_t:t\in\krtree}$
such that for all $t\in\krtree$, the three listed
properties of the lemma hold.

Fix a coloring $c:[\kappa]^m\to\mu$  where $\mu<\kappa$.  
Define $d:[\krtree]^m\to\mu$ by $d(z)=c(\sigma^{-1}[\varphi_0[z]])$.  

Apply Shelah's
Theorem \ref{thm.SHL} to the restriction of $d$
to $m$-element antichains to obtain
a strong embedding $e$ and a node $w$ such that $e$ preserves
$\prec$ on $\cone(w)$ and $d(e[a])=d(e[b])$ for all 
$\prec$-similar $m$-element antichains $a$ and $b$ of $\cone(w)$.

Now $D_w:=\ran(\varphi_w)$ is a strongly diagonal set 
and $\prec$ is $D_w$-vip.
Let $D=e[\varphi_w[S]]$.  Since $\varphi_w[S]$ is a subset of $D_w$,
$\varphi_w[S]$ is a strongly diagonal set since $D_w$ is.
Also,  $\prec$ is $\varphi_w[S]$-vip level order on the downwards
closure under initial segments of $\varphi_w[S]$, 
since $\prec$ is a $D_w$-vip level order.
Since $e$ is a strong embedding, $D$ is
a strongly diagonal set.  Since $e$ preserves $\prec$ on $\cone(w)$
and $\varphi_w[S]\subseteq D_w\subseteq\cone(w)$, it follows that 
the restriction of $\prec$ to the downward closure of $D$
is a $D$-vip order.  Hence for all $x\in[D]^m$, 
the ordered similarity type $(\clp(x),\prec_x)$ is
a vip $m$-type.  

Finally let $K:=\sigma^{-1}[\varphi_0[D]]=
\sigma^{-1}\circ \varphi_0\circ e\circ\varphi_w\circ \sigma[\kappa]$.  
Note that  $D=\varphi_0^{-1}[\sigma[K]]$, so for all $m$-element
subsets $u$ of $K$, the image, $x=\varphi_0^{-1}[\sigma[K]][u]\subseteq D$,
is a strongly diagonal set, and $(\clp(x),\prec_x)$
is a vip $m$-type.

Since $\varphi_0$, $e$ and $\varphi_w$ are all pnp maps, so is
their composition.  By the First Translation Theorem
\ref{thm.1.translate}, this mapping is an isomorphism
of $\mathbb{G}$ into itself. 

We claim that $c$ is constant on $m$-element subsets of $K$
whose images under $\varphi_0^{-1}\circ\sigma$ are $\prec$-similar.
Consider two such $m$-element subsets $u,v$ of $K$.
Let $a'$ and $b'$ be the subsets of $\kappa$ for which
$\sigma^{-1}\circ \varphi_0\circ e\circ\varphi_w\circ \sigma[a']=u$
and 
$\sigma^{-1}\circ \varphi_0\circ e\circ\varphi_w\circ \sigma[b']=v$.
Let $a=\varphi_w[\sigma[a']]$ and $b=\varphi_w[\sigma[b']]$.
Then $a$ and $b$ are $m$-element strongly diagonal subsets
of $\cone(w)$.  Notice that $u=\sigma^{-1}[\varphi_0[e[a]]]$
and $v=\sigma^{-1}[\varphi_0[e[b]]]$.  Thus $e[a]$ and $e[b]$
are $\prec$-similar.  Since $e$ is a strong embedding
which preserves $\prec$ on $\cone(w)$, it follows that
$a$ and $b$ are $\prec$-similar.  By the application of
Shelah's Theorem above, $d(e[a])=d(e[b])$.  From the 
definition of $d$ it follows that 
$c(u)=c(\sigma^{-1}[\varphi_0[e[a]]]=d(e[a])$ and 
$c(v)=c(\sigma^{-1}[\varphi_0[e[b]]]=d(e[b])$.
Consequently, $c(u)=c(v)$.

Since the image under $\varphi_0^-1\circ\sigma$ of any
$m$-element subset of $K$ is similar to a  vip $m$-type 
and any two $\prec$-similar subsets receive the
same color from $c$, the $m$-element subsets of $K$
are colored with at most $t_m^+$ colors, so the
theorem follows.
\end{proof}

\section{Lower bounds for Rado graphs}\label{sec.lower.bnd.2}
The computation of lower bounds for Rado graphs is a bit
more complicated than the computation for $\kappa$-dense linear orders.
We reduce the problem by showing for suitable $\kappa$
that if $D\subseteq\krtree$ is the range 
of a pnp diagonalization with level harmony and $\prec$
is a $D$-vip level order, then every vip $m$-type is
realized as $(\clp(x),\prec_x)$ for some $x\subseteq D$.
This theorem is the companion to Theorem \ref{thm:Qk}.
Its proof uses a pnp diagonalization with level harmony 
in place of a semi-strong embedding.

\begin{theorem}\label{thm:Gk}
  Suppose that $\kappa$ is a cardinal which is measurable
  in the generic
  extension obtained by adding $\lambda$ Cohen subsets of $\kappa$,
  where $\lambda\to(\kappa)^{6}_{2^{\kappa}}$. 
Further suppose $f:\krtree\to\krtree$ 
is a pnp diagonalization with level harmony, $D:=f[\krtree]$,
and $\prec$ is a $D$-vip order of the levels of $\krtree$. 
Then every vip $m$-type $(\tau,\dless)$ 
is realized as $(\clp(x),\prec_x)$ for some $x\subseteq D$.
\end{theorem}


\begin{proof}
For $t\in \pre \alpha2$, $i=0,1$ and $\delta=0,1$, 
define well-orderings $\prec^{i,\delta}_t$
on $\pre\alpha2$ as follows. Let $\beta^0_t=
\lg(\hat{f}(t))$ and
set $\beta^1_t=\lg({f}(t))$.
Then $s\prec^{i,\delta}_t s'$ if and only if 
$f(s\cat\langle \delta\rangle)\restrict\beta^i_t
\prec
f(s'\cat\langle \delta\rangle)\restrict\beta^i_t$.

Let $\prec^{\,\prime}$ be any \compact well-ordering of the 
levels of $\krtree$.
Call a triple $\set{s,s',t}$ \emph{local} if 
$\lg(s)=\lg(s')=\lg(t)$, $s\lex s'$, $t\prec^{\,\prime} s$,
and $t\prec^{\,\prime} s'$, 
and $s\meet s'\not\subseteq t$.  

Let $d$ be a coloring of the triples of $\krtree$  defined as follows:
if $\set{s,s',t}$ is not local, let $d^{\,i,\delta}(\set{s,s',t}):=2$ and 
otherwise set $d^{\,i,\delta}(\set{s,s',t}):=
\bool{s\prec^{\,\prime} s'\iff s\prec^{i,\delta}_t s'}$.
For $b=\set{s,s',t}\in[\pre\alpha2]^3$, define $d(b):=
(d^{\,0,0}(b),
d^{\,0,1}(b),
d^{\,1,0}(b),
d^{\,1,1}(b))$.

Apply Shelah's Theorem \ref{thm.SHL} to $d$ and $\prec^{\,\prime}$
to obtain a strong embedding $e:\krtree\to\krtree$ and
an element $w$ so that
for triples from $T:=e[\cone(w)]$, the coloring depends only
on the $\prec^{\,\prime}$-ordered similarity type of the triple.  

Then two local triples $\set{s,s',t}$ and $\set{u,u',v}$ of $T$
are colored the same if and only if for all $i=0,1$ and $\delta=0,1$,
\begin{equation*}
s\prec_{t}^{i,\delta}s'\iff u\prec_{v}^{i,\delta}u'.
\end{equation*}
 
Thus for $t\in T$, the orderings $\prec_t^{i,\delta}$ must always 
agree with one of  $\prec^{\,\prime}$ and its converse on $T$
on pairs $\{ s,s'\}\subseteq T$ with $\{s,s',t\}$ local and
$s\prec^{\,\prime} s'$.  Similarly, they must always agree
with one of  $\prec^{\,\prime}$ and its converse on $T$
on pairs $\{ s,s'\}\subseteq T$ with $\{s,s',t\}$ local and
$s'\prec^{\,\prime} s$.  Since $\prec$ is a well-order,
all of the orderings $\prec_t^{i,\delta}$ are also well-orders.  
Thus they always agree with $\prec^{\,\prime}$.

Let $(\tau,\dless)$ be an arbitrary vip $m$-type.
Let $L$ be the set of leaves of $\tau$.
Then $\tau$ is a subtree of $\prele{2m-2}2$ and every level
of $\prele{2m-2}2$ has exactly one element of $\closure{L}$.
Extend $\dless$ defined on $\tau$ to 
$\dless^*$ defined on all
of $\prele{2n}2$ in such a way that the extension is still a 
$\closure{L}$-vip order.
 
Apply
Lemma \ref{lem.compact.embed} 
to get an order preserving strong embedding $j$
of $(\prele{2m-2}2,\dless^*)$ into $(\cone(w),\prec^{\,\prime})$.

Let $\seq{t_\ell:\ell\le 2m-2}$ enumerate the elements of $\closure{L}$
in increasing order of length.  Note that $\lg(t_\ell)=\ell$.
For $\ell\le 2m-2$, define
\[
\beta_\ell:=
\begin{cases}
\lg(\hat{f}(e(j(t_\ell))))&\text{if $t_\ell\notin L$,}\\
\lg(f(e(j(t_\ell))))&\text{if $t_\ell\in L$.}\\
\end{cases}
\]

Finally define $\rho:\tau\to\krtree$
by recursion on $\ell\le 2m-2$.
For $\ell=0$, let $\rho(\emptyset)=\hat{f}(e(j(\emptyset)))$.
For $\ell>0$, consider three cases for elements of
$\tau\cap{}^\ell 2$.  If $t_\ell\in L$,
let $\rho(t_\ell)=f(e(j(t_\ell)))$.  
If $t_\ell\notin L$, let $\rho(t_\ell)=\hat{f}(e(j(t_\ell)))$.  
Note that in both these cases,  $\beta_\ell=\lg(\rho(t_\ell))$. 
If $s\in\tau\setminus\closure{L}$ has length $\ell$, then there
is a unique immediate successor in $\tau$, $s\cat\langle\delta\rangle$.
In this case, let
$\rho(s)=f(e(j(s))\cat\langle\delta\rangle)\restrict\beta_\ell$. 

Since $j$ sends $\dless^*$-increasing pairs to
$\prec^{\,\prime}$-increasing pairs
and $e$ is a $\prec^{\,\prime}$-order preserving strong embedding, 
their composition sends sends $\dless^*$-increasing pairs to
$\prec^{\,\prime}$-increasing pairs.  
Since for $v_\ell=e(j(t_\ell))\in T$, the order
$\prec^{\,\prime}$ agrees with $\prec^i_{v_\ell}$ on $T\cap \pre
\gamma2$ where $\gamma=\lg(v_\ell)$, it follows that $\rho$
sends $\dless^*$-increasing pairs to $\prec$-increasing pairs.

Since $\hat{f}$ preserves extension and
lexicographic order and $\hat{f}(s)\subsetneq f(s)$, 
$\rho$ preserves extension and lexicographic order.
By construction $\rho$ sends levels to levels, meets to meets and 
leaves to leaves.   Let $x=\rho[L]$ be the image under
$\rho$ of the leaves of $\tau$.  By construction, $x\subseteq\ran(f)$. 
Also $(\clp(x),<_x)=(\tau,\dless)$, as required.  

Since $(\tau,\dless)$ was arbitrary, the theorem follows.
\end{proof}

We would like to define a coloring of the $m$-tuples of $\kappa$
using $t_m^+$ colors
such that for any $H\subseteq\kappa$ with $(H,E\restrict H)$
isomorphic to our Rado graph $\mathbb{G}$, every color
is the color of some $m$-tuple from $H$.

By Lemma \ref{lem.Rado.increasing}, there is an $<$-increasing
map $h:\kappa\to H$ such that the graph $(h[\kappa],E\restrict h[\kappa])$
is isomorphic to $\mathbb{G}$.  By Lemma \ref{lem.Rado.iso.1},
$g=\sigma\circ h\circ\sigma^{-1}$ is a pnp map.
By the previous theorem, every vip $m$-type can be
realized in the range of a pnp diagonalization with
level harmony.  By the Translation Theorem,
a pnp map gives rise to an induced subgraph of
the Rado graph which is isomorphic to the whole graph.

Our plan is to define a coloring $c$ using $t_m^+$
as follows. 
Apply Lemma \ref{lem.poly.diag} to obtain a pre-$S$-vip
order $\prec$ and a sequence $\seq{\varphi_t:t\in\krtree}$
such that for all $t\in\krtree$, the three listed
properties of the lemma hold.  In particular,
$\varphi_0$ is a pnp diagonalization into $S$
with level harmony, $D_0:=\varphi_0[\krtree]$
is a strongly diagonal set and $\prec$ is a $D$-vip order.
Enumerate the vip $m$-types as 
$\seq{(\tau_j,\dless_j):j<t_m^+}$; define $d$ on 
$m$-element subsets of $\pre\kappa2$
by $d(x)=j$ if $(\clp(\varphi_0[x]),\prec_{\clp(\varphi_0[x])})=
(\tau_j,\dless_j)$; and set $c(b)=d(\sigma[b])$. 

Then for every isomorphic copy $H$ of the Rado graph, 
the associated pnp map $\varphi_0\circ g$ described above has the 
property that different colors of $m$-element subsets 
of $K$ are distinguished by the ordered similarity types 
of elements of the range of $\varphi_0\circ g$.  
Thus it is enough  to show there is a pnp diagonalization $f$ 
with level harmony whose range is a subset of $\varphi_0\circ g$, and 
use Theorem \ref{thm:Gk} with $f$ to show that all the
colors appear in every isomorphic copy.
 
We show that for any pnp map $g$ with range a strongly diagonal set
there is a pnp diagonalization $f$ with level harmony with
$\ran(f)\subseteq \ran(g)$ in Theorem \ref{thm.thin.inside} below.
We work by successive approximation, using the next lemma.

\begin{lemma}\label{lem.curly.D}
Let $\mathcal{D}$ be the collection of 
pnp maps whose domain is $\pre\kappa2$ and  
whose range is a strongly diagonal subset of $\pre\kappa2$.
Then $\mathcal{D}$ is closed under composition
as are the following subfamilies:
\begin{enumerate}
\item maps in $\mathcal{D}$ 
with meet regularity;
\item maps in $\mathcal{D}$ 
with meet regularity that preserve lexicographic order of incomparable pairs; 
\item maps in $\mathcal{D}$ 
that are polite on strongly diagonal sets;
\item maps in $\mathcal{D}$ 
that are polite.
\end{enumerate}
Moreover, if $g$ is in one of these families and $e$ is a strong
embedding, then $e\circ g$ is also in that family.
\end{lemma}

\begin{proof}
Apply Lemma \ref{lem.preserve.polite} and 
use the definitions of terms to see that the various families 
are closed under composition. 
Since strong embeddings preserve all the properties in
question, and carry strongly diagonal sets to strongly
diagonal sets, their composition with any function in
one of the various families is also in that family.
\end{proof}

\begin{lemma}\label{lem.harmony}
Suppose $w\in \pre\kappa2$ and  $S$ is cofinal and transverse.
Further suppose $f$ is polite to strongly diagonal sets and has
range a strongly diagonal set $D$ with $\closure{D}\subseteq
S\cap\cone(w)$ and $g$ is a pnp
diagonalization which has level harmony. Then $f\circ g$
is a pnp diagonalization into $S\cap\cone(w)$
which has level harmony.
\end{lemma}

\begin{proof}
First notice that $f$ and $g$ are in $\mathcal{D}$,
so by Lemma \ref{lem.curly.D}, their composition is
a polite pnp map whose range is a strongly diagonal set.
Since $\qle$ agrees with the lexicographic order on
incomparable pairs, $g$ preserves $\qle$ and sends
all pairs to incomparable pairs, the composition
$f\circ g$ preserves $\qle$.  Since $g$ is injective
with range a strongly diagonal set and $f$ is
a pnp map, the composition is injective.
Thus, since the range of $f\circ g$ is a subset of
the range of $f$, the composition, $f\circ g$, is 
a pnp diagonalization into $S\cap\cone(w)$.

For notational convenience, let $h=f\circ g$. Then 
$
%
\hat{h}(s):=h(s\cat\langle 0\rangle)\meet h(s\cat\langle 1\rangle)
=f(g(s\cat\langle 0\rangle))\meet f(g(s\cat\langle 1\rangle))
$.

\begin{claim}\label{claim.hat.self}
For any $s\in \pre\kappa2$, $\hat{h}(s)=
h(s\cat\langle 1-\delta\rangle)\meet h(s)\subsetneq h(s)$,
where $\delta=g(s)(\hat{g}(s))$.
\end{claim}

\begin{proof}
Since $g$ has level harmony, $\hat{g}(s)\subsetneq g(s)$.
Set $\delta:=g(s)(\lg(\hat{g}(s))$.
Since $s\cat\langle 0\rangle\lex s\cat\langle 1\rangle$
and $g$ preserves lexicographic order, it follows that
$g(s\cat\langle 0\rangle)\lex g(s\cat\langle 1\rangle)$,
so $g(s\cat\langle \delta\rangle)(\lg(\hat{g}(s))=\delta$.

Since $f$ is polite, by meet regularity,
$f(g(s\cat\langle 1-\delta\rangle))
\meet
f(g(s\cat\langle \delta\rangle))
=
f(g(s\cat\langle 1-\delta\rangle))
\meet
f(g(s))$.  That is, $\hat{h}(s)=h(s\cat\langle 1-\delta\rangle)\meet h(s)$.
\end{proof}

\begin{claim}\label{claim.hat.lex}
The function $\hat{h}$ preserves lexicographic order.
\end{claim}

\begin{proof} 
Suppose $s\lex t$.  Since $g$ has level harmony, $\hat{g}(s)
\lex \hat{g}(t)$, $\hat{g}(s)\subsetneq g(s)$  and
$\hat{g}(t)\subsetneq g(t)$.  Hence $\hat{g}(s)\meet\hat{g}(t)=
g(s)\meet g(t)$.  Thus 
$\lg(g(s)\meet g(t))<
\lg(g(s\cat\langle 0\rangle)
\meet 
g(s\cat\langle 0\rangle))$
and 
$\lg(g(s)\meet g(t))<
\lg(g(t\cat\langle 0\rangle)
\meet 
g(t\cat\langle 0\rangle))$.
Since $f$ preserves meet length order, it follows that
$\lg(h(s)\meet h(t))<
\lg(\hat{h}(s))$
and
$\lg(h(s)\meet h(t))<
\lg(\hat{h}(t))$.

Since $f$ preserves lexicographic
order, $h(s)\lex h(t)$.  By the Claim \ref{claim.hat.self},
$\hat{h}(s)\subsetneq h(s)$ and 
$\hat{h}(t)\subsetneq h(t)$.
Thus $\hat{h}(s)\lex\hat{h}(t)$.
\end{proof}

\begin{claim}\label{claim.hat.diff.length} 
For all $s$ and $t$, $\lg(s)<\lg(t)$ implies 
$\lg(h(s))<\lg(\hat{h}(t))$.
\end{claim}

\begin{proof}
Suppose $\lg(s)<\lg(t)$.  Then $\lg(g(s))<\lg(\hat{g}(t))
=lg(g(s\cat\langle 0\rangle)\meet g(s\cat\langle 1\rangle)$.
Since $g(s)\meet g(s)=g(s)$, by preservation of meet length
order by $f$, it follows that $\lg(h(s)<\lg(\hat{h}(t))$.
\end{proof}

\begin{claim}\label{claim.hat.same.length} 
For all $s$ and $t$, $\lg(s)=\lg(t)$ implies 
$\lg(\hat{h}(s))<\lg(h(t))$.
\end{claim}

\begin{proof}
Suppose $\lg(s)=\lg(t)$.  Then $\lg(\hat{g}(s))<\lg(g(t))$.
Argue as in the previous claim: by preservation of meet length
order by $f$, it follows that $\lg(\hat{h}(s)<\lg(h(t))$.
\end{proof}

\begin{claim}\label{claim.hat.extend}
The function $\hat{h}$ preserves extension.
\end{claim}

\begin{proof} 
Suppose $s\subsetneq t$.  Then $\hat{g}(s)\subsetneq \hat{g}(t)$,
since $g$ has level harmony.  Since $f$ satisfies preservation of
meet length order, $\lg(\hat{h}(s))<\lg(\hat{g}(s))$.
By Claim \ref{claim.hat.self}, $\hat{h}(t)\subsetneq h(t)$, so
$\lg(\hat{h}(s))<\lg(h(t))$.  Thus to show $\hat{h}(s)\subseteq
\hat{h}(t)$, it is enough to show $\hat{h}(s)\subseteq h(t)$.

If $h(t)(\lg(\hat{h}(s)))=1$,
then $\hat{h}(s)\subseteq h(t)$, since the range of $h$
is strongly diagonal.  
If $t$ is one of $s\cat\langle 0\rangle$ and $s\cat\langle 1\rangle$,
then  $\hat{h}(s)\subseteq h(t)$ by definition of $\hat{h}(s)$.

So suppose $h(t)(\lg(\hat{h}(s)))=0$ and $\lg(t)>\lg(s)+1$.  
Since $g$ has level harmony, $\hat{g}(s)\subsetneq g(s)$
and $\hat{g}(s)\subseteq\hat{g(t)}\subsetneq g(t)$.
By Claim \ref{claim.hat.self}, 
$\hat{h}(s)=h(t\cat\langle 1-\delta\rangle)
\meet h(s)$ where 
$\delta=g(\lg(\hat{g}(s))$.
For notational convenience, let $\alpha=\lg(\hat{g}(s))$.
For the first subcase, suppose $g(t)(\alpha)\ne g(s)(\alpha)$.
Then $g(t)(\alpha)=g(s\cat\langle 1-\delta\rangle)(\alpha)$,
so $g(s)\meet g(s\cat\langle 1-\delta\rangle)=g(s)\meet g(t)$.
Since the range of $g$ is strongly diagonal and $f$ is polite,
it follows that 
$f(g(s))\meet f(g(s\cat\langle 1-\delta\rangle)=f(g(s))\meet f(g(t))$,
so $\hat{h}(s)=h(s)\meet h(t)\subseteq h(t)$.
For the second subcase, in which $g(t)(\alpha)= g(s)(\alpha)$,
interchange the roles of $g(s)$ and $g(s\cat\langle 1-\delta\rangle)$.
The parallel argument concludes with the inclusion
$\hat{h}(s)=h(s\cat\langle 1-\delta\rangle)\meet h(t)\subseteq h(t)$.
\end{proof}
Now the lemma follows from the claims.
\end{proof}

The next lemma gives an inequality, for pnp maps, which compares
lengths of meets of images with  lengths of images of meets.

\begin{lemma}\label{lem.meet.localization}
Suppose $g:\krtree\to \krtree$ is a pnp map and
$\set{x,u,v}$ is a three element strongly diagonal set with
$x\meet u=x\meet v\subsetneq u\meet v$.
Then $\lg(g(x)\meet g(u))\le\lg(g(x\meet u))$
and $\lg(g(x)\meet g(v))\le\lg(g(x\meet u))$.
\end{lemma}

\begin{proof}
Let $\alpha:=\lg(u\meet x)=\lg(u\meet v)$ and
set $\beta:=\lg(g(x\meet u))=\lg(g(x\meet v)$.
Since $g$ is a pnp map, $\beta<\lg(g(x))$
and $g(x)(\beta)=x(\alpha)$.
Similarly,
$\beta<\lg(g(u))$
and $g(u)(\beta)=u(\alpha)$.
Also, $\beta<\lg(g(v))$
and $g(v)(\beta)=v(\alpha)$.

Since $x\meet u=x\meet v$, 
it follows that $x(\alpha)\ne u(\alpha)=v(\alpha)$.
Consequently, $g(x)(\beta)\ne g(u)(\beta)=g(v)(\beta)$.
Thus $\lg(g(x)\meet g(u))\le\beta$ and
$\lg(g(x)\meet g(v))\le\beta$, so the lemma follows.
\end{proof}

Next we show how to use Shelah's Theorem \ref{thm.SHL}
to obtain a pnp map with meet regularity from 
a pnp map whose range is a strongly diagonal set.

\begin{lemma}\label{lem.meet.reg}
Suppose that $\kappa$ is a cardinal which is measurable
in the generic extension obtained by adding $\lambda$ 
Cohen subsets of $\kappa$, where $\lambda\to(\kappa)^{6}_{2^{\kappa}}$.  
Further suppose $g:\krtree\to\krtree$ is a pnp map whose range is a 
strongly diagonal set. 
Then there is a pnp map $f$ which satisfies meet regularity 
and whose range is a subset of the range of $g$.
\end{lemma}

\begin{proof}
Apply Lemma \ref{lem.poly.diag} to obtain a pre-$S$-vip
order $\prec$ and a sequence $\seq{\varphi_t:t\in\krtree}$
such that for all $t\in\krtree$, the three listed
properties of the lemma hold.  
Say a three element subset $\set{b_0,b_1,b_2}$ listed in increasing
lexicographic order \emph{witnesses meet regularity for $g$}
if it is diagonal and either 
($b_0\meet b_1=b_0\meet b_2$ and 
$g(b_0)\meet g(b_1)=g(b_0)\meet g(b_2)$)
or
($b_2\meet b_0=b_2\meet b_1$ and 
$g(b_2)\meet g(b_0)=g(b_2)\meet g(b_1)$).
Say a three element subset $\set{b_0,b_1,b_2}$ listed in increasing
lexicographic order \emph{refutes meet regularity for $g$}
if it is diagonal, 
$b_0\meet b_1=b_0\meet b_2$ implies
$g(b_0)\meet g(b_1)\ne g(b_0)\meet g(b_2)$
and
$b_2\meet b_0=b_2\meet b_1$ implies
$g(b_2)\meet g(b_0)\ne g(b_2)\meet g(b_1)$.

Define a coloring $d'$ on three element subsets of $\pre\kappa2$
by $d'(x)=0$ if $x$ is strongly diagonal and witnesses meet 
regularity for $g$, $d'(x)=1$ if $x$ is strongly diagonal and
refutes meet regularity for $g$ and $d'(x)=2$ otherwise.
Apply Shelah's Theorem \ref{thm.SHL} to $d'$ and $\prec$
to obtain a strong embedding $e$ and $w\in\pre\kappa2$
such that $e$ preserves $\prec$ on $\cone(w)$ and
for all three element subsets $x$ of $T=e[\cone(w)]$,
the value of $d'(x)$ depends only on the $\prec$-ordered
similarity type of $x$. 

Let $\psi=\varphi_w$. Then $\psi$ is 
a pnp diagonalization into $S\cap\cone(w)$.
We claim that $f:=g\circ e\circ\psi$ is a pnp map
which satisfies meet regularity and 
whose range is a subset of the range of $g$.

Since $e\circ\psi$ is a pnp map which witnesses meet regularity
and whose range is a strongly diagonal set, 
it is enough to show that every strongly diagonal
subset $x\subseteq T$ witnesses meet regularity for $g$. 

Let $(\tau,\dless)$ be an arbitrary vip $m$-type.
Enumerate the leaves of $\tau$ in increasing order 
of length as $a_0, a_1,a_2$.  We must show that 
$(\tau,\dless)$ witnesses meet regularity for $g$.

\begin{claim}\label{claim.x.short} 
If $\lg(a_0)=1$, then every $x\subseteq T$ with
ordered similarity type $(\tau,\dless)$ witnesses meet
regularity for $g$.     
\end{claim}

\begin{proof}
Notice that since that since $\lg(a_0)=1$ and $\tau$ has an
elements of the meet closure of the leaves of lengths
$0,1,2,3,4$, we must have
$\lg(a_0)=1<\lg(a_1\meet a_2)$.

Let $\nu$ be a cardinal larger than $\lg(g(e(\psi(a_0))))$.
For $\alpha<\nu$,
let $b_\alpha$ be the sequence of length $2\alpha+3$ if $\alpha<\omega$
and of length $\gamma+2n+1$ if $\alpha=\gamma+n$ for some limit
ordinal $\gamma\ge\omega$ and $n<\omega$ such that
\[
b_\alpha(\eta)=
\begin{cases}
a_1(0),&\text{if $\eta=0$,}\\
a_1(1),&\text{if $\eta=1$,}\\
a_2(2),&\text{if $0<\eta<\lg(b_\alpha)-1$ even,}\\
a_1(2),&\text{if $\eta=\lg(b_\alpha)-1$,}\\
a_2(3),&\text{otherwise.}\\
\end{cases}
\]
Note that the length of $b_\alpha$ is always odd and
at least $3$, and that for odd $\eta\ge 3$, $b_\alpha(\eta)=a_2(3)$.
For $\beta<\alpha$, by construction, $b_\beta\meet b_\alpha=
b_\beta\restrict(\lg(b_\beta)-1)$.  
Using the above calculations, the careful reader may check
that for $\beta<\alpha$, one has $\clp(\set{a_0,b_\beta,b_\alpha})=\tau$.
Notice that $a_0\meet b_\beta=a_0\meet b_\alpha=\emptyset$.
For $\beta<\nu$, let $s_\beta=e(\psi(b_\beta))$, set $t_0=e(\psi(a_0))$  
and let $x(\beta,\alpha)=\set{t_0,s_\beta,s_\alpha}$.

By Lemma \ref{lem.diag.to.diag}, since $e\circ\psi$ is a polite pnp map,  
we have 
$\clp(x(\beta,\alpha))=\clp(\set{a_0,b_\beta,b_\alpha})=\tau$.
Since there is only one vip order on $\tau$, it follows that
$(\clp(x(\beta,\alpha)),\dless)=(\tau,\dless)$.

Since each $x(\beta,\alpha)$ is a strongly diagonal set,
it either witnesses or refutes meet regularity for $g$.
Since each $x(\beta,\alpha)$ is a subset of $T$, either
they all witness or all refute meet regularity for $g$.

Since $e\circ\psi$ satisfies meet regularity,
$t_0\meet s_\beta=t_0\meet s_\alpha$. 
Suppose each $x(\beta,\alpha)$ refutes meet regularity.
Then $g(t_0)\meet g(s_\beta)\ne g(t_0)\meet
g(s_\alpha)$ for $\beta<\alpha<\nu$. 
Since $\nu$ is a cardinal larger than $\lg(g(e(\psi(a_0))))=\lg(g(t_0))$,
by the Pigeonhole Principle, there are $\beta_0<\alpha_0$
with $g(t_0)\meet g(s_\beta)= g(t_0)\meet
g(s_\alpha)$.  This contradiction shows that 
each $x(\beta,\alpha)$ witnesses meet regularity, so
by choice of $d'$, $e$ and $w$, the claim follows.
\end{proof}

\begin{claim}\label{claim.x.long} 
If $\lg(a_0\meet a_1)=1$, then every $x\subseteq T$ with
ordered similarity type $(\tau,\dless)$ witnesses meet
regularity for $g$.     
\end{claim}

\begin{proof}
Notice that since that since $\lg(a_0\meet a_1)=1$ and $\tau$ has an
elements of the meet closure of the leaves of lengths
$0,1,2,3,4$, we must have
$\lg(a_0\meet a_1)=1<\lg(a_0)<\lg(a_1)<\lg(a_2)$.
As in the previous case, there is only one vip level order on $\tau$.

Let $z=e(\psi(a_0))$ and let 
$\nu=\left(2^{|\lg(g(z))|}\right)^+$.
We proceed much as in the previous case, except we start 
by constructing $\seq{b_\alpha:\alpha<\nu}$ and
$\seq{c_\alpha:\alpha<\nu}$ such that for all $\beta<\alpha$,
$c_\alpha\meet b_\beta=\emptyset=c_\alpha\meet b_\alpha$.
Then we set $t_\alpha=e(\psi(c_\alpha))$, $s_\alpha=e(\psi(b_\alpha))$,
and 
$x(\beta,\alpha)=\set{s_\beta,s_\alpha,t_\alpha}$.
In the previous case we had only $t_0$, so in this case we will
use the Pigeonhole Principle to select a collection of $\alpha$'s
for which the behavior of $g$ on $t_\alpha$ is sufficiently
uniform to reduce this case to one like the previous one.

For $\alpha=\gamma+n$, where $\gamma=0$ or $\gamma$ limit, 
let $b_\alpha$ be the sequence of length $\gamma+2n+2$ such that
\[
b_\alpha(\eta)=
\begin{cases}
a_1(0),&\text{if $\eta=0$,}\\
a_1(1),&\text{if $\eta<\lg(b_\alpha)-1$ odd,}\\
a_1(2),&\text{if $0<\eta<\lg(b_\alpha)-1$ even,}\\
a_0(1),&\text{if $\eta=\lg(b_\alpha)-1$.}\\
\end{cases}
\]
Note that the length of $b_\alpha$ is always even and its Cantor
normal form has an even finite part that is at least $2$.
Also, for positive even $\eta$, $b_\alpha(\eta)=a_1(2)$.
For $\beta<\alpha$, by construction, $b_\beta\meet b_\alpha=
b_\beta\restrict(\lg(b_\beta)-1)$.  
For $\alpha<\eta$, let $c_\alpha$ be the sequence of length
$\lg(b_\alpha)+1$ such that $c_\alpha(0)=a_2(0)$, $c_\alpha(1)=a_2(1)$,
$c_\alpha(\lg(b_\alpha))=a_2(3)$, and 
for all $\eta$ with $1<\eta<\lg(b_\alpha)$, $c_\alpha(\eta)=a_2(2)$. 
Then for all $\alpha,\beta<\nu$, $b_\beta\meet c_\alpha=\emptyset$.
Using these definitions and calculations,
the careful reader may check
that for $\beta<\alpha$, one has $\clp(\set{b_\beta,b_\alpha,c_\alpha})=\tau$.

Since $e\circ\psi$ is a polite pnp map,  
by Lemma \ref{lem.diag.to.diag}, we have
$\clp(x(\beta,\alpha))=\clp(\set{b_\beta,b_\alpha,c_\alpha})=\tau$.
Since there is only one vip order on $\tau$, it follows that
$(\clp(x(\beta,\alpha)),\dless)=(\tau,\dless)$.

Since $e\circ\psi$ satisfies meet regularity, 
$t_\alpha\meet s_\beta=t_\alpha\meet s_\alpha$ for all $\beta<\alpha<\nu$.

By Lemma \ref{lem.meet.localization}, 
$\lg(g(t_\alpha) \meet g(s_\beta)) \le \lg(g(e(\psi(c_\alpha\meet
b_\alpha)))) =\lg(g(z))$.  Similarly, 
$\lg(g(t_\alpha)\meet g(s_\beta))\le \lg(g(z))$.

Set $\gamma=\lg(g(z))$.
Apply the Pigeonhole Principle, to find a sequence $p$ of length
$\gamma$ so that the collection $A=A_p=\set{\alpha:
  g(t_\alpha)\restrict\gamma =p}$ has cardinality $\nu$.

Since each $x(\beta,\alpha)$ is a strongly diagonal set,
it either witnesses or refutes meet regularity for $g$.
Since each $x(\beta,\alpha)$ is a subset of $T$, either
they all witness or all refute meet regularity for $g$.
Suppose each $x(\beta,\alpha)$ refutes meet regularity.
Then $g(t_\alpha)\meet g(s_\beta)\ne g(t_\alpha)\meet
g(s_\alpha)$ for $\beta<\alpha<\nu$. 

In particular, for all $\beta<\alpha$ from $A$, one has
\[
\begin{array}{rl}
g(t_\beta)\meet g(s_\beta)&=p\meet g(s_\beta)\\
&=g(t_\alpha)\meet g(s_\beta)\\
&\ne g(t_\alpha)\meet g(s_\alpha).\\
\end{array}
\]
Now we have a contradiction that $\set{g(t_\alpha)\meet g(s_\alpha):
\alpha\in A}$ is a set of $\nu$ distinct initial segments of $p$
and $\lg(p)=\gamma<\nu$.  This contradiction shows that 
each $x(\beta,\alpha)$ witnesses meet regularity, so
by choice of $d'$, $e$ and $w$, the claim follows.
\end{proof}

\begin{claim}\label{claim.x.med} 
If $\lg(a_0\meet a_2)=1$, then every $x\subseteq T$ with
ordered similarity type $(\tau,\dless)$ witnesses meet
regularity for $g$.     
\end{claim}

\begin{proof}
Notice that since that since $\lg(a_0\meet a_2)=1$ and $\tau$ has an
elements of the meet closure of the leaves of lengths
$0,1,2,3,4$, we must have
$\lg(a_0\meet a_2)=1<\lg(a_0)<\lg(a_1)<\lg(a_2)$.
As in the previous case, there is only one vip level order on $\tau$.

Let $z=e(\psi(a_0))$ and let 
$\nu=\left(2^{|\lg(g(z))|}\right)^+$.
We proceed much as in the previous case, except we start 
by constructing $\seq{b_\beta:\beta<\nu}$ and
$\seq{c_\beta:\alpha<\nu}$ such that for all $\beta<\alpha$,
$c_\beta\meet b_\beta=\emptyset=c_\beta\meet b_\alpha$.
Then we set $t_\beta=e(\psi(c_\beta))$, $s_\beta=e(\psi(b_\beta))$,
and 
$x(\beta,\alpha)=\set{s_\beta,t_\beta,s_\alpha}$.

For $\beta=\gamma+n$, where $\gamma=0$ or $\gamma$ limit, 
let $b_\beta$ be the sequence of length $\gamma+4n+2$ such that
for $\eta=\zeta+4k+\ell$ where $\zeta=0$ or $\zeta$ limit, 
$k<\omega$ and $\ell<4$,
\[
b_\beta(\eta)=
\begin{cases}
a_2(\ell),&\text{if $\eta<\lg(b_\beta)-1$,}\\
a_0(1),&\text{if $\eta=\lg(b_\beta)-1$.}\\
\end{cases}
\]
Note that the length of $b_\beta$ is always even and its Cantor
normal form has an even finite part that is at least $2$.
For $\beta<\alpha$, by construction, $b_\beta\meet b_\alpha=
b_\beta\restrict(\lg(b_\beta)-1)$, and $b_\alpha(\lg(b_\beta)-1)=a_2(1)$.
For $\beta<\nu$, let $c_\beta$ be the sequence of length
$\lg(b_\beta)+1$ such that $c_\beta(0)=a_1(0)$, 
$c_\beta(\lg(b_\beta)-1)=a_1(1)$,
$c_\beta(\lg(b_\beta))=a_1(2)$, and 
for all $\eta$ with $0<\eta<\lg(b_\beta)-1$, $c_\beta(\eta)=0$. 
Then for all $\beta\le\alpha<\nu$, $b_\beta\meet c_\alpha=\emptyset$.
Using these definitions and calculations,
the careful reader may check
that for $\beta<\alpha$, one has $\clp(\set{b_\beta,c_\beta,b_\alpha})=\tau$.

The rest of the proof in this case parallels that of the
previous one and the details are left to the reader.
\end{proof}

Since $S$ is cofinal and transverse, 
its intersection with $\cone(p_1)$ is cofinal above $p_1$ and transverse.  
Thus by Lemma \ref{lem.cofinal.dense}, $(S\cap\cone(w),\qle)$
is $\kappa$-dense.  Then $U:=T(S\cap\cone(w))$ is an almost
perfect tree by Lemma \ref{lem.nearly.perfect}.  

\begin{claim}\label{claim.x.from.U}  
If $\lg(a_1\meet a_2)=1$ and $p_0,p_1$ are incomparable elements
of $U$ such that $p_0\meet p_1$ a densely splitting point of
$S\cap\cone(w)$ and $p_0\lex p_1$ if and only if $a_0\lex a_1$, 
then there are $c\in U\cap\cone(p_0)$ and $b,b'\in U\cap\cone(p_1)$
with $\max(\lg(p_0),\lg(p_1))<\lg(b\meet b')<\lg(b)<\lg(b')$ such that  
$(\set{c,b,b'},\prec)$ has the same ordered similarity type as
$(\tau,\dless)$.
\end{claim}

\begin{proof}
By Lemma \ref{lem.1.almost.perfect}, $U$ is almost perfect.
Use that fact to find a densely splitting node $q_1$ of  $U$
extending $p_1$ of length longer than $\max(\lg(p_0),\lg(p_1))$.

Let $C(S\cap\cone(w))$ be the set of all limit ordinals $\alpha>0$ such that
every $t\in U\cap {}^{\alpha>}2$ has proper extensions
in both $S\cap {}^{\alpha>}2$ and $\mathcal{W}(S)\cap {}^{\alpha>}2$
By Lemma\ref{lem.nearly.perfect}, $C(S\cap\cone(w))$
is closed and unbounded in $\kappa$.

Let $\lambda\in C(S\cap\cone(w))$ 
be a regular cardinal greater than $\lg(q_1)$. 
Use the fact that $U$ is almost perfect to find
an extension $c$ of $p_0$ in $S$ of length greater than $\lambda$.

Define $h:\prele\lambda2\to U\cap\cone(q_1)$ by recursion as follows
and prove by induction that for all $s\in\pre\lambda2$, $h(s)$
is a densely splitting point of $U\cap\cone(q_1)\cap\pre\lambda2$.
To start the recursion, let 
$h(\emptyset)=q_1$.  If $h(s)$ has been defined and $\delta<2$,
let $h_0(s\cat\langle\delta\rangle)$ be an extension of
$h(s)\cat\langle\delta\rangle$ of length less than $\lambda$
which is a densely splitting node of $S\cap\cone(w)$ and
let  $h(s\cat\langle\delta\rangle)$ be an extension of
$h_0(s)\cat\langle 1-\delta\rangle$ of length less than $\lambda$
which is a densely splitting node of $S\cap\cone(w)$.
If $\gamma<\lambda$ is a limit ordinal and $h(s\restrict \beta)$
has been defined for all $\beta<\gamma$, let $h(s)$ be
an extension of
$\bigcup\set{h(s)\restrict\beta:\beta<\gamma}$ 
of length less than $\lambda$
which is a densely splitting node of $S\cap\cone(w)$.
If $h(s\restrict \beta)$
has been defined for all $\beta<\lambda$, let $h(s)$ be
$\bigcup\set{h(s)\restrict\beta:\beta<\lambda}$. 
Since $U$ is an almost perfect tree and for each $s$ of
limit length $\gamma$, the node
$\bigcup\set{h(s)\restrict\beta:\beta<\gamma}$ is
an even-handed limit of densely splitting points,
the recursion is well-defined.

Set $A:=h[\pre\lambda2]$.  Since $\lambda$ is regular,
it follows that $A\subseteq\pre\lambda2$.  Also, $|A|=2^\lambda$.
 
Let $\xi=\lg(c)$ and let $h^*:A\to U\cap\pre\xi2$ be an injection
with $s\subseteq h^*(s)$ for all $s\in A$.  Such an injection 
exists since $U$ is almost perfect.  Define $\prec^*$ on $A$
by $s\prec^* t$ if and only if $h^*(s)\prec h^*(t)$.
Notice that $\prec^*$ is a well-order, since $\prec$ is a well order.  

Since $h$ preserves lexicographic order, we can find 
$r_0\in A$ such that both $\set{s\in A:s\lex r_0}$
and $\set{t\in A:r_0\lex t}$ have cardinality $2^\lambda$. 
Thus there are $s_0$ and $t_0$ with $s_0\lex r_0\lex t_0$
and $r_0\prec^* s_0$, $r_0\prec^* t_0$.  That is there
are pairs from $A$ where the lexicographic order and
the $\prec^*$-order agree and where they disagree.

Choose $r_1,r_2\in A$ such that $r_1\lex r_2$ if and only
if $a_1\lex a_2$ and $r_1\prec^* r_2$ if and only if
$a_1\restrict 2\dless a_2\restrict 2$.

Use the fact that $U$ is almost perfect to find be an element
$b$  of $S$ extending $h^*(a_1)$ of
length greater than $\xi=\lg(c)$ and to find an element $b'$
of $S$ extending $h^*(a_2)$ of length greater than $\lg(b)$.
Then $(\set{c,b,b'},\prec)$ has the same ordered similarity
type as $(\tau,\dless)$.
\end{proof}

\begin{claim}\label{claim.x.special} 
If $\lg(a_1\meet a_2)=1$, then every $x\subseteq T$ with
ordered similarity type $(\tau,\dless)$ witnesses meet
regularity for $g$.     
\end{claim}

\begin{proof}
Notice that since that since $\lg(a_1\meet a_2)=1$ and $\tau$ has an
elements of the meet closure of the leaves of lengths
$0,1,2,3,4$, we must have
$\lg(a_1\meet a_2)=1<\lg(a_0)<\lg(a_1)<\lg(a_2)$.

The proof of this claim shares similarities with the proofs
of Claims \ref{claim.x.short}, \ref{claim.x.med} and
\ref{claim.x.long}, 
once the construction of the family
of sets $x(\beta,\alpha)$ has been completed.

Let $p_{0}$ and $p_{1}$ be incomparable elements of $U$
such that $p_{0}\meet p_{1}$ 
is a densely splitting node of $S\cap\cone(w)$
and $p_{0}\lex p_{1}$ if and only if $a_0\lex a_1$.

Set $z=e(p_{0}\meet p_{1})$ and let 
$\nu=\left(2^{|\lg(g(z))|}\right)^+$.

Define by recursion sequences 
$\seq{p_{\delta,\alpha}:\alpha<\nu}$ for $\delta<2$,
$\seq{c_\alpha:\alpha<\nu}$
$\seq{b_\alpha:\alpha<\nu}$ and
$\seq{b'_\alpha:\alpha<\nu}$ as follows.
To start the recursion,
let  $p_{0,0}=p_0$ and $p_{1,0}=p_1$.
If $p_{0,\alpha}$ and $p_{1,\alpha}$ have been defined with
$p_{0}\subseteq p_{0,\alpha}$ and
$p_{1}\subseteq p_{1,\alpha}$, then 
apply Claim \ref{claim.x.from.U}
to $p_{0,\alpha}$ and $p_{1,\alpha}$ to
obtain $c_\alpha\in S\cap\cone(p_{0,\alpha})$
and $b_\alpha,b'_\alpha\in S\cap\cone(p_{1,\alpha})$ 
such that $(\set{c_\alpha,b_\alpha,b'_\alpha},\prec)$
has the same ordered similarity type as $(\tau,\dless)$
and the following inequalities hold:
\[
\max(\lg(p_{0,\alpha}),\lg(p_{1,\alpha}))<
\lg(b_\alpha\meet b'_\alpha)<\lg(b_\alpha)<\lg(b'_\alpha).
\]
Note that $c_\alpha\meet b_\alpha=p_{0}\meet p_{1}$.

If $\alpha>0$ and
$p_{0,\beta}$, $p_{1,\beta}$, 
$c_\beta$,
$b_\beta$ and
$b'_\beta$
have all been defined for $\beta<\alpha$,  
then set 
$p_{0,\alpha}:=\bigcup\set{c_\beta:\beta<\alpha}$ and 
$p_{1,\alpha}:=\bigcup\set{b'_\beta:\beta<\alpha}$.

By construction, the sequence $\mathcal{S}=\seq{b_\alpha:\alpha<\nu}$
is increasing in length.  Fix attention on a specific
pair $\beta<\alpha<\nu$.
Since $s_\alpha\in \cone(p_{1,\alpha})$, it follows from the
definition of $p_{1,\alpha}$ that $b'_\beta\subseteq s_\alpha$.
Thus $(\set{c_\beta,b_\beta,b'_beta},\prec)$
and $(\set{c_\beta,b_\beta,b_\alpha},\prec)$ have the
same $\prec$-ordered similarity type, namely $(\tau,\dless)$.
Furthermore, $c_\beta\meet b_\beta=p_0\meet p_1$ and
$c_\beta\meet b_\alpha=p_0\meet p_1$.

For $\alpha<\nu$, let $t_\alpha=e(c_\alpha)$ and let
$s_\alpha=e(b_\alpha)$. For $\beta<\alpha<\nu$,
let $x(\beta,\alpha)=\set{t_\beta,s_\beta,s_\alpha}$.
Since $e$ preserves meets, it follows that
$e(c_\beta)\meet e(b_\beta)=e(p_0\meet p_1)=z$ and
$e(c_\beta)\meet e(b_\alpha)=e(p_0\meet p_1)=z$.
Since $e$ is a  strong embedding which preserves
$\prec$ on $\cone(w)$, and $x(\beta,\alpha)\subseteq\cone(w)$,
it follows that $(x(\beta,\alpha),\prec)$ and 
$(\set{c_\beta,b_\beta,b_\alpha},\prec)$ have the
same $\prec$-ordered similarity type, namely $(\tau,\dless)$. 

The rest of the proof in this case parallels those of 
Claims \ref{claim.x.med} and \ref{claim.x.long}
and the details are left to the reader.
\end{proof}

By Claims \ref{claim.x.short}, \ref{claim.x.med}, 
\ref{claim.x.long} and \ref{claim.x.special}, we have
shown that $(\tau,\dless)$ witnesses meet regularity for $g$.
Since $(\tau,\dless)$ was an arbitrary vip $m$-type,
every such type witness meet regularity for $g$.
Thus $f=g\circ e\circ\psi$ is a pnp map with meet regularity whose
range is a subset of the range of $g$.
\end{proof}

\begin{lemma}\label{lem.lex.preserve}
Suppose that $\kappa$ is a cardinal which is measurable
in the generic extension obtained by adding $\lambda$ 
Cohen subsets of $\kappa$, where $\lambda\to(\kappa)^{6}_{2^{\kappa}}$.  
Further suppose $g:\krtree\to\krtree$ is a pnp map which satisfies
meet regularity and whose range is a strongly diagonal set. 
Then there is a pnp map $f$ which satisfies meet regularity and 
preserves lexicographic order
and whose range is a subset of the range of $g$.
\end{lemma}

\begin{proof}
Apply Lemma \ref{lem.poly.diag} to obtain a pre-$S$-vip
order $\prec$ and a sequence $\seq{\varphi_t:t\in\krtree}$
such that for all $t\in\krtree$, the three listed
properties of the lemma hold.  
Let $\tau_0$ be the $2$-type whose leaves are $\langle 0\rangle$
and $\langle 1,0\rangle$.
Let $\tau_0$ be the $2$-type whose leaves are $\langle 1\rangle$
and $\langle 0,0\rangle$.
Let $\tau_2$ be the $2$-type whose leaves are $\langle 0\rangle$
and $\langle 1,1\rangle$.
Let $\tau_3$ be the $2$-type whose leaves are $\langle 1\rangle$
and $\langle 0,1\rangle$.
Then for all $x\in[\pre\kappa2]^2$, the set $\clp(x)$ is one of
$\tau_0$, $\tau_1$, $\tau_2$ and $\tau_3$.

Define a coloring $\ddprime$ on pairs from $\pre\kappa2$ by
$\ddprime(x)=j$ where $\clp(g[x])=\tau_j$.

Apply Shelah's Theorem \ref{thm.SHL} to $\ddprime$ and $\prec$
to obtain a strong embedding $e$ and $w\in\pre\kappa2$
such that $e$ preserves $\prec$ on $\cone(w)$ and
for all pairs $x$ of $T=e[\cone(w)]$,
the value of $\ddprime(x)$ depends only on the $\prec$-ordered
similarity type of $x$. 

Let $\psi=\varphi_w$. Then $\psi$ is 
a pnp diagonalization into $S\cap\cone(w)$.
We claim that $f:=g\circ e\circ\psi\circ g\circ e\circ\psi$ is a pnp map
which satisfies meet regularity and preserves lexicographic order.
Clearly the  range of $f$ is a subset of the range of $g$.

Since $e\circ\psi$ is a pnp map which witnesses meet regularity
and whose range is a strongly diagonal set, 
by Lemma \ref{lem.meet.reg}, the functions $g\circ e\circ\psi$
and $f$ satisfy meet regularity.

Since $g$ is a pnp map, for any pair $x$ from $e[\cone(w)]$
with $\clp(x)=\tau_0$ or $\clp(x)=\tau_1$, we must have
$\clp(g(x))=\tau_0$ or $\clp(g(x))=\tau_1$.
Similarly, for any pair $x$ from $e[cone(w)]$
with $\clp(x)=\tau_1$ or $\clp(x)=\tau_2$, we must have
$\clp(g(x))=\tau_1$ or $\clp(g(x))=\tau_2$.

Say $g$ \emph{sends $i$ to $j$}
if for all incomparable pairs $x$ from $e[\cone(w)]$ with $\clp(x)=\tau_i$
one has $\clp(g[x])=\tau_j$.

\begin{claim}
Either (a) $g$ sends $0$ to $1$ and $1$ to $0$, or 
(b) $g$ sends $0$ to $0$ and $1$ to $1$.
\end{claim}

\begin{proof}
Consider $a=\langle 0,0,0,0\rangle$, $b=\langle 0,1\rangle$
and $c=\langle 1,0,0\rangle$.  Then $\set{a,b,c}$ is diagonal,
and $c\meet a=\emptyset=c\meet b$.

Now $\clp(\set{a,c})=\tau_1$ and $\clp(\set{b,c})=\tau_0$.
Let $s=e(\psi(a))$, $t=e(\psi(b))$ and $u=e(\psi(c))$.
Set $x:=\set{s,u}$ and $y=\set{t,u}$.
Since $e\circ\psi$ is a polite pnp map, 
by Lemma \ref{lem.diag.to.diag}, we have 
$\clp(x)=\tau_1$ and $\clp(y)=\tau_0$.

Since $g\circ e\circ\psi$ satisfies meet regularity,
$g(u)\meet g(s)=g(u)\meet g(t)$.  That is, $g(s)$
and $g(t)$ are on the same side of $g(u)$.
Let $\gamma=\lg(g(u)\meet g(s))$.
Then $g(s)(\gamma)=g(t)(\gamma)=1- g(u)(\gamma)$.
If $g(s)(\gamma)=0$ then $\clp(g[x])=\tau_1$ and $\clp(g[y])=\tau_0$,
so $g$ sends $1$ to $1$ and $0$ to $0$.
Otherwise $g(s)(\gamma)=1$. In this case $\clp(g[x])=\tau_0$ and 
$\clp(g[y])=\tau_1$,
so $g$ sends $1$ to $0$ and $0$ to $1$.
\end{proof}

\begin{claim} Either
(a) $g$ sends $2$ to $3$ and $3$ to $2$, or 
(b) $g$ sends $2$ to $2$ and $3$ to $3$.
\end{claim}

\begin{proof}
The proof parallels that of the previous claim using
$a=\langle 1,1,1,1\rangle$, $b=\langle 1,0\rangle$
and $c=\langle 0,0,1\rangle$.  The details are left
to the reader.
\end{proof}

Notice that  $e\circ\psi$ sends pairs from $\pre\kappa2$
to incomparable pairs from $\cone(w)$.  For incomparable
pairs, $e\circ\psi$ preserves the similarity type by
Lemma \ref{lem.diag.to.diag}.
Hence by the above claims, 
for any incomparable pair $z$ from $\pre\kappa2$
with $\clp(z)=\tau_i$, we have $\clp(f[z])=\tau_i$.
Therefore $f$ is a pnp map which satisfies meet
regularity and preserves lexicographic order
and has range a subset of the range of $g$. 
\end{proof}

\begin{theorem}\label{thm.thin.inside}
Suppose that $\kappa$ is a cardinal which is measurable
in the generic
extension obtained by adding $\lambda$ Cohen subsets of $\kappa$,
where $\lambda\to(\kappa)^{6}_{2^{\kappa}}$.  
If $g:\krtree\to\krtree$ is a pnp map whose range is a strongly diagonal set
whose meet closure is a subset of $S=\sigma[\kappa]$, 
then there is a pnp diagonalization $f:\krtree\to \ran(g)$ into $S$
with level harmony.
\end{theorem}

\begin{proof}
Apply Lemma \ref{lem.poly.diag} to obtain a pre-$S$-vip
order $\prec$ and a sequence $\seq{\varphi_t:t\in\krtree}$
such that for all $t\in\krtree$, the three listed
properties of the lemma hold.  
Apply Lemma \ref{lem.meet.reg} to $g$ to obtain a pnp map $g_0$
which satisfies meet regularity and whose range is strongly diagonal 
subset of the range of $g$.  Apply Lemma \ref{lem.lex.preserve} to $g_0$
to obtain a pnp map $h$ which satisfies meet regularity and
preserves lexicographic order and whose range is a strongly diagonal
subset of $\ran(g_0)\subseteq \ran(g)$.
 
Define a coloring $d^*$ on three element antichains $b=\set{b_0,b_1,b_2}$
listed in increasing order of length as follows: $d^*(b)=0$ 
if $\lg(b_0)<\lg(b_1\meet b_2)$ and 
$\lg(h(b_0))<\lg(h(b_1)\meet h(b_2))$; 
$d^*(b)=1$ 
if $\lg(b_0)<\lg(b_1\meet b_2)$ and 
$\lg(h(b_0))\ge\lg(h(b_1)\meet h(b_2))$; and
$d^*(b)=2$ otherwise.

Apply Shelah's Theorem \ref{thm.SHL} to $d^*$ and $\prec$
to obtain a strong embedding $e$ and $w\in\pre\kappa2$
such that $e$ preserves $\prec$ on $\cone(w)$ and
for all pairs $x$ of $T=e[\cone(w)]$,
the value of $d^*(x)$ depends only on the $\prec$-ordered
similarity type of $x$. 

Let $\psi=\varphi_w$. Then $\psi$ is 
a pnp diagonalization into $S\cap\cone(w)$ with level harmony.
We claim that $f':=h\circ e\circ\psi$ is a pnp map
which satisfies meet regularity and preserves lexicographic order.
Clearly the  range of $f'$ is a subset of the range of $g$.

Since $e\circ\psi$ is a pnp map which witnesses meet regularity
and preserves lexicographic order, and whose range is a strongly 
diagonal set, 
by Lemma \ref{lem.lex.preserve}, the function $f'$
satisfies meet regularity and preserves lexicographic order.

Our next goal is to prove $f'$ also preserves meet length order.
The claim below is a preliminary step toward that goal.

\begin{claim}\label{claim.low.meet}
If $b=\set{b_0,b_1,b_2}\subseteq e[\cone(w)]$ is a strongly diagonal
set listed in increasing order of length and 
$\lg(b_0)<\lg(b_1\meet b_2)$, then $d^*(b)=0$ and
$\lg(h(b_0))<\lg(h(b_1)\meet \lg(h(b_2))$. 
\end{claim}

\begin{proof}
Assume toward a contradiction that 
$b=\set{b_0,b_1,b_2}\subseteq e[\cone(w)]$ is a strongly diagonal
set listed in increasing order of length with
$\lg(b_0)<\lg(b_1\meet b_2)$ and $d^*(b)\ne 0$. 
Then $d^*(b)=1$.  Let $\tau=\clp(b)$ and 
list the elements of $a$ in increasing order of length
as $\tau=\set{a_0,a_1,a_2}$.  Then $\lg(a_0)$.

Let $z=e(\psi(a_0))$ and let $\nu$ be an uncountable
cardinal greater than $2^{|\lg(h(z))|}$.
Using the proof of Claim \ref{claim.x.short} as a guide,
construct a sequence $\seq{b_\alpha:\alpha<\nu}$
such that for all $\alpha<\nu$, $a_0\meet b_\alpha=\emptyset$
and for all $\beta<\alpha<\nu$,
$\clp(\set{a_0,b_\beta,b_\alpha})=\tau$
and $\lg(b_\beta)<\lg(b_\alpha)$.
Set 
$s_\alpha:=e(\psi(b_\alpha))$ for all $\alpha<\nu$.
For $\beta<\alpha<\nu$, let $x(\beta,\alpha):=\set{z,s_\beta,s_\alpha}$.
Thus for all $\beta<\alpha<\nu$, $\clp(x(\beta,\alpha)=\tau$.

There is only one vip order $\dless$ on $\tau$.
Since $\psi$ is a pnp diagonalization into $S\cap\cone(w)$,  
$\prec$ is a pre-$S$-vip order, and $e$ preserves $\prec$
on $\cone(w)$, it follows that 
$(\clp(x(\beta,\alpha),\prec_{x(\beta,\alpha)})=(\tau,\dless)$,
since both are vip $3$-types and $\clp(x(\beta,\alpha)=\tau$.

Since $d^*$ takes the same value on all triples from $e[\cone(w)]$
of the same $\prec$-ordered similarity type, it follows that
$\lg(h(z))\ge \lg(h(s_\beta)\meet h(s_\alpha))$ for all 
$\beta<\alpha<\nu$.  Since the range of $h$ 
is a strongly diagonal set, the inequalities must be strict.

Since $\nu$ is an uncountable cardinal greater than $2^{|\lg(h(z))|}$,
by the Pigeonhole Principle, there are $\beta<\gamma$ with
$h(t_\beta)\restrict\lg(h(z))=h(t_\gamma)\restrict\lg(h(z))$.
Thus we have reached the contradiction that
$h(t_\beta)\meet h(t_\gamma)$ must have length short and
greater than or equal to $\lg(h(z))$.  Thus the claim follows.
\end{proof}

\begin{claim}\label{claim.trio.meet}
If $\lg(x\meet y)<\lg(u\meet v)$ and $|\set{x,y,u,v}|=3$,
then $\lg(f'(x)\meet f'(y))<\lg(f'(u)\meet f'(v))$.
\end{claim}

\begin{proof}
If $x=y$, then the claim follows from Claim \ref{claim.low.meet}.
So assume $x\ne y$.  Then
one of $x$ and $y$ is either $u$ or $v$, so $x\meet y\subsetneq
u\meet v$.  Let $z$ be the unique one of $x$ and $y$ which
is not in $\set{u,v}$.  Then $z\meet u=z\meet v=x\meet y$.  By meet
regularity, $f'(z)\meet f'(u)=f'(z)\meet f'(v)=f'(x)\meet f'(y)$.  
Let $\beta=\lg(f'(x)\meet f'(y)$.
Since $f'$ preserves lexicographic order, 
$f'(z)(\beta)\ne f'(u)(\beta)=f'(v)(\beta)$.  It follows that
$\lg(f'(u)\meet f'(v))>\beta=\lg(f'(x)\meet f'(y))$.
\end{proof}

\begin{claim}\label{claim.pine}
If $x\meet y\subsetneq u\meet v$ and $|\set{x,y,u,v}|=4$,
then $\lg(f'(x)\meet f'(y))<\lg(f'(u)\meet f'(v))$.
\end{claim}

\begin{proof}
If $x\meet y=x$ or $x\meet y=y$, then the claim follows
from Claim \ref{claim.low.meet}.  So assume $x\ne x\meet y\ne y$.
Since $x\meet y$ is a proper initial segment of $u\meet v$,
either 
$(x\meet y)\cat\langle 0\rangle\subseteq u\meet v$
or
$(x\meet y)\cat\langle 1\rangle\subseteq u\meet v$.
Consequently, either
$x\meet y=x\meet u=x\meet v$ or 
$x\meet y=y\meet u=y\meet v$.
In the first case, the claim follows from Claim \ref{claim.trio.meet}
applied to $\set{x,u,v}$ and in the second case, it follows 
from Claim \ref{claim.trio.meet}
applied to $\set{y,u,v}$.
\end{proof}

\begin{claim}\label{claim.scales}
If $\lg(x\meet y)<\lg(u\meet v)$, $x\meet v\not\subseteq u\meet v$ 
and $|\set{x,y,u,v}|=4$, then $\lg(f'(x)\meet f'(y))<\lg(f'(u)\meet f'(v))$.
\end{claim}

\begin{proof}
Since $\lg(x\meet y)<\lg(u\meet v)$ and $x\meet v\not\subseteq u\meet
v$, it follows that $x\meet y$ and $u\meet v$ are incomparable.
By  Lemma \ref{lem.meet.localization}, $\lg(f'(x)\meet f'(y))
\le lg(f'(x\meet y))$.  

If one of $u$ and $v$ is an initial segment of the other.
Since $e\circ\psi$ is a pnp map, $e(\psi(x\meet y))$ is shorter
than both $e(\psi(u))$ and $e(\psi(v))$.
If neither is an initial segment of the other, the
$e(\psi(x\meet y))$ is shorter than
than both $e(\psi(u))$ and $e(\psi(v))$ since $e\circ\psi$
preserves meet length order.
Thus by Claim \ref{claim.low.meet} applied to
$\set{e(\psi(x\meet y)),e(\psi(u)),e(\psi(v))}$,
$\lg(f'(x\meet y))<\lg f'(u)\meet \lg(f'(v))$ and the claim follows.
\end{proof}

\begin{claim}\label{claim.meet.length}
The function $f'$ preserves meet length order.
\end{claim}

\begin{proof}
Let $x,y,u,v$ be arbitrary with $\lg(x\meet y)<\lg(u\meet v)$.
Consider the set $\set{x,y,u,v}$.  It must have at least
two elements and at most four.  If it has three or four
elements, then $\lg(f'(x)\meet f'(y))<\lg(f'(u)\meet f'(v))$
by one of Claims \ref{claim.trio.meet}, \ref{claim.pine} and
\ref{claim.scales}.  If it has two elements,
say $z$ and $w$ with $\lg(z)<\lg(w)$, then there are
only three possible meets, listed here in increasing
order of length: $z\meet w$, $z\meet z$ and $w\meet w$.
Since $f'$ is a pnp map, $\lg(f'(z))<\lg(f'(w))$.  Since the range
of $f'$ is strongly diagonal set, $f'(z)$ and $f'(w)$ are incomparable,
so $f'(z)\meet f'(w)$ is a proper initial segment of both
$f'(z)$ and $f'(w)$.  Thus $f'$ preserves meet length order.
\end{proof}

Since $f'$ satisfies meet regularity and preserves both
lexicographic order and meet length order, it is polite.
We saw above that $f'$ is a pnp map whose range is a
subset of the range of $g$.  Recall that $\psi$ is a
pnp diagonalization with level harmony into $S\cap\cone(w)$.
Thus by Lemma \ref{lem.harmony}, the function $f=g\circ f'$
is a pnp diagonalization into $S$ with $\ran(f)\subseteq
\ran(g)\subseteq S$. 
\end{proof}

Note that for $m=2$, there are four $m$-types and each of
them admits a single vip order.  Any copy of an
uncountable Rado $\mathbb{G}$ has an induced subgraph
which is a countable Rado graph.  Since Laflamme, Sauer
and Vuksanovic \cite{LSVpreprint} have shown that
these four types must appear in translations of 
every induced subgraph of the countable Rado graph which
is itself isomorphic to the countable Rado graph,
it follows from their work that
$\mathbb{G}_\kappa\nrightarrow(\mathbb{G}_\kappa)^2_{<\omega,r^+_2-1}$.
For larger values of $m$, we use Shelah's Theorem.

\begin{question} Suppose $\kappa$ is an uncountable cardinal
with $\kappa^{<\kappa}=\kappa$, $\mathbb{G}_\kappa$ is a
$\kappa$-Rado graph and $2<m<\omega$.
Does $\mathbb{G}\nrightarrow(\mathbb{G})^m_{<\omega,r^+_m-1}$?
That is, does the lower bound hold even when $\kappa$ does not
satisfy the hypothesis of Shelah's Theorem?   
\end{question}

\begin{theorem}\label{thm.lower.bound.2}
Let $m\ge 3$ be a natural number and 
  suppose that $\kappa$ is a cardinal which is measurable
  in the generic
  extension obtained by adding $\lambda$ Cohen subsets of $\kappa$,
  where $\lambda\to(\kappa)^{6}_{2^{\kappa}}$.
Further suppose $\mathbb{G}_\kappa$ is a $\kappa$-Rado graph.
Then for $r^+_m$ equal to the number of vip $m$-types,
$\mathbb{G}_\kappa$ 
satisfies 
$$\mathbb{G}_\kappa\nrightarrow(\mathbb{G}_\kappa)^m_{<\omega,r^+_m-1}.$$
\end{theorem}

\begin{proof}
Apply Lemma \ref{lem.poly.diag} to obtain a pre-$S$-vip
order $\prec$ and a sequence $\seq{\varphi_t:t\in\krtree}$
such that for all $t\in\krtree$, the three listed
properties of the lemma hold.  
Enumerate the vip $m$-types:
$(\tau_0,\dless_0)$, $(\tau_1,\dless_1)$, \dots , $(\tau_{r-1},\dless_{r-1})$.

Define a coloring  $d:[\krtree]^m\to r$
by $d(x)=i$ if $(\clp(\varphi_0[x]),\prec_{\varphi_0[x]})=(\tau_i,\dless_i)$.
Then $c:[\kappa]^m\to r$ defined by $c(b)=d(\sigma[b])$ is a coloring
of the $m$-element subsets of $\kappa$ with $r$ colors.  

To prove the theorem, we must show that every isomorphic copy
of $\mathbb{G}_\kappa$ contains sets of every color.  Toward that end,
let $K\subseteq \kappa$  be an arbitrary set such
that $(K,E\restrict K)\cong (\kappa,E)$ and let $\rho:\kappa\to K$
be the isomorphism and let $j<r$ be arbitrary.  
By Lemmas \ref{lem.Rado.increasing} and
\ref{thm.1.translate}, we may assume that $\sigma\circ\rho\circ\sigma^{-1}$
is a pnp map with domain $S$.  Hence 
$g:=\varphi_0\circ\sigma\circ\rho\circ\sigma^{-1}\circ\varphi_0$ 
is a pnp map with
domain $\krtree$ whose range is a strongly diagonal subset of $S$ 
and $\closure{\ran(g)} \subseteq S$.

Apply Theorem \ref{thm.thin.inside} to $g$ obtain a pnp
diagonalization $f:\krtree\to\ran(g)$ into $S$ with level harmony.
Since $\prec$ is a pre-$S$-vip order, it is also a $D$-vip order
for $D=\ran(f)$.  

By Theorem \ref{thm:Gk} applied to $f$ and $D$, 
there is some $y\subseteq D$ such that $(\tau_j,\dless_j)=
(\clp(y),\prec^y)$.
Recall that $\ran(f)\subseteq \ran(g)$.  
Let $a$ be such that $g[a]=x$.  Then
$b:=\rho[\sigma^{-1}[\varphi_0[a]]]$ is well-defined
and by definition of $\sigma$ and $\rho$, we can see that $b\subseteq
K$.  Let $x=\sigma[b]$.  Then $\varphi_0[x]=g[a]=y$ by definition of $g$
and choice of $a$.  So $d(x)=j$, by definition of $d$.  Thus
$c(b)=d(\sigma[b])=j$.  Since $K$ and $j$ were arbitrary,
every isomorphic copy of $\mathbb{G}_\kappa$ contains sets of every color.
\end{proof}

Figure \ref{fig.Rplusvalue} summarizes 
the  calculation from \cite{Jeanfincomb} of values of
$r_m^+$ for $m\le 5$.  A comparison with $r_m$, the number
of $m$-types, is also included, where $r_m$ is the critical
value for finite colorings of $m$-tuples of 
the countable Rado graph.

\begin{figure}[t]
\[
\begin{array}{|r|r|r|}
\hline
m&r_m^+&r_m\\
\hline
1&1&1\\
\hline
2&4&4\\
\hline
3&128&112\\
\hline
4&26,368&12,352\\
\hline
5&41,932,288&4,437,760\\
\hline
\end{array}
\]
\caption{Some small values of $r_m^+$ and $r_m$.}\label{fig.Rplusvalue}
\end{figure}
\bibliographystyle{plain}%

\end{document}